\input amstex
\documentstyle{amsppt}
\magnification=\magstep1 \NoRunningHeads

\topmatter

\title
 Odometer actions of the Heisenberg group
\endtitle

\author
Alexandre I. Danilenko  and Mariusz Lemanczyk
\endauthor
\email
alexandre.danilenko@gmail.com
\endemail

\address
 Institute for Low Temperature Physics
\& Engineering of National Academy of Sciences of Ukraine, 47 Lenin Ave.,
 Kharkov, 61164, UKRAINE
\endaddress
\email alexandre.danilenko\@gmail.com
\endemail

\address
Faculty of Mathematics and Computer Science, Nicolaus Copernicus
University, ul. Chopina 12/18, 87-100 Toru\'n, Poland
\endaddress
\address
 and Institute of Mathematics, Polish Academy of Sciences
ul. \'Sniadeckich 8, 00-950 Warsaw, Poland
\endaddress
\email
mlem\@mat.uni.torun.pl
\endemail

\thanks
The second named author was supported in part by  Narodowe Centrum Nauki grant DEC-2011/03/B/ST1/00407.
\endthanks

\abstract
Let $H_3(\Bbb R)$ denote the 3-dimensional real  Heisenberg group.
Given a family of lattices $\Gamma_1\supset\Gamma_2\supset\cdots$ in it, let $T$ stand for the associated uniquely ergodic $H_3(\Bbb R)$-{\it odometer}, i.e. the inverse limit of the $H_3(\Bbb R)$-actions by rotations on the homogeneous spaces $H_3(\Bbb R)/\Gamma_j$, $j\in\Bbb N$.
The decomposition of the underlying Koopman unitary representation  of $H_3(\Bbb R)$ into a countable direct sum of irreducible components is explicitly described.
The ergodic 2-fold self-joinings of $T$ are found.
It is shown that in general, the $H_3(\Bbb R)$-odometers are neither isospectral nor spectrally determined.
\endabstract

\endtopmatter

\document

\NoBlackBoxes

\head 0. Introduction
\endhead

Let $T=(T_g)_{g\in G}$ be an ergodic  measure preserving action of a locally compact second countable group $G$ on a standard probability space $(X,\goth B,\mu)$.
Denote by $U_T = (U_T (g))_{g\in G} $ the associated Koopman unitary representation of $G$ in $L^2(X, \mu)$:
$$
U_T(g)f:=f\circ T_g^{-1}, \qquad f\in L^2(X,\mu).
$$
Suppose first that $G$ is Abelian.
If $U_T$ is a direct countable sum of 1-dimensional unitary sub-representations (generated by the $U_T$-eigenfunctions) then $T$ is said {\it to have a pure point spectrum}.
In 1932,  J. von Neumann \cite{Ne}  developed a theory of such actions in the case $G=\Bbb R$.
We highlight three
main aspects of this theory:
\roster
\item"(A1)" {\it isospectrality}: two ergodic flows with pure point spectrum are isomorphic if and only if  the associated Koopman unitary representations are unitarily equivalent,
\item"(A2)" {\it classification by simple algebraic invariants}:
the ergodic flows with pure point spectrum considered up to isomorphism are in one-to-one correspondence with the countable subgroups in  $\widehat{\Bbb R}$ which is the dual of $\Bbb R$,
\item"(A3)"  {\it structure}:
if an ergodic flow has pure point spectrum then it is isomorphic to a flow by rotations on a compact metric Abelian group endowed with the Haar measure.
\endroster
We also draw attention to a special subclass---apparently the simplest one---of flows with pure point {\it rational} spectrum.
They are precisely those flows which admit a structure of inverse limit of transitive flows.

The  results analogous to (A1)--(A3) hold for the general Abelian $G$ with similar proofs.
The non-Abelian case---considered in this paper---is more intricate.
G.~Mackey in \cite{Ma3} extended the concept  of pure point spectrum to actions of non-Abelian groups in the following way: $T$ has a  pure point spectrum if $U_T$ is a direct sum of countably many finite dimensional unitary representations of $G$.
He established a  structure for these actions: an ergodic action $T$ has pure point spectrum if and only if it is isomorphic to  a $G$-action by  rotations on a homogeneous space of a compact group.
However, in general, the $G$-actions with  pure point spectrum are not isospectral even in the case of finite $G$.
G.~Mackey refers to \cite{To} for a counterexample.
(See also a discussion in \cite{Le--We, Section~6}.)
Hence no classification for them is obtained.

In the present paper we consider  the case where $G$ is the 3-dimensional real Heisenberg group $H_3(\Bbb R)$ which is apparently
 the `simplest' non-Abelian nilpotent  connected Lie group.
Moreover, we single out a special class of actions of $H_3(\Bbb R)$ which we call {\it odometers}.
They are inverse limits of transitive $H_3(\Bbb R)$-actions on homogeneous spaces by lattices in $H_3(\Bbb R)$.
In this connection, we note that for   discrete finitely generated groups $G$, the $G$-odometers were considered by M.~Cortez  and S.~Petit in \cite{CoPe} in the context of topological dynamics.
We define $G$-odometers for arbitrary   locally compact second countable groups $G$ and study them as measure preserving dynamical systems.
Thus, by construction, the $H_3(\Bbb R)$-odometers are counterparts of the $\Bbb Z$-actions with pure point {\it rational} spectrum.
However, by taking into account definitions,  the Heisenberg odometers  are neither actions with pure point spectrum nor weakly mixing  (the Cartesian square of a Heisenberg odometer is not ergodic).

Our purpose here is to investigate whether   von Neumann's theory of flows with pure point spectrum extends (or partially extends)  to  the Heisenberg odometers.

We compute explicitly the spectrum of a Heisenberg odometer $T$, i.e. we describe the decomposition of $U_T$ into the direct countable sum of  its irreducible components and we calculate their multiplicities in terms of the underlying sequence of lattices in $H_3(\Bbb R)$ (Theorem~5.2).
These components consist of two families:  1-dimensional unitary representations which occur in $U_T$ with multiplicity 1 and infinitely dimensional unitary representations.
 In the non-degenerate case (i.e. where the underlying $\Bbb R^2$-odometer is not transitive) every infinitely dimensional  irreducible unitary representation  that occurs in $U_T$  has there  infinite multiplicity.
 To prove this result we use Kirillov's orbit method \cite{Ki} and Howe-Richardson spectral multiplicity  formula \cite{Ho}, \cite{Ri}.
The 1-dimensional family  is parameterized by a so-called  {\it off-rational}\footnote{See Definition~5.3.} subgroup in $\Bbb R^2$ and the infinite dimensional family is parameterized by an off-rational subgroup in $\Bbb R$.
Besides, there is a certain relation between these two subgroups.
We prove that conversely,  given two  off-rational subgroups in $\Bbb R^2$ and $\Bbb R$ connected with this relation, there is  a Heisenberg odometer whose Koopman representation is determined by these subgroups (Proposition~5.5).
We show by example that the $H_3(\Bbb R)$-odometers are not isospectral (Example~5.11).

Furthermore, we  explicitly describe all ergodic 2-fold self-joinings of the Heisenberg odometers in Theorem~6.5.
We provide an example of  an ergodic 2-fold self-joining of a transitive Heisenberg odometer which is neither transitive nor ``odometric'' at all (Example 6.4).
We note that this is in marked contrast  to the properties of locally compact transitive systems with pure point spectrum: each ergodic 2-fold self-joining of such a system is transitive itself (and hence has a pure point spectrum).

Since our definition of the Heisenberg odometers is ``structural'', it seems natural to investigate 
whether there is an equivalent  spectral definition (as in the Abelian case).
Thus we 
come to the following problem.

\roster
\item"(A4)" {\it Spectral determinacy:}
if an $H_3(\Bbb R)$-action is spectrally  equivalent to a Heisenberg odometer, is it isomorphic to a Heisenberg odometer? 
\endroster 
Answering (A3), we prove the spectral determinacy of the subclass of transitive Heisenberg odometers and, more generally, the subclass of Heisenberg odometers with transitive underlying $\Bbb R^2$-odometers (Theorem~7.1).
However,  the entire class of Heisenberg odometers is not spectrally determined.
Nevertheless, we show that if the maximal spectral type of an $H_3(\Bbb R)$-action $R$ coincides with the maximal spectral type of an $H_3(\Bbb R)$-odometer then $R$ has the structure of
a compact skew product extension of an $H_3(\Bbb R)$-action with pure point spectrum.
We obtain a criterion for when  $R$ is conjugate to a $H_3(\Bbb R)$-odometer in terms of the corresponding ``extending'' cocycle (Theorem~7.7).
Using this criterion plus the orbit theory of amenable group actions we construct an example of a probability preserving $H_3(\Bbb R)$-system which generates  the same (up to the unitary equivalence) Koopman unitary representation as an  $H_3(\Bbb R)$-odometer but which is not isomorphic to any Heisenberg odometer (Example~7.11).

\comment
It worth noting that even in the case of $\Bbb Z$-actions, the spectral determinacy is only known for the following two classes of systems: the systems with pure point spectrum and the Gaussian-Kronecker systems (see a recent survey \cite{Le}).
\endcomment

We also consider $H_3(\Bbb Z)$-odometers and compare their properties with the properties of $ H_3(\Bbb R)$-odometers.
In this connection it is interesting to note that $H_3(\Bbb Z)$ is not of type $I$ and hence neither the unitary dual  $\widehat{H_3(\Bbb Z)}$ is a standard Borel space nor the spectral theorem (i.e. a unique decomposition of an arbitrary unitary representation of  $H_3(\Bbb Z)$ into a direct  integral of irreducible representations) holds.
However for the class of $H_3(\Bbb Z)$-odometers the situation is different.
Every  such odometer has a pure point spectrum in the sense of \cite{Ma3} and  the decomposition of the corresponding Koopman unitary representations of $H_3(\Bbb Z)$ into irreducible components is well defined.
Such decompositions for a subclass of {\it normal} free $H_3(\Bbb Z)$-odometers were explicitly  computed  in \cite{Li-Ug}.
We note that the normal $H_3(\Bbb Z)$-odometers are isospectral.

The outline of the paper is as follows.
In Section~1 we discuss basic properties of $H_3(\Bbb R)$, compute the group of its automorphisms and classify the lattices in $H_3(\Bbb R)$.
The odometer actions of locally compact groups are introduced in Section~2.
Some general structural properties and the freeness of odometers are studied there.
In Section~3 a convenient criterion for the freeness of Heisenberg odometers is found.
In Section~4, for a transitive $H_3(\Bbb R)$-odometer, the decomposition of the Koopman unitary representation $U_T$ into irreducible unitary representations of $H_3(\Bbb R)$ is described.
A similar decomposition for a general $H_3(\Bbb R)$-odometer is obtained in Section~5.
The corresponding measure of maximal spectral type sits on the union of two countable subgroups (in $\Bbb R^2$ and $\Bbb R$).
All admissible pairs of such subgroups are described in Section~5.
The fact that the $H_3(\Bbb R)$-odometers are not isospectral is also proved in Section~5.
Section~6 is devoted to a description of the ergodic 2-fold self-joinings of $H_3(\Bbb R)$-odometers.
In Section~7 we investigate  spectral determinacy of the $H_3(\Bbb R)$-odometers.
In Section~8 we consider the $H_3(\Bbb Z)$-odometers.
Concluding remarks and open problems are discussed in the final Section~9.

{\it Acknowledgements.} We thank the referee for the useful remarks.

\head 1. Preliminaries on $H_3(\Bbb R)$
\endhead

We recall that  $H_3(\Bbb R)$ consists  of $3\times 3$ upper triangular matrices of the form
$$
    \pmatrix
 1 & a & c\\ 0 & 1 & b\\ 0 & 0 & 1\\
\endpmatrix,
$$
where $a,b,c$ are arbitrary reals.
The Heisenberg group endowed with the natural topology is a connected, simply-connected nilpotent Lie group.
We now let
$$
a(t):=
\pmatrix
 1 & t & 0\\ 0 & 1 & 0\\ 0 & 0 & 1\\
\endpmatrix,
\quad
b(t):=
\pmatrix
 1 & 0 & 0\\ 0 & 1 & t\\ 0 & 0 & 1\\
\endpmatrix,
\quad
c(t):=\pmatrix
 1 & 0 & t\\ 0 & 1 & 0\\ 0 & 0 & 1\\
\endpmatrix.
$$
Then $\{a(t)\mid t\in\Bbb R\}$, $\{b(t)\mid t\in\Bbb R\}$ and $\{ c(t)\mid t\in\Bbb R\}$ are three closed one-parameter subgroups in $H_3(\Bbb R)$.
The last  is the center of $H_3(\Bbb R)$.
Every element  $g$ of $H_3(\Bbb R)$  can be written uniquely as the product $g=c(t_3)b(t_2)a(t_1)$ for some $t_1,t_2,t_3\in\Bbb R$.
We also note that $[a(t_1),b(t_2)]:=a(t_1)b(t_2)a(t_1)^{-1}b(t_2)^{-1}=c(t_1t_2)$ and the commutator of $H_3(\Bbb R)$ equals the center of $H_3(\Bbb R)$.

The subgroups $H_{2,a}:=\{a(t_1)c(t_3)\mid t_1,t_3\in\Bbb R\}$
and $H_{2,b}:=\{b(t_2)c(t_3)\mid t_2,t_3\in\Bbb R\}$ are both Abelian,  normal and closed in $H_3(\Bbb R)$.
The corresponding group extensions
$$
\align
&0\to H_{2,a}\to H_3(\Bbb R)\to H_3(\Bbb R)/H_{2, a}\to 0\quad \text{and}\\
&0\to H_{2,b}\to H_3(\Bbb R)\to H_3(\Bbb R)/H_{2, b}\to 0
\endalign
$$
 both split.
 This implies that $H_3(\Bbb R)$ is isomorphic to  the semidirect product $\Bbb R^2\rtimes_B~\Bbb R$, where the homomorphism $B:\Bbb R\to \text{GL}_2(\Bbb R)$ is given by $B(t):=\pmatrix 1&t\\
0&1\endpmatrix$, $t\in\Bbb R$.
The subgroups $H_{2,a}$ and $H_{2,b}$ are automorphic in $H_3(\Bbb R)$,
i.e. there is an isomorphism $\theta$ of $H_3(\Bbb R)$ with $\theta(H_{2,a})=H_{2,b}$.
We define $\theta$ by setting $\theta(a(t)):=b(t)$, $\phi(b(t)):=a(t)$ and $\theta(c(t)):=c(-t)$ for all $t\in\Bbb R$.
To put it another way,
$$
\theta
\pmatrix
 1 & a & c\\ 0 & 1 & b\\ 0 & 0 & 1\\
\endpmatrix=
\pmatrix
 1 & b & ab-c\\ 0 & 1 & a\\ 0 & 0 & 1\\
\endpmatrix.
$$
We call $\theta$ the {\it flip} in  $H_3(\Bbb R)$.
We note that $\theta^2=\text{id}$.

The set of unitarily equivalent classes of irreducible (weakly continuous) representations of
$H_3(\Bbb R)$ is called the {\it unitary dual} of
$H_3(\Bbb R)$.
It is denoted by $\widehat {H_3(\Bbb R)}$.
The irreducible unitary representations of $H_3(\Bbb R)$ are well known.
They consist (up to  unitary equivalence) of  a family of 1-dimensional representations $\pi_{\alpha,\beta}$, $\alpha,\beta\in\Bbb R$, and a family of infinite dimensional representations $\pi_\gamma$, $\gamma\in\Bbb R\setminus\{0\}$, as follows \cite{Ki}:
$$
\align
\pi_{\alpha,\beta}(c(t_3)b(t_2)a(t_1))&:=e^{2\pi i(\alpha t_1+\beta t_2)} \quad\text{and }\\
(\pi_\gamma(c(t_3)b(t_2)a(t_1))f)(x)&:=e^{2\pi i\gamma (t_3+t_2x)}f(x+t_1),\quad f\in L^2(\Bbb R,\lambda_{\Bbb R}).
\endalign
$$
Thus we can identify $\widehat{H_3(\Bbb R)}$  with the disjoint union $\Bbb R^2\sqcup \Bbb R^*$.
We recall that there is a natural Borel structure on the  unitary dual of each locally compact second countable group \cite{Ma2}.
 In the case of the Heisenberg group this Borel $\sigma$-algebra coincides with the standard $\sigma$-algebra of Borel subsets in  $\Bbb R^2\sqcup\Bbb R^*$.

Given an arbitrary unitary representation $U=(U(g))_{g\in H_3(\Bbb R)}$
of $H_3(\Bbb R)$ in a separable Hilbert space $\Cal H$, there are  a measure $\sigma_U$ on $\widehat{H_3(\Bbb R)}$ (i.e. two measures $\sigma_U^{1,2}$ on $\Bbb R^2$ and $\sigma_U^3$ on $\Bbb R^*$) and a map $l_U:\widehat{H_3(\Bbb R)}\to\Bbb N\cup\{\infty\}$  (i.e. two maps $l_U^{1,2}:\Bbb R^2\ni (x,y)\mapsto l_U^{1,2}(x,y)\in\Bbb N\cup\{\infty\}$
 and
$l_U^3:\Bbb R^*\ni z\mapsto l_U^3(z)\in\Bbb N\cup\{\infty\}$) such that the following decompositions hold (up to the unitary equivalence):
$$
\align
\Cal H&= \int^{\oplus}_{\Bbb R^2}\,\bigoplus_{j=1}^{l_U^{1,2}(\alpha,\beta)}\Bbb C\,d\sigma_U^{1,2}(\alpha,\beta)
\oplus\int_{\Bbb R^*}^\oplus\bigoplus_{j=1}^{l_U^3(\gamma)}L^2(\Bbb R,\lambda_\Bbb R)\,d\sigma_U^3(\gamma)\quad\text{and} \\
U(g)&=\int^{\oplus}_{\Bbb R^2}\,\bigoplus_{j=1}^{l_U^{1,2}(\alpha,\beta)} \pi_{\alpha,\beta}(g)\,d\sigma_U^{1,2}(\alpha,\beta)
\oplus\int_{\Bbb R^*}^\oplus\bigoplus_{j=1}^{l_U^3(\gamma)}\pi_\gamma(g)\,d\sigma_U^3(\gamma).
\endalign
$$
The equivalence class of  $\sigma_U$ is called the {\it maximal spectral type} of $U$.
The map $l_U$ is called the  {\it multiplicity function} of $U$.
The essential range of $l_U$ is called the {\it set of spectral multiplicities} of $U$.
The maximal spectral type
and the multiplicity function of $U$ ($\sigma_U$-mod 0)
are both determined uniquely by the unitary equivalence class of $U$.

Below, in this section  we:
\roster
\item"(a)" describe explicitly the group of all continuous automorphisms of the Heisenberg group and
\item"(b)" classify the lattices in $H_3(\Bbb R)$ up to automorphism of $H_3(\Bbb R)$.
\endroster
Although these results are not new, we think it is easier to reprove them here than to find them in the literature and ``adjust'' to our notation.

From now on we denote by $p$ the natural projection
$$
H_3(\Bbb R)\ni g\mapsto p(g):=(t_1,t_2)\in\Bbb R^2
$$
whenever $g=c(t_3)b(t_2)a(t_1)$.

Let $\widetilde\theta$ be an automorphism of $H_3(\Bbb R)$.
Of course, $\widetilde\theta$ preserves  the center of $H_3(\Bbb R)$.
Hence there exists $\tau\in\Bbb R^*$ such that
$$
\widetilde\theta(c(t))=c(\tau t), \quad t\in\Bbb R.\tag1-1
$$
Since the center of $H_3(\Bbb R)$ is the kernel of $p$, $\widetilde\theta$ passes through $p$, i.e. there is a matrix $A=\pmatrix \xi_1 &\eta_1\\ \xi_2 & \eta_2 \endpmatrix\in  GL_2(\Bbb R)$
such that $p\circ\widetilde\theta=A\circ p$.
Therefore we can write $\widetilde\theta(a(t))$ as
$$
\widetilde\theta(a(t))=c(f(t))b(\xi_2t)a(\xi_1t),\quad t\in\Bbb R,
$$
for some continuous function $f:\Bbb R\to\Bbb R$.
Since $\widetilde\theta(a(t_1+t_2))=\widetilde\theta(a(t_1))\widetilde\theta(a(t_2))$,
it follows that
$$
f(t_1+t_2)=f(t_1)+f(t_2)+\xi_1\xi_2t_1t_2
$$
for all $t_1,t_2\in\Bbb R$.
Therefore $f(t)=\frac{\xi_1\xi_2}2t^2+\xi t$ for some real
$\xi\in\Bbb R$.
Thus,
$$
\widetilde\theta(a(t))=c\bigg(\frac{\xi_1\xi_2}2t^2+\xi t\bigg)
b(\xi_2t)a(\xi_1t),\quad t\in\Bbb R.\tag1-2
$$
In a similar way,
$$
\widetilde\theta(b(t))=c\bigg(\frac{\eta_1\eta_2}2t^2+\eta t\bigg)
b(\eta_2t)a(\eta_1t),\quad t\in\Bbb R,\tag1-3
$$
for some $\eta \in\Bbb R$.
Since $\widetilde\theta$ is a homomorphism of $H_3(\Bbb R)$, we  have
$$
\widetilde\theta(a(t_1))\widetilde\theta(b(t_2))=
\widetilde\theta(a(t_1)b(t_2))=\widetilde\theta(b(t_2)a(t_1)c(t_1t_2))=
\widetilde\theta(b(t_2))\widetilde\theta(a(t_1))\widetilde\theta(c(t_1t_2))
$$
for all $t_1,t_2\in\Bbb R$.
Applying \thetag{1-2} and \thetag{1-3} and then \thetag{1-1}, we obtain
$$
c(\xi_1\eta_2t_1t_2)=c(\eta_1\xi_2t_1t_2+\tau t_1t_2),
$$
which yields $\tau=\det A$.
Thus, the three parameters $A,\xi,\eta$ determine $\widetilde\theta$ completely.
Conversely, given $A\in GL_2(\Bbb R)$ and  $\xi,\eta\in\Bbb R$, the formulas \thetag{1-1}--\thetag{1-3} with $\tau=\det A$ determine completely an automorphism of $H_3(\Bbb R)$.
We will denote it by $\theta_{A,\xi,\eta}$.

Thus, we have proved the following proposition.

\proclaim{Proposition 1.1}
\rom{Aut}$(H_3(\Bbb R))=\{\theta_{A,\xi,\eta}\mid A\in GL_2(\Bbb R),\ \xi,\eta\in\Bbb R\}$.
\endproclaim

We also note that $N:=\{\theta_{I,\xi,\eta}\mid \xi,\eta\in\Bbb R\}$
is a normal subgroup in \rom{Aut}$(H_3(\Bbb R))$.
It is isomorphic naturally to $\Bbb R^2$.
In fact, $N$ is the subgroup of inner automorphisms of $H_3(\Bbb R)$.
We obtain a short exact sequence
$$
\{1\}\leftarrow GL_2(\Bbb R)\leftarrow
\text{\rom{Aut}}(H_3(\Bbb R))\leftarrow N\leftarrow \{1\}.\tag1-4
$$
Moreover, a  direct calculation shows  that
the homomorphism
$$
GL_2(\Bbb R)\ni A\mapsto\theta_{A,0,0}\in \text{\rom{Aut}}(H_3(\Bbb R))
$$
is a  cross-section of the natural projection in \thetag{1-4}.

\proclaim{Corollary 1.2}
 $\text{\rom{Aut}}(H_3(\Bbb R))$ is isomorphic to the semidirect product $\Bbb R^2\rtimes GL_2(\Bbb R)$, where
 the corresponding action of $GL_2(\Bbb R)$ on $\Bbb R^2$ is given by $Av:=\det A \cdot (A^*)^{-1}v$.
\endproclaim

We now describe the structure of lattices in the Heisenberg group and classify them up to  group automorphism.
Recall that a lattice is  a discrete subgroup of finite covolume.
In the case of a simply connected nilpotent Lie group, every lattice is cocompoct \cite{Ra}.
Fix a lattice $\Gamma$  in $H_3(\Bbb R)$.
There is a real $\xi_\Gamma>0$\footnote{If $\xi_\Gamma=0$ then the intersection of $\Gamma$ with the center of $H_3(\Bbb R)$ is trivial.
This yields  that $\Gamma$ is Abelian and hence there is a line (i.e. a 1-dimensional subspace) $\Cal L$ in $\Bbb R^2$ such that $\Gamma$ is contained in a normal closed subgroup $H:=\{c(t)b(t_2)a(t_1)\mid (t_1,t_2)\in\Cal L, t\in\Bbb R\}$ of $H_3(\Bbb R)$.
Therefore the quotient group $H_3(\Bbb R)/H$  being a quotient space of $H_3(\Bbb R)/\Gamma$ is of finite Haar measure, i.e. compact. However, it is straightforward to verify that $H_3(\Bbb R)/H$ is isomorphic to $\Bbb R$, a contradiction.} such that
$$
\Gamma\cap\{c(t)\mid t\in\Bbb R\}=\{c(m\xi_\Gamma)\mid m\in\Bbb Z\}.
$$
The central extension
$$
\{0\}\leftarrow\Bbb R^2\overset{p}\to\leftarrow H_3(\Bbb R)\overset{c}\to\leftarrow\Bbb R\leftarrow \{0\}
$$
induces a short exact sequence
$$
\{0\}\longleftarrow p(\Gamma)\overset{p}\to\longleftarrow\Gamma\overset{c}\to\longleftarrow
\xi_\Gamma\Bbb Z\longleftarrow\{0\}.
$$
We note that $p(\Gamma)$ is a lattice in $\Bbb R^2$.\footnote{Indeed, the flow $F:\Bbb R\times H_3(\Bbb R)/\Gamma\ni (t,g\Gamma)\mapsto c(t)g\Gamma\in H_3(\Bbb R)/\Gamma$ is periodic with period $\xi_\Gamma$. Hence the Borel space $O$ of $F$-orbits is standard. On the other hand, $O$ is isomorphic as a Borel space to $H_3(\Bbb R)/p^{-1}(p(\Gamma))$. The later space is, in turn, isomorphic to the Borel space $\Bbb R^2/p(\Gamma)$. Hence $\Bbb R^2/p(\Gamma)$ endowed with the natural quotient Borel structure is standard. Hence $p(\Gamma)$ is closed in $\Bbb R^2$ by \cite{Ma2}.}
Therefore
 there is a matrix $A\in GL_2(\Bbb R)$ such that $p(\theta_{A,0,0}(\Gamma))=Ap(\Gamma)=\Bbb Z^2$.
The commutator subgroup $[\Gamma,\Gamma]$ is of a finite index $k_\Gamma>0$ in $p^{-1}(0)\cap\Gamma$.

It is easy to see that $k_{\theta(\Gamma)}=k_\Gamma$ for each $\theta\in\text{Aut}(H_3(\Bbb R))$.
Let $\zeta$ and $\eta$ be the smallest non-negative  reals such that $c(\zeta)b(0)a(1)\in\theta_{A,0,0}(\Gamma)$ and $c(\eta)b(1)a(0)\in\theta_{A,0,0}(\Gamma)$.
Then $\theta_{I,-\zeta,-\eta}(c(\zeta)a(1))=a(1)$ and
$\theta_{I,-\zeta,-\eta}(c(\eta)a(1))=b(1)$.
We let $\widetilde\Gamma:=\theta_{I,-\zeta,-\eta}(\theta_{A,0,0}(\Gamma))$.
Then $a(1),b(1)\in\widetilde\Gamma$.
Hence $[a(1),b(1)]=c(1)\in\widetilde\Gamma$.
On the other hand, $[a(1),b(1)]=c(k_{\widetilde\Gamma}\xi_{\widetilde\Gamma})
=c(k_{\Gamma}\xi_{\widetilde\Gamma})$.
It follows that $\xi_{\widetilde\Gamma}=1/k_\Gamma$.
The elements $a(1),b(1)$ and $c(k_\Gamma^{-1})$ generate the lattice $\widetilde\Gamma$.
 Thus we  have proved the following proposition.

\proclaim{Proposition 1.3}
Given a lattice $\Gamma$ in $H_3(\Bbb R)$, there is an automorphism $\theta$ of $H_3(\Bbb R)$ such that
$$
\theta(\Gamma)=\left\{\pmatrix 1&l&\frac n{k_\Gamma}\\
0&1&m\\
0&0&1
\endpmatrix \Bigg|\, l,m,n\in\Bbb Z\right\}.
$$
Hence two lattices $\Gamma_1$ and $\Gamma_2$ in $H_3(\Bbb R)$ are automorphic
if and only if $k_{\Gamma_1}=k_{\Gamma_2}$.
Two lattices $\Gamma_1$ and $\Gamma_2$ in $H_3(\Bbb R)$ are conjugate
if and only  if $k_{\Gamma_1}=k_{\Gamma_2}$ and  $p(\Gamma_1)=p(\Gamma_2)$.
\endproclaim

\proclaim{Corollary 1.4}
If an element $g\in H_3(\Bbb R)$ commutes with each element of a lattice $\Gamma$ in $H_3(\Bbb R)$ then $g$ belongs to the center of $H_3(\Bbb R)$.
\endproclaim

The following  lemma will be used in the proof of  Theorem~5.7 which is one of the main results of this paper.

\comment

\proclaim{Lemma 1.4?}
Let $\Gamma,\Gamma_0,\Gamma_1$ and $\Gamma_2$ be four lattices in $H_3(\Bbb R)$ such that  $\Gamma_1\cup\Gamma_2\subset\Gamma_0\subset\Gamma$.
If $\Gamma_1$ and $\Gamma_2$ are conjugate in $H_3(\Bbb R)$
and every element of $p(\Gamma_1)$ has a root of order $k_{\Gamma_{0}}\cdot\#(p(\Gamma_0)/p(\Gamma_1))$ in $p(\Gamma)$
then  $\Gamma_1$ and $\Gamma_2$  are also conjugate in $\Gamma$.
\endproclaim
\demo{Proof}
In view of Proposition~1.3, we may assume without loss of generality that $\Gamma_0=\{c(n/k_{\Gamma_0})b(m)a(l)\mid n,m,l\in\Bbb Z\}$.
In particular, $p(\Gamma_0)=\Bbb Z^2$.
Since $\Gamma_1$ and $\Gamma_2$ are conjugate, $p(\Gamma_1)=p(\Gamma_2)$.
Since $p(\Gamma_1)\subset\Bbb Z^2$, there is a matrix $A\in GL_2(\Bbb R)$ with integer entries such that  $A(\Bbb Z^2)=p(\Gamma_1)$.
Reducing $A$ to the Smith normal form we find two positive integers $q_1$ and $q_2$ with $q_1 | q_2$ and two matrices $L,R\in GL_2(\Bbb Z)$ such that
$$
L^{-1}AR=\pmatrix
q_1 & 0\\
0 & q_2
\endpmatrix.
$$
It  follows from the last statement of Proposition~1.3 that the lattices $\theta_{L,0,0}(\Gamma_0)$
and $\Gamma_0$ are conjugate.
Therefore there are $\xi,\eta\in\Bbb R$ such that
$\theta_{L,\xi,\eta}(\Gamma_0)=\Gamma_0$.
Replacing the quadruple $\Gamma,\Gamma_0,\Gamma_1,\Gamma_2$ with the quadruple
$\theta_{L,\xi,\eta}(\Gamma), \Gamma_0,   \theta_{L,\xi,\eta}(\Gamma_1), \theta_{L,\xi,\eta}(\Gamma_2)$, we will assume that
$$
p(\Gamma_1)=p(\Gamma_2)=\{(lq_1, mq_2)\mid m,l\in\Bbb Z\}\subset\Bbb Z^2.
$$
Let $t_{1,j},t_{2,j}$
be the smallest positive reals such that
$c(t_{1,j})a(q_1)\in \Gamma_j$ and
$c(t_{2,j})b(q_2)\in \Gamma_j$, $j=1,2$.
Since $\Gamma_1\cup\Gamma_2\subset\Gamma_0$,
it follows that $t_{1,j},t_{2,j}\in k^{-1}_{\Gamma_0}\Bbb Z$.
We now let
$$
\xi_1:=\frac{t_{1,1}-t_{1,2}}{q_1}\quad\text{and}\quad \xi_2:=\frac{t_{2,2}-t_{2,1}}{q_2}.
$$
Since $q_1q_2= \#(p(\Gamma_0)/p(\Gamma_1))$, it follows from the condition of the lemma that there is $t_0\in\Bbb R$ such that an element $\gamma:=c(t_0)a(\xi_2)b(\xi_1)$ belongs to
$\Gamma$.
It is straightforward to verify that $\gamma\Gamma_1\gamma^{-1}=\Gamma_2$.
\qed
\enddemo

\endcomment

\proclaim{Lemma 1.5}
Let $\Gamma$ and $\Lambda$ be two lattices in $H_3(\Bbb R)$ such that $\Lambda\subset\Gamma$.
Let $p(\Gamma)=A(\Bbb Z^2)$ for a matrix $A\in GL_2(\Bbb R)$.
Then for each positive real $\xi$ such that $\xi/\xi_\Gamma\in\Bbb N$, $\xi_\Lambda/\xi\in\Bbb N$ and $|\det A|/\xi\in\Bbb N$, there is a lattice $\Theta$ such that $\Gamma\supset\Theta\supset\Lambda$, $p(\Theta)=p(\Gamma)$
and $\xi_\Theta=\xi$.
\endproclaim

\demo{Proof}
In view of Proposition 1.3, we may assume without loss of generality that
$\Gamma=\{c(j_3/k)b(j_2)a(j_1)\mid j_1,j_2,j_3\in\Bbb Z\}$ for some $k\in\Bbb N$.
It follows from the condition of the lemma that
there are positive integers $k_1,k_2,k_3$
such that
$k=k_1k_2$, $\xi=1/k_2$ and $\xi_\Lambda=k_3/k_2$.
Let $\Theta$ denote the group generated by $\Lambda$ and the lattice
$\{c(j_3/k_2)b(j_2)a(j_1)\mid j_1,j_2,j_3\in\Bbb Z\}$.
Then $\Theta$ is a lattice in $H_3(\Bbb R)$.
Of course,
$\Gamma\supset\Theta\supset\Lambda$ and $p(\Theta)=p(\Gamma)$.
It is easy to see that
$$
\xi_\Theta\Bbb Z=\xi_\Lambda\Bbb Z+\frac 1{k_2}\Bbb Z=\frac 1{k_2}\Bbb Z,
$$
i.e. $\xi_\Theta=\xi$.
\qed
\enddemo

\remark{Remark 1.6} It is straightforward to verify that
for each lattice $\Gamma$ in $H_3(\Bbb R)$, we have
  $k_\Gamma=\frac{|\det A|}{\xi_\Gamma}$, where the matrix $A$ is defined by $p(\Gamma)=A(\Bbb Z^2)$.
  Indeed, since $[\Gamma,\Gamma]=\{c(k_\Gamma\xi_\Gamma m)\mid m\in \Bbb Z\}$,
  it follows that $[\gamma_1,\gamma_2]= c(\pm k_\Gamma\xi_\Gamma)$ for each pair of elements $\gamma_1,\gamma_2\in\Gamma$ whose images $p(\gamma_1)$ and $p(\gamma_2)$ in $\Bbb R^2$ generate the entire lattice $p(\Gamma)$.
\endremark

\comment
Let $\widetilde\Gamma_1\subset\Gamma_1$
and $\widetilde\Gamma_2\subset\Gamma_2$ be four lattices in $H_3(\Bbb R)$.
Suppose that $\Gamma_1=g\Gamma_2g^{-1}$ and $\widetilde\Gamma_1=\widetilde g\widetilde\Gamma_2\widetilde g^{-1}$ for some $g,\widetilde g\in H_3(\Bbb R)$.
Then there is $\gamma\in\Gamma_2$ such that
$\widetilde\Gamma_1=g\gamma\widetilde\Gamma_2\gamma^{-1} g^{-1}$.
\endcomment

\remark{Remark 1.7} There is a natural action of Aut$(H_3(\Bbb R))$ on $\widehat{H_3(\Bbb R)}$.
It is not difficult  to compute that $\pi_{\alpha,\beta}\circ\theta_{A,\xi,\eta}=\pi_{\xi_1\alpha+\xi_2\beta, \eta_1\alpha+\eta_2\beta}$ and $\pi_\gamma\circ\theta_{A,\xi,\eta}\sim\pi_{\gamma\cdot\det A}$.
\endremark

\head 2. Odometer actions of locally compact groups
\endhead

In this section we introduce the concept of  odometer action for locally compact second countable groups.
We study some general structural properties of these actions and discuss when they are free.

Let $G$ and $K$ be  locally compact second countable groups.
Let $T=(T_g)_{g\in G}$ be an action of $G$  on a standard Borel space $Z$.
A Borel map $c:G\times Z\to K$ is called a cocycle of $T$ if
$$
c(g_1g_2,z)=c(g_1,T_{g_2}z)c(g_2,z)
$$
for all $g_1,g_2\in G$ and $z\in Z$.
Two cocycles $c,c':G\times Z\to K$ are {\it cohomologous} if there is a Borel map
$d:Z\to K$ such that
$d(gz)c(g,z)d(z)^{-1}=c'(g,z)$ for all  $g\in G$ and $z\in Z$.
Given a cocycle $c$, we can define an action $T^c=(T^c_g)_{g\in G}$ of $G$ on $Z\times K$ by setting
$$
T^c_g(z,k)=(T_gz,c(g,z)k).
$$
It is called the $c$-skew product extension of $T$.
If $c$ is cohomologous to another cocycle $c'$ of $T$ then the actions $T^c$ and $T^{c'}$ of $G$ are isomorphic.
If $T$ preserves a $\sigma$-finite measure $\mu$ then $T^c$ preserves the $\sigma$-finite measure $\mu\times\lambda_K$, where $\lambda_K$ denotes a left Haar measure on $K$.
If $T^c$ is ergodic then $c$ is called {\it ergodic}.

Let $N$ be a closed  subgroup of $G$.
Let the two groups be unimodular.
Fix a Borel cross-section $s:G/N\to G$ of the natural projection $G\to G/N$ such that $s(N)=1$.
Then a cocycle $h_s:G\times G/N\to N$  of the natural $\sigma$-finite (Haar) measure preserving $G$-action on $G/N$ by rotations is well defined by the formula
$h_s(g,x):=s(gx)^{-1}gs(x)$.
The cohomology class of $h_s$ does not depend on the choice of $s$.
Below we will call $h_s$  a {\it choice-cocycle}.
It is easy to see that the $h_s$-skew product extension is isomorphic to the action of $G$ on itself by left translations.

\remark{Remark 2.0} Given a cocycle $c:G\times G/N\to K$,
a continuous group homomorphism $\phi:N\to K$ is well defined by
$\phi(n):=c(n,N)$.
Consider a mapping $d:G/N\ni x\mapsto d(x):=c(s(x),N)\in K$.
Then we have
$$
\align
c(g_1,g_2N)&=c(g_1g_2,N)c(g_2,N)^{-1}\\
&=
c(s(g_1g_2N)h_s(g_1g_2,N),N)c(s(g_2N)h_s(g_2,N),N)^{-1}\\
&=d(g_1g_2N)\phi(h_s(g_1g_2,N)h_s(g_2,N)^{-1})d(g_2N)^{-1}\\
&=d(g_1g_2N)\phi(h_s(g_1,g_2N))d(g_2N)^{-1}.
\endalign
$$
In other words,  $c$ is cohomologous to the cocycle $\phi\circ h_s$.
Thus, the cohomology class of $c$ is determined completely by $\phi$ (or, more precisely, by the conjugacy class of $\phi$).
\endremark

\remark{Remark 2.1}
If $K$ is Abelian then the essential range\footnote{For the definition of the essential range of a cocycle we refer to \cite{Sc}.} $E(c)$ of $c$ is the closure of the subgroup $\phi(N)$ in $K$.
Indeed, this follows from the fact that $h_s$
is ergodic (in fact, the corresponding $h_s$-skew product is transitive).
\endremark

Let $\Gamma_1\supset \Gamma_2\supset\cdots$ be a nested sequence of lattices in $G$.
Then we consider a projective sequence of homogeneous $G$-spaces
$$
G/\Gamma_1\leftarrow G/\Gamma_2\leftarrow\cdots.\tag2-1
$$
All arrows are $G$-equivariant and onto.
Denote by $X$ the projective limit of this sequence.
Then $X$ is a locally compact second countable $G$-space.
Indeed,  $G/\Gamma_1$ is locally compact and every arrow in \thetag{2-1}
is a finite-to-one mapping.\footnote{We note that $\Gamma_n$ is of finite index in $\Gamma_1$.}
 The $G$-action is minimal and uniquely ergodic.
The only  invariant probability measure $\mu$ on $X$ is the projective limit
of the probability Haar measures on $G/\Gamma_n$.

\definition{Definition 2.2}
We call the dynamical system $(X,\mu,G)$ a $G$-{\it odometer}\footnote{Thus the $G$-odometers are ergodic. Cf. \cite{CoPe}, where  the case of  discrete and finitely generated $G$ was under consideration.}.
If $\Gamma_n$ is normal in $\Gamma_1$ for each $n$ then we call the corresponding $G$-odometer {\it normal}.
\enddefinition

We  note
that $X$ is compact if and only if each $\Gamma_n$ is {\it uniform}, i.e. co-compact in~$G$.

To show a `product structure' on $X$ we
 denote by $s_j:\Gamma_{j-1}/\Gamma_j\to \Gamma_{j-1}$ a Borel cross-section of the natural projection $\Gamma_{j-1}\to Y_j:=\Gamma_{j-1}/\Gamma_j$ such that $s_j(\Gamma_j)=1$.
 For consistency of the notation we set $\Gamma_0=G$.
Let $h_j:\Gamma_{j-1}\times Y_j\ni(g,y)\mapsto h_j(g,y)\in\Gamma_j$ stand for the corresponding choice-cocycle.
Then
$gs_j(y)=s_j(gy)h_j(g,y)$ for all $g\in\Gamma_{j-1}$ and $y\in Y_j$.
Denote by $\psi_j$ a mapping
$$
G/\Gamma_1\times\Gamma_1/\Gamma_2\times\cdots\times\Gamma_{j-1}/\Gamma_j
\ni (y_1,\dots,y_j)\mapsto s_1(y_1)s_2(y_2)\cdots s_j(y_j)\Gamma_j\in G/\Gamma_j
$$
for $j>1$ and $\psi_1:=\text{id}$.
It is easy to see that each $\psi_j$ is a Borel isomorphism.
We obtain the following  infinite  commutative diagram
$$
\CD
G/\Gamma_1 @<<< G/\Gamma_2 @<<< G/\Gamma_3 @<<<\cdots\\
@A{\psi_1}AA @A{\psi_2}AA @A{\psi_3}AA  \cdots\\
Y_1 @<<< Y_1\times Y_2@<<<
Y_1\times Y_2\times Y_3@<<<\cdots\\
\endCD,
$$
where the horizontal arrows in the lower lines are natural projections.
Therefore the sequence of mappings $(\psi_j^{-1})_{j=1}^\infty$
generates a Borel isomorphism $\psi$ of $X:=\projlim_{j}G/\Gamma_j$
onto $Y:=Y_1\times Y_2\times\cdots$.
Since $X$ is a $G$-space, $Y$ is also a $G$-space.
We thus obtain the following proposition.

\proclaim{Proposition 2.3} The $G$-space $X$ is Borel isomorphic to the infinite Cartesian product
$$
Y=G/\Gamma_1\times\Gamma_1/\Gamma_2\times\Gamma_2/\Gamma_3\times\cdots.
$$
in such a way that
$$
g\cdot(y_1,y_2,y_3,\dots)=(gy_1,h_1(g,y_1)y_2, h_2(h_1(g,y_1),y_2)y_3,\dots)
$$
for all $(y_1,y_2\dots)\in Y$ and $g\in G$.
\endproclaim

We recall that given a lattice $\Gamma$ in $G$ and a measure preserving action $V=(V_\gamma)_{\gamma\in\Gamma}$ of $\Gamma$ on a standard probability space $(Z,\kappa)$, there is a natural {\it induced} measure preserving $G$-action \cite{Ma1} on the product space $(G/\Gamma\times Z,\lambda\times\kappa)$:
$$
g(x,z):=(gx,V_{h_s(g,x)}z),
$$
where $\lambda$ is the probability Haar measure on $G/\Gamma$ and $h_s:G\times G/\Gamma\to \Gamma$ is the 1-cocycle corresponding to a cross-section $s:G/\Gamma\to G$.
We denote this action by Ind$_{\Gamma}^G(V)$.
The isomorphism class of Ind$_{\Gamma}^G(V)$ does not depend on a particular choice of $s$.

Let  $K$ be a compact metric group and let $\alpha:\Gamma\times Z\to K$ be a cocycle of $V$.
Then we can define an {\it induced} cocycle $\widetilde\alpha:G\times(G/\Gamma\times Z)\to K$ of Ind$_\Gamma^G(V)$ by setting
$$
\widetilde\alpha(g,(x,z)):=\alpha(h_s(g,x),z).
$$
The next lemma is about  basic properties of induced cocycles.
They follow easily from the definitions
and we state this lemma  without proof.

\proclaim{Lemma 2.4}
\roster
\item"\rom{(i)}"
\rom{Ind}$_\Gamma^G(V^\alpha)=(\text{\rom{Ind}}_\Gamma^G(V))^{\widetilde\alpha}$.
In particular, $ \widetilde\alpha$ is ergodic if and only if $\alpha$ is ergodic.
\item"\rom{(ii)}"
Two cocycles $\alpha,\beta$ of $V$ are cohomologous if and only if the induced cocycles $\widetilde\alpha$ and $\widetilde\beta$ of \rom{Ind}$_\Gamma^G(V)$ are cohomologous.
\item"\rom{(iii)}"
For each cocycle $\delta$ of
 \rom{Ind}$_\Gamma^G(V)$, there exists a cocycle $\alpha$ of $V$ such that $\delta$ is cohomologous to $\widetilde\alpha$.
 \item"\rom{(iv)}"
  Let $\goth F$ be a factor of $V$ and let $\widetilde{\goth F}:=\text{\rom{Ind}}_\Gamma^G(\goth F)$ denote the
  corresponding factor of \rom{Ind}$_\Gamma^G(V)$.
  Then $\widetilde\alpha$ is  $\widetilde{\goth F}$-measurable 
  if and only if $\alpha$ is  $\goth F$-measurable.
 \endroster
\endproclaim

The next statement follows easily from Proposition~2.3.

\proclaim{Corollary 2.5}
Let $T$ be a $G$-odometer associated with a sequence
$\Gamma_1\supset\Gamma_2\supset\dots$  of lattices in $G$.
Let $V$ denote  the $\Gamma_1$-odometer associated with the sequence
$\Gamma_2\supset \Gamma_3\supset\cdots$ of lattices in $\Gamma_1$.
Then $\text{\rom{Ind}}^G_{\Gamma_1}(V)=T$.
\endproclaim

We now find some conditions under which a $G$-odometer  is free.
Let $X$ denote the space of an  odometer $T$.
Every point $x\in X$ can be written as an infinite sequence
$x=(g_1\Gamma_1,g_2\Gamma_2,\dots)$ for some $g_1\in G$, $g_2\in g_1\Gamma_1$, $g_3\in g_2\Gamma_2$, \dots.
Denote by $S_x$ the stability group of $T$ at $x$.

\proclaim{Proposition 2.6}
Let $x=(g\Gamma_1,g\gamma_1\Gamma_2,g\gamma_1\gamma_2\Gamma_3,\dots)\in X$
for some $g\in H_3(\Bbb R)$ and $\gamma_n\in
\Gamma_n$, $n\in\Bbb N$.
Then
$S_x=g \bigg(\bigcap_{n=1}^\infty
\gamma_1\cdots\gamma_{n-1}\Gamma_n \gamma_{n-1}^{-1}\cdots
\gamma_1^{-1}\bigg)  g^{-1}$.
Therefore  $T$ is free if and only if
$$
\bigcap_{n=1}^\infty
\gamma_1\cdots\gamma_{n-1}\Gamma_n\gamma_{n-1}^{-1}\cdots
\gamma_1^{-1}=\{1\}
$$
for a.e. $x$, i.e. for a.e. sequence $(\gamma_1\Gamma_2,\gamma_2\Gamma_3,\dots)\in\prod_{n\in\Bbb N}\Gamma_n/\Gamma_{n+1}$, where the infinite product space is endowed with the infinite product measure obtained from   equi-distributions on the factors.
Hence $T$ is free if and only if the $\Gamma_1$-odometer associated with the sequence $\Gamma_2\supset\Gamma_3\supset\cdots$ is free.
The restriction of $T$ to the center of $G$ is   free if and only if the intersection of  the center with the group $\bigcap_{n=1}^\infty\Gamma_n$ is trivial.
\endproclaim
\demo{Proof}
Since
$T_hx=(hg\Gamma_1,hg\gamma_1\Gamma_2,hg\gamma_1\gamma_2\Gamma_3,\dots)$, it follows that
 $T_hx=x$ if and only if $h\in\bigcap_{n=1}^\infty g\gamma_1\cdots\gamma_{n-1}\Gamma_n \gamma_{n-1}^{-1}\cdots\gamma_1^{-1}g^{-1}$.
\qed
\enddemo

\proclaim{Corollary 2.7} Let $(X,\mu,T)$ be the odometer associated with a sequence $\Gamma_1\supset\Gamma_2\supset\cdots$ of lattices in $G$.
\roster
\item"\rom{(i)}"
If $\Gamma_n$ is normal in $\Gamma_1$ for each $n$, then $T$ is free if and only if $\bigcap_{n=1}^\infty\Gamma_n=\{1\}$.
\item"\rom{(ii)}"
Let $\widetilde\Gamma_n:=\bigcap_{\gamma\in\Gamma_1}\gamma\Gamma_n\gamma^{-1}$. Then $\widetilde\Gamma_n$ is a subgroup of finite index in $\Gamma_n$.
It is normal in $\Gamma_1$.
Moreover,
$\Gamma_1\supset \widetilde\Gamma_2\supset\widetilde\Gamma_3\supset\cdots$.
Denote by $(\widetilde X,\widetilde\mu,\widetilde T)$ the corresponding normal $G$-odometer.
We call it the normal cover of $(X,\mu,T).$
The sequence of natural $G$-equivariant maps $G/\widetilde \Gamma_n\to G/\Gamma_n$ generates  a
factor map $(\widetilde X,\widetilde\mu,\widetilde T)\to (X,\mu,T)$.
Thus every odometer is a factor of a  normal odometer.
 Moreover, if $\bigcap_{n=1}^\infty\Gamma_n=\{1\}$ then
 $\bigcap_{n=1}^\infty\widetilde\Gamma_n\subset \bigcap_{n=1}^\infty\Gamma_n$ and the corresponding normal cover is free\footnote{See also \cite{CoPe} for a analogous claim for a subclass of discrete $G$.}.
Furthermore, $\bigcap_{n=1}^\infty\widetilde\Gamma_n\cap\{c(t)\mid t\in\Bbb R\}=\bigcap_{n=1}^\infty\Gamma_n\cap\{c(t)\mid t\in\Bbb R\}$.
\endroster
\endproclaim

\remark{Remark 2.8}
If there is a F{\o}lner sequence $(F_n)_{n=1}^\infty$ in $G$ such that
$G=\bigsqcup_{\gamma\in\Gamma_n}F_{n-1}\gamma$
and $F_{n}=\bigsqcup_{\gamma\in C_n}F_{n-1}\gamma$
for some finite subset $C_n\subset\Gamma_n$ and each $n>0$ then the corresponding $G$-odometer is the $(C,F)$-action associated with $(F_{n-1},C_n)_{n=1}^\infty$ (see \cite{Da3}).
Every $(C,F)$-action is free.
We note that $\Gamma_n=C_n\Gamma_{n+1}$ and $c\Gamma_{n+1}\cap c'\Gamma_{n+1}=\emptyset$ whenever $c\ne c'\in C_n$.
Hence $\# C_n$ is the index of $\Gamma_{n+1}$ in $\Gamma_n$.
\endremark

\head 3. Freeness for  Heisenberg odometers
\endhead
In this section we consider in detail the problem of when $H_3(\Bbb R)$-odometers are free.
 Some criteria of the freeness are found.
 A couple of counterexamples are also given here.

\proclaim{Theorem 3.1}
Let $T$ be the $H_3(\Bbb R)$-odometer associated with a sequence $\Gamma_1\supset\Gamma_2\supset\cdots$ of lattices in $H_3(\Bbb R)$.
Then $T$ is free if and only if
$\{c(t)\mid t\in\Bbb R\}\cap\bigcap_{n=1}^\infty\Gamma_n=\{1\}$.
\endproclaim
\demo{Proof}
Let $(X,\goth B,\mu)$ denote  the (standard probability) space of $T$.
Suppose that
$\{c(t)\mid t\in\Bbb R\}\cap\bigcap_{n=1}^\infty\Gamma_n=\{1\}$.
It follows from Proposition~2.6 that  the restriction of $T$ to the center of $H_3(\Bbb R)$ is free.
Let $\goth S$ stand for the space of closed subgroups in  $H_3(\Bbb R)$.
It is a standard Borel space when endowed with the Fell topology \cite{Fe}.
Denote by $S_x$ the stability group of $T$ at a point $x\in X$.
If follows from Proposition~2.6 that $S_x$ is countable.
The map $X\ni x\mapsto S_x\in\goth S$ is Borel \cite{AGH}.
It is easy to see that
$$
S_{T_gx}=gS_xg^{-1}\quad\text{ for all }\quad g\in H_3(\Bbb R)\quad\text{ and }x\in X.\tag3-1
$$
We let $\widetilde S_x:=S_x\cap\{c(t)\mid t\in\Bbb R\}$.
It follows from \thetag{3-1} that the Borel map $X\ni x\mapsto \widetilde S_x\in\goth S$ is invariant under $T$ and hence constant almost everywhere.
Since the restriction of $T$ to the center of $H_3(\Bbb R)$ is free, we obtain that
$\widetilde S_x$  is trivial for almost all $x$.
 According to \cite{GoSi} and \cite{CRa}\footnote{The stability groups of an ergodic probability preserving action of a simple connected nilpotent Lie group are conjugate (in a Borel way) over a conull subset.},  there is  $S\in\goth S$ and a Borel map $\phi:X\to G$ such that
$$
\phi(x)S_x\phi(x)^{-1}=S\quad\text{ for almost all $x$}.
\tag3-2
$$
Since $S$ is countable and closed and  $S\cap\{c(t)\mid t\in\Bbb R\}=\{1\}$, it follows that either $S$ is trivial or there are a countable subgroup $D$
 in $\Bbb R$,
 a homomorphism $\varphi:D\to\Bbb R$ and three reals $\theta_1,\theta_2,\theta_3$ such that
$\theta_1^2+\theta_2^2\ne 0$ and\footnote{We note that $D$ is not necessarily closed in $\Bbb R$.
 Consider, for instance, the closed subgroup $S=\{c(n)a(m\tau+n)\mid m,n\in\Bbb Z\}$ in $H_3(\Bbb R)$, where $\tau$ is an irrational number.}
$$
S=\{c(\varphi(d)\theta_3)b(d\theta_2)a(d\theta_1)\mid d\in D\}.
$$
Suppose that $S$ is non-trivial.
We note that \thetag{3-2} defines $\phi$ up to multiplying with an element of the normalizer $N(S)$ of $S$ in $H_3(\Bbb R)$.
We recall that $N(S):=\{g\in H_3(\Bbb R)\mid gSg^{-1}=S\}$.
It is straightforward to verify that
$$
N(S)=\{c(t_1\theta_3)b(t_2\theta_2)a(t_2\theta_1)\mid t_1,t_2\in\Bbb R\}.
$$
Indeed, an element $c(r_3)b(r_2)a(r_1)$ of $H_3(\Bbb R)$ belongs to $N(S)$
if and only if the commutator $[c(r_3)b(r_2)a(r_1),c(\varphi(d)\theta_3)b(d\theta_2)a(d\theta_1)]$ belongs to $S$ for all $d\in D$.
Thus $c(r_1d\theta_2-r_2d\theta_1)\in S$.
This is possible only if $r_1\theta_2-r_2\theta_1=0$, i.e. the vector $(r_1,r_2)\in\Bbb R^2$ belongs to the line generated by $(\theta_1,\theta_2)$, as desired.
We note that $N(S)$ is a closed normal subgroup of $H_3(\Bbb R)$ and $H_3(\Bbb R)/N(S)=\Bbb R$.
Moreover, the natural projection can be written in the following way
$$
H_3(\Bbb R)/N(S)\ni gN(S)=b(-t\theta_1)a(t\theta_2)N(S)\mapsto t\in\Bbb R,
\tag3-3
$$
where the real $t$ is uniquely determined by the coset $gN(S)$.
We now let $\widetilde\phi(x):=\phi(x)N(S)$.
Then $\widetilde\phi$ is a well defined (mod 0) Borel map from $X$ to $H_3(\Bbb R)/N(S)=\Bbb R$.
Substituting \thetag{3-2} into \thetag{3-1} we deduce  that
$$
\phi(T_gx)^{-1}S\phi(T_gx)=g^{-1}\phi(x)^{-1}S\phi(x)g
$$
and hence
$$
\phi(x)g\phi(T_gx)^{-1}\in N(S)\tag3-4
$$
for a.a. $x$.
Therefore
$-\widetilde\phi(x)+\widetilde\phi(T_gx)\in gN(S)$
for a.a. $x$ for each $g\in H_3(\Bbb R)$.
It follows that $\widetilde\phi$ is  invariant under $(T_g)_{g\in N(S)}$.
Denote by $q:(X,\mu)\to (Y,\nu)$ the
$(T_g)_{g\in N(S)}$-ergodic decomposition.
Then there are an ergodic flow $V=(V_t)_{t\in\Bbb R}$ on $(Y,\nu)$ and a function $\psi:Y\to\Bbb R$ such that $V_t\circ q=q\circ T_{b(-t\theta_1)a(t\theta_2)}$, $\widetilde\phi=\psi\circ q$ and $-\psi+\psi\circ V_t=t$ for all $t\in\Bbb R$.
The latter equality follows from \thetag{3-3} and \thetag{3-4}.
Of course, such a function $\psi$ does not exist.
Thus we obtain a contradiction.
 Hence $S$ is trivial, i.e. $T$ is free.

 The converse follows from Proposition~2.6 and from the obvious remark that if $T$ is free then the restriction of $T$ to the center of $H_3(\Bbb R)$ is also free.
\qed
\enddemo

The next claim follows from Proposition~2.6, Theorem~3.1 and Corollary~2.7(ii).

\proclaim{Corollary 3.2}
\roster
\item"\rom{(i)}"
A Heisenberg odometer is free if and only if its restriction to the center is free.
\item"\rom{(ii)}"
A Heisenberg odometer is free if and only if its normal cover is free.
\endroster
\endproclaim

It follows from the next example that the condition $\{c(t)\mid t\in\Bbb R\}\cap\bigcap_{n=1}^\infty\Gamma_n=\{1\}$
does not imply that
$\bigcap_{n=1}^\infty\Gamma_n=\{1\}$.

\example{Example 3.3}
Let $\Gamma_n:=\{c(n!i_3)b(n!i_2)a(i_1)\mid i_1,i_2,i_3\in\Bbb Z\}$.
Then $\Gamma_n$ is a lattice in $H_3(\Bbb R)$, $\Gamma_1\supset\Gamma_2\supset\cdots$ and
$\{c(t)\mid t\in\Bbb R\}\cap\bigcap_{n=1}^\infty\Gamma_n=\{1\}$.
On the other hand,
$\bigcap_{n=1}^\infty\Gamma_n=\{a(i_1)\mid i_1\in\Bbb Z\}$.
\endexample

Consider now the case of normal Heisenberg odometers.

\proclaim{Proposition 3.4}
If $T$ is a normal free Heisenberg odometer then
 $\bigcap_{n=1}^\infty p(\Gamma_n)=\{0\}$.
 \endproclaim
\demo{Proof}
Suppose that $\bigcap_{n=1}^\infty p(\Gamma_n)\ne \{0\}$.
  Then there is $g\in H_3(\Bbb R)\setminus\{c(t)\mid t\in\Bbb R\}$ and a sequence $(r_n)_{n=1}^\infty$, $r_n\in \Bbb R$, such that
  $gc(r_n)\in \Gamma_n$ for all $n\in\Bbb N$.
  Since $g$ does not belong to the center of $H_3(\Bbb R)$,
  there is $\gamma$ in $\Gamma_1$ such that $g$ and $\gamma$ do not commute (see Corollary~1.4).
  Therefore $\gamma g\gamma^{-1}=gc(t)$ for some $t\ne 0$.
  This yields $gc(r_n)c(t)=\gamma gc(r_n)\gamma^{-1}\in \gamma\Gamma_n\gamma^{-1}=\Gamma_n$ and hence $c(t)\in\Gamma_n$ for all $n$.
  This contradicts to the fact that $\bigcap_{n=1}^\infty\Gamma_n=\{1\}$.
  \qed
 \enddemo

It may seem that the condition $\bigcap_{n=1}^\infty p(\Gamma_n)=\{0\}$ is satisfied for every odometer with $\bigcap_{n=1}^\infty\Gamma_n=\{1\}$.
 In the next example we show that this is not true.
Hence the corresponding odometer is not normal.

\example{Example 3.5}
Let $k_1:=1$ and $k_n:=k_{n-1}(k_{n-1}+1)$.
Let
$$
\Gamma_n:=\{c(k_n^2i_3+k_ni_1)b(k_n^2i_2)a(i_1))\mid i_1,i_2,i_3\in\Bbb Z\}.
$$
Then $\Gamma_n$ is a lattice in $H_3(\Bbb R)$,
$\Gamma_1\supset\Gamma_2\supset\cdots$ and $\bigcap_{n=1}^\infty\Gamma_n=\{1\}$.
On the other hand, $p(\Gamma_n)=
\pmatrix
1 & 0  \\
0 & k_n^2
\endpmatrix
\Bbb Z^2$ and hence $\bigcap_{n=1}^\infty p(\Gamma_n)=\pmatrix 1 & 0\\ 0 &0\endpmatrix \Bbb Z^2\ne \{0\}$.
\endexample

\head 4.
Spectral analysis~for~homogeneous~actions~of~$H_3(\Bbb R)$
on nil-manifolds
\endhead

In this section we review
\roster
\item"(a)" Kirillov's orbit method to get a description of the unitary dual  for the Heisenberg group \cite{Ki} and
\item"(b)" the Howe-Richardson formula for the multiplicities of the irreducible components in spectral decomposition of  homogeneous $H_3(\Bbb R)$-spaces \cite{Ho}, \cite{Ri}.
\endroster

Given a lattice $\Gamma$ in $H_3(\Bbb R)$, the quotient space $H_3(\Bbb R)/\Gamma$ is called a {\it nil-manifold}.
Consider the quasi-regular representation $U$ of $H_3(\Bbb R)$ in $L^2(H_3(\Bbb R)/\Gamma)$, i.e. the Koopman representation of $H_3(\Bbb R)$ corresponding to the natural action of this group by rotations on $H_3(\Bbb R)/\Gamma$.
Since $\Gamma$ is uniform, it follows that there is a countable family $\Theta$ of irreducible unitary representations of $H_3(\Bbb R)$ and positive integers $m_W$, $W\in\Theta$,
such that
$$
U=\bigoplus_{W\in\Theta}\bigoplus_{j=1}^{m_{W}} W.\tag4-1
$$
The family  $\Theta$ (or, rather their unitary equivalence classes) and the map $\Theta\ni W\mapsto m_W$
are defined uniquely by $U$.
Our purpose is to describe them explicitly.
Partition $\Theta$ into two subsets $\Theta_1$ and $\Theta_\infty$ consisting of the 1-dimensional representations and the $\infty$-dimensional representations of $H_3(\Bbb R)$ respectively.

If $\pi_{\alpha,\beta}\in\Theta_1$ for some $(\alpha,\beta)\in\Bbb R^2$, then there is $f\in L^2(H_3(\Bbb R)/\Gamma)$ such that $f\circ T_g^{-1}=e^{2\pi i(\alpha t_1+\beta t_2)}f$ for all $g\in\Bbb H_3(\Bbb R)$, where $g=c(t_3)b(t_2)a(t_1)$.
In particular, $f$ is invariant under the action of the center.
We recall (see the footnote No 3) that the natural quotient map
$$
H_3(\Bbb R)/\Gamma\to H_3(\Bbb R)/p^{-1}(p(\Gamma))=\Bbb R^2/p(\Gamma)
$$
 is the $(T_{c(t)})_{t\in\Bbb R}$-ergodic decomposition.
Hence $f$ is constant on the fibers of this quotient map.
Moreover, $f$ is an eigenfunction for the $\Bbb R^2$-action by rotations
on $\Bbb R^2/p(\Gamma)$.
The converse is also true.
Indeed, if $A\in\text{GL}_2(\Bbb R)$ then $(A^*)^{-1}\Bbb Z^2$ are the eigenvalues
of the $\Bbb R^2$-actions by rotations on $\Bbb R^2/A(\Bbb Z^2)$.
If $p(\Gamma)=A(\Bbb Z^2)$ for some matrix $A\in\text{GL}_2(\Bbb R)$ then
we denote by $p(\Gamma)^*$ the {\it dual lattice} $(A^*)^{-1}\Bbb Z^2$ in $\Bbb R^2$.
It is easy to see that the dual lattice is well defined, i.e. does not depend on the choice of $A$.
We now obtain that
$$
\Theta_1=\{\pi_{\alpha,\beta}\mid (\alpha,\beta)\in  p(\Gamma)^*\}.
 $$
Since $T$ is ergodic (it is even transitive), it follows that the corresponding multiplicities are all $1$, i.e. $m_W=1$ whenever $W\in\Theta_1$.

Now let us calculate  $\Theta_\infty$.
Fix a $p(\Gamma)$-fundamental parallelepiped $\Delta$ in $\Bbb R^2$.
Consider the map
$$
\Psi:H_3(\Bbb R)/\Gamma\ni g\Gamma\mapsto ((t_1,t_2)+p(\Gamma),t_3+\xi_\Gamma\Bbb Z)\in
H_3(\Bbb R)/p^{-1}(p(\Gamma))\times\Bbb R/\xi_\Gamma\Bbb Z,
$$
where the reals $t_1,t_2,t_3$ are determined by the formula $g\in c(t_3)b(t_2)a(t_1)\Gamma$ with $t_1,t_2\in\Delta$.
This map conjugates the $H_3(\Bbb R)$-action on the nilmanifold
$X:=H_3(\Bbb R)/\Gamma$
with the  $\Bbb R/\xi_\Gamma\Bbb Z$-group extension of the
 $H_3(\Bbb R)$-action  on  the quotient space $Y:=H_3(\Bbb R)/p^{-1}(p(\Gamma))$.
We let $\eta_n(w)=e^{2\pi i nw/\xi_\Gamma}$, $n\in\Bbb R$, $w\in[0,\xi_\Gamma)=\Bbb R/\xi_\Gamma\Bbb Z$.
Since the group characters $\{\eta_n\mid n\in\Bbb Z\}$ form an orthonormal basis  in $L^2(\Bbb R/\xi_\Gamma\Bbb Z)$, we
obtain a decomposition of $L^2(X)$ into the following orthogonal sum
$$
L^2(X)=\bigoplus_{n\in\Bbb Z} L^2(Y)\otimes\eta_n.\tag4-2
$$
It is easy to see that every subspace $L^2(Y)\otimes\eta_n$
is invariant under $U$.
We note that if $\Psi(g\Gamma)=(y,t_3+\xi_\Gamma\Bbb Z)$ for some $g\in H_3(\Bbb R)$, $y\in Y$ and $t_3\in\Bbb R$
then $\Psi(c(t)g\Gamma)= (y,t+t_3+\xi_\Gamma\Bbb Z)$ for all  $t\in\Bbb R$.
This implies that
 $$
U(c(t))(f\otimes\eta_n)=e^{2\pi i n t/\xi_\Gamma}f\otimes\eta_n
$$
for all $f\in L^2(Y)$, $n\in\Bbb Z$.
Hence $n/\xi_\Gamma\in\Theta_\infty$ and
$$
U\restriction (L^2(Y)\otimes\eta_n)\subset\bigoplus_{1}^{m_{\pi_{n/\xi_\Gamma}}}\pi_{n/\xi_\Gamma} \tag4-3
$$
if $n\ne 0$.
Of course,
 $U\restriction(L^2(Y)\otimes 1)=\bigoplus_{W\in \Theta_1}W$.
 Therefore we deduce from \thetag{4-1}, \thetag{4-2} and \thetag{4-3} that
$$
\Theta_\infty=
\{\pi_{\gamma}\mid \gamma\in  \xi_\Gamma^{-1}\Bbb Z, \gamma\ne 0\}
$$
and
$$
U\restriction (L^2(Y)\otimes\eta_n)=\bigoplus_{1}^{m_{\pi_{n/\xi_\Gamma}}}\pi_{n/\xi_\Gamma}
$$
if $n\ne 0$.
Thus we have computed the spectrum of $U$.

\proclaim{Theorem 4.1} We have $\Theta=
\{\pi_{\alpha,\beta}\mid (\alpha,\beta)\in  p(\Gamma)^*\}\cup
\{\pi_{\gamma}\mid \gamma\in  \xi_\Gamma^{-1}\Bbb Z, \gamma\ne 0\}$.
\endproclaim

Our next aim is to compute  the spectral  multiplicities $m_W$ when $W\in \Theta_\infty$.
For that we will apply Kirillov's orbit method \cite{Ki}
and the Howe-Richardson multiplicity formula (see \cite{Ho}, \cite{Ri} and an earlier work \cite{Mo}).
In view of Proposition~1.3 we may assume
without loss of generality that
$$
\Gamma=\left\{\pmatrix 1&l&\frac n{k}\\
0&1&m\\
0&0&1
\endpmatrix \Bigg|\, l,m,n\in\Bbb Z\right\}
$$
for some $k\in\Bbb N$ which is fixed.
By Theorem~4.1, $\pi_\gamma\in\Theta$ if and only if
$\gamma= kn$ with $n\in\Bbb Z\setminus\{0\}$.
Thus our purpose is to calculate the multiplicity $m_{\pi_{ kn}}$ for each $n\in\Bbb Z\setminus\{0\}$.

The Lie algebra of $H_3(\Bbb R)$ is
$$
\goth h_3(\Bbb R):=\left\{\pmatrix
 0 & a & c\\ 0 & 0 & b\\ 0 & 0 & 0\\
\endpmatrix \Bigg| \,a,b,c\in\Bbb R\right\}.
$$
We endow it with the natural  topology.
Then the exponential map $\exp:\goth h_3(\Bbb R)\to H_3(\Bbb R)$ is a homeomorphism.
We choose a special {\it strong Malcev} basis
$$
X_1:=
\pmatrix 0&1&0\\
0&0&0\\
0&0&0
\endpmatrix,\
X_2:=
\pmatrix 0&0&0\\
0&0&1\\
0&0&0
\endpmatrix, \
X_3:=
\pmatrix 0&0&\frac 1k\\
0&0&0\\
0&0&0
\endpmatrix
$$
for  $\goth h_3(\Bbb R)$.
Fix an ideal $\goth m:=\Bbb R X_3+\Bbb R X_2$ of $\goth h_3(\Bbb R)$.
Then $\exp (\goth m)=H_{2,b}$.
We consider the dual space $\goth h_3(\Bbb R)^*$ to $\goth h_3(\Bbb R)$
as the space of lower triangular $3\times 3$-matrices (modulo the subspace of upper triangular matrices including the diagonal).
The canonical duality form is given via the trace of the product, i.e.
$$
\goth h_3(\Bbb R)\times\goth h_3(\Bbb R)^*\ni
\left(
\pmatrix 0&d_1&d_3\\
0&0&d_2\\
0&0&0
\endpmatrix,
\pmatrix *&*&*\\
z_1&*&*\\
z_3&z_2&*
\endpmatrix
\right)
\mapsto z_1d_1+z_2d_2+z_3d_3.
$$
We note that then
$$
X_1^*=
\pmatrix *&*&*\\
1&*&*\\
0&0&*
\endpmatrix,\
X_2^*=
\pmatrix *&*&*\\
0&*&*\\
0&1&*
\endpmatrix, \
X_3^*=
\pmatrix *&*&*\\
0&*&*\\
k&0&*
\endpmatrix
$$
is the dual basis for $\goth h_3(\Bbb R)^*$.
Then $\goth m^*=\Bbb R X_2^*+\Bbb R X_3^*\subset\goth h_3(\Bbb R)^*$.
It is not difficult now to compute explicitly the co-adjoint action Ad$^*$ of $H_3(\Bbb R)$ on $\goth h_3(\Bbb R)^*$:
$$
\text{Ad}^*(c(t_3)b(t_2)a(t_1))
\pmatrix
*&*&*\\
z_1&*&*\\
z_3&z_2&*
\endpmatrix
=
\pmatrix
*&*&*\\
z_1+t_2z_3&*&*\\
z_3&z_2-t_1z_3&*
\endpmatrix.
$$
Hence the Ad$^*(H_3(\Bbb R))$-orbit $\Omega(l)$ of a linear functional
$l=\pmatrix
*&*&*\\
z_1&*&*\\
z_3&z_2&*
\endpmatrix\in \goth h_3(\Bbb R)^*$
with $z_3\ne 0$
is $l+\Bbb R X_2^*+\Bbb R X_1^*$.
Recall that according to  Kirillov's orbit method \cite{Ki},
$\widehat{H_3(\Bbb R)}$ is in a one-to-one correspondence with Ad$^*(H_3(\Bbb R))$-orbits in $\goth h_3(\Bbb R)^*$.\footnote{The case $z=0$ corresponds to the family of 1-dimensional representations.}
More precisely,
the irreducible unitary representation of $H_3(\Bbb R)$ corresponding to $\Omega(l)$ is induced from the following 1-dimensional representation (character) $\tau$ of $\exp(\goth m)=H_{2,b}$:
$$
\tau(\exp X):=e^{2\pi i\langle X, l\rangle}, \quad X\in\goth m.
$$
It is straightforward to verify that Ind$^{H_3(\Bbb R)}_{H_{2,b}}\tau$
of $H_3(\Bbb R)$ is unitarily equivalent to $\pi_{z_3}$.

We now explain the Howe-Richardson formula for $m_{\pi_{kn}}$.
For that we let
$$
\Cal W:=\{l\in\Bbb ZX_3^*+\Bbb ZX_2^*\mid l\restriction (\Bbb RX_3)\ne 0\}.
$$
To put it another way,
$$
\Cal W=
\left\{
\pmatrix *&*&*\\
0&*&*\\
kn&m&*
\endpmatrix\Bigg|\, n,m\in\Bbb Z, n\ne 0
\right\}.
$$
Given $l\in\goth m^*\subset\goth h_3(\Bbb R)^*$ and $g\in H_3(\Bbb R)$,
we let Ad$^*_0(g)l:=\text{Ad}^*(g)l\restriction\goth m\in\goth m^*$.
We now see that if $\gamma=c(t_3)b(t_2)a(t_1)\in\Gamma$ then
$$
\text{Ad}^*_0(\gamma)
\pmatrix *&*&*\\
0&*&*\\
kn&m&*
\endpmatrix=
\pmatrix *&*&*\\
0&*&*\\
kn&m-t_1nk&*
\endpmatrix.
$$
Since $t_1\in\Bbb Z$, it follows that Ad$^*_0$  is  an action of $\Gamma$ on   $\Cal W$.
By \cite{Ri} and \cite{Ho},
$m_{\pi_{z_3}}$ equals the number of
Ad$^*_0(\Gamma)$-orbits in the intersection $\Omega(l)\cap\Cal W$.
Hence
$$
m_{\pi_{z_3}}=
\cases
|n|k & \text{if }z_3=  kn, \,n\in\Bbb Z\setminus\{0\},\\
0 & \text{otherwise.}
\endcases
$$
We deduce from this formula, Proposition~1.3 and Theorem~4.1
the following theorem.

\proclaim{Theorem 4.2}
    Let $\Gamma$ be an arbitrary lattice in $H_3(\Bbb R)$.
Then
$$
m_{\pi_{\gamma}}=
\cases
|n|k_\Gamma & \text{if }\gamma=  n\xi_\Gamma^{-1}, \ n\in\Bbb Z\setminus\{0\},\\
0 & \text{otherwise.}
\endcases
$$
Therefore
$$
U=\bigoplus_{(\alpha,\beta)\in p(\Gamma)^*} \pi_{\alpha,\beta}\oplus
\bigoplus_{0\ne n\in\Bbb Z}\bigoplus_1^{|n|k_\Gamma}\pi_{n\xi_{\Gamma}^{-1}}.
$$
\endproclaim

In the next statement we show that the homogeneous actions of the Heisenberg group  on nilmanifolds are  isospectral.

\proclaim{Corollary 4.3}
Let $\Gamma$ and $\Gamma'$ be two lattices in $H_3(\Bbb R)$.
Denote by $T$ and $T'$ the corresponding measure preserving actions of $H_3(\Bbb R)$ on the homogeneous spaces $H_3(\Bbb R)/\Gamma$ and $H_3(\Bbb R)/\Gamma'$ respectively.
The following are equivalent:
\roster
\item"\rom{(i)}"
 $T$ and $T'$ are isomorphic.
\item"\rom{(ii)}"
$p(\Gamma)=p(\Gamma')$ and $k_{\Gamma}=k_{\Gamma'}$.
 \item"\rom{(iii)}"
$p(\Gamma)=p(\Gamma')$ and $\xi_{\Gamma}=\xi_{\Gamma'}$.
\item"\rom{(iv)}"
The Koopman representations of $H_3(\Bbb R)$ generated by $T$ and $T'$ are unitarily equivalent.
\item"\rom{(v)}"
$T$ and $T'$ have the same spectrum (i.e. the same  maximal spectral type).
\endroster
\endproclaim

\demo{Proof}
It is well known (and easy to see) that $T$ and $T'$ are isomorphic if and only if  $\Gamma$ and $\Gamma'$ are conjugate.
Therefore (i) and (ii) are equivalent according to the second claim of Proposition~1.3.
It follows from Remark~1.6 that  (ii) is equivalent (iii).
Of course, (i) implies (iv).
It follows from Theorem~4.2 that (iv) implies (iii).
 The equivalence of  (iv) and  (v) follows  also from Theorem~4.2.
\qed
\enddemo

\head 5.  Heisenberg odometers
\endhead

In this section we describe the spectral decomposition of Heisenberg odometers into irreducible components.
A system of algebraic invariants (consisting of two countable Abelian groups with a relation between them) is associated to every $H_3(\Bbb R)$-odometer.
It is shown that this system of invariants is complete for a certain f-isomorphism  relation (to be defined below) on the class of Heisenberg odometers.
We also show that the  f-isomorphism   is weaker then the isomorphism, i.e. the class of $H_3(\Bbb R)$-odometers is not isospectral.

Let $\Gamma_1\supset\Gamma_2\supset\cdots$ be a nested sequence of lattices in $H_3(\Bbb R)$.
Denote by $(X,\mu,T)$ the associated $H_3(\Bbb R)$-odometer.
Let $(Y,\nu)$ stand for the space of $(T_{c(t)})_{t\in\Bbb R}$-ergodic components and let $f:X\to Y$ stand for the corresponding projection.
Then an $\Bbb R^2$-action $V=(V_{t_1,t_2})_{(t_1,t_2)\in\Bbb R^2}$
is well defined by the formula $V_{t_1,t_2}f(x):=f(T_{b(t_2)a(t_1)}x)$.
We call it the {\it underlying $\Bbb R^2$-odometer.}\footnote{It is the $\Bbb R^2$-odometer associated with the sequence $(p(\Gamma_j))_{j=1}^\infty$ of lattices $p(\Gamma_j)$ in $\Bbb R^2$.}

\definition{Definition 5.1}
We say that $T$ is non-degenerate if one of the following equivalent conditions is satisfied:
\roster
\item"(i)"
The underlying $\Bbb R^2$-odometer is non-transitive.
\item"(ii)"
The subgroup $\bigcup_{j=1}^\infty p(\Gamma_j)^*$ is not closed in $\Bbb R^2$.
\item"(iii)"
The sequence $(p(\Gamma_j))_{j=1}^\infty$ does not stabilize, i.e.
for each $j>0$ there is $j_1>j$ such that $p(\Gamma_{j})\ne p(\Gamma_{j_1})$.
\endroster
\enddefinition

\proclaim{Theorem 5.2}
Let $U$ stand for the Koopman unitary representation of $H_3(\Bbb R)$ generated by a  Heisenberg odometer $T$.
\roster
\item"\rom(i)"
If $T$ is non-degenerate then
$$
U=\bigoplus_{(\alpha,\beta)\in \bigcup_{j=1}^\infty p(\Gamma_j)^*} \pi_{\alpha,\beta}\oplus
\bigoplus_{0\ne\gamma\in \bigcup_{j=1}^\infty \xi_{\Gamma_j}^{-1}\Bbb Z}\bigoplus_1^{\infty}\pi_\gamma.
$$
\item"\rom(ii)" If there is $l>0$ such that $p(\Gamma_{j})= p(\Gamma_{l})$ for all $j\ge l$ then
$$
U=\bigoplus_{(\alpha,\beta)\in  p(\Gamma_{l})^*} \pi_{\alpha,\beta}\oplus
\bigoplus_{0\ne\gamma\in \bigcup_{j=l}^\infty \xi_{\Gamma_j}^{-1}\Bbb Z}\bigoplus_1^{m(\gamma)}\pi_\gamma,
$$
where $m(\gamma):=|\gamma|\xi_{\Gamma_j}k_{\Gamma_j}$ for each  $\gamma\in\xi_{\Gamma_j}^{-1}\Bbb Z$, $j\ge l$.
\endroster
\endproclaim
\demo{Proof}
(i) Consider
the following nested family
$$
L^2(H_3(\Bbb R)/\Gamma_1)\subset L^2(H_3(\Bbb R)/\Gamma_2) \subset\cdots
$$
of $U$-invariant subspaces in $L^2(X)$.
Their union $\bigcup_{j=1}^\infty
L^2(H_3(\Bbb R)/\Gamma_j)$ is dense in $L^2(X)$.
Therefore Theorem~4.2 yields that
$$
U=\lim_{j\to\infty}\bigg(\bigoplus_{(\alpha,\beta)\in p(\Gamma_j)^*} \pi_{\alpha,\beta}\oplus
\bigoplus_{0\ne n\in\Bbb Z}
\bigoplus_1^{|n|k_{\Gamma_j}}\pi_{n\xi_{\Gamma_j}^{-1}}\bigg).
$$
Thus we see that $U$ decomposes into a direct sum of 1-dimensional and infinite dimensional irreducible unitary representations.
Of course, the multiplicity of every 1-dimensional irreducible sub-representation $\pi_{\alpha,\beta}$ of $U$ is $1$.
Our purpose now is to compute the multiplicity of each infinite dimensional irreducible sub-representation $\pi_\gamma$, $0\ne\gamma\in \bigcup_{j=1}^\infty \xi_{\Gamma_j}^{-1}\Bbb Z$, of $U$.
There are $j\in\Bbb N$ and $n\in\Bbb Z\setminus\{0\}$ such that
$$
\gamma=\frac{n
}{\xi_{\Gamma_j}}=\frac{n
q
}{\xi_{\Gamma_{j+1}}}, \text{ where }q:=\frac{\xi_{\Gamma_{j+1}}}{\xi_{\Gamma_{j}}}.
$$
Then the multiplicity of $\pi_{\gamma}$
in $U\restriction L^2(H_3(\Bbb R)/\Gamma_{j})$ is $|n|k_{\Gamma_j}$ and the multiplicity of  $\pi_{\gamma}$
in $U\restriction L^2(H_3(\Bbb R)/\Gamma_{j+1})$ is $|n|qk_{\Gamma_{j+1}}$.
Let $p(\Gamma_j)=A_j(\Bbb Z^2)$ for some matrix $A_j\in GL_2(\Bbb R)$, $j\in\Bbb N$.
It follows from Remark 1.6 that
$$
\frac{k_{\Gamma_{j+1}}}{k_{\Gamma_j}}
=
\frac{|\det A_{j+1}|}{|\det A_j|}
\cdot\frac{\xi_{\Gamma_{j}}}{\xi_{\Gamma_{j+1}}}
.
$$
Then
$$
|n|qk_{\Gamma_{j+1}}=|n|\frac{\xi_{\Gamma_{j+1}}}{\xi_{\Gamma_j}}
k_{\Gamma_j}
\frac{|\det A_{j+1}|}{|\det A_j|}
\cdot\frac{\xi_{\Gamma_{j}}}{\xi_{\Gamma_{j+1}}}=
|n|
k_{\Gamma_j}
\frac{|\det A_{j+1}|}{|\det A_j|}.
$$
Since $T$ is non-degenerated, we may assume without loss of generality that $p(\Gamma_j)\ne p(\Gamma_{j+1})$ and hence $\frac{|\det A_{j+1}|}{|\det A_j|}\ge 2$ for each $j$.
Hence by Theorem~4.2,
the multiplicity of $\pi_{\gamma}$
in $U\restriction L^2(H_3(\Bbb R)/\Gamma_{j+1})$ is at least two times greater
than the multiplicity of  $\pi_{\gamma}$
in $U\restriction L^2(H_3(\Bbb R)/\Gamma_{j})$.
It follows that the multiplicity of $\pi_\gamma$ in $U$ is infinite.

(ii) is proved in a similar way.
It is enough to note that
$$
\frac{k_{\Gamma_{j+1}}}{k_{\Gamma_j}}=
\frac{\xi_{\Gamma_{j}}}{\xi_{\Gamma_{j+1}}}
$$
for all $j\ge l$ and hence $m(\gamma)$ is well defined.
\qed
\enddemo

\definition{Definition 5.3}
 A subgroup $S$ in $\Bbb R^m$ is {\it off-rational} if its closure $\overline S$ is co-compact in $\Bbb R^m$ and there are a subgroup $Q\subset\Bbb Q^m$ and a matrix $A\in GL_m(\Bbb R)$  such that
$S=AQ$.
\enddefinition

We now associate to $S$ an off-rational  subgroup $\tau(S)$ in $\Bbb R$.
Since $S$ is off-rational, there is a  sequence of matrices  $A_j\in GL_m(\Bbb  R)\cap M_m(\Bbb Z)$ such that
$A_1^{-1}\Bbb Z^m\subset A_2^{-1}\Bbb Z^m\subset\cdots$ and $\bigcup_{j=1}^\infty A_j^{-1}\Bbb Z^m=Q$ and hence
$S=\bigcup_{j=1}^\infty AA_j^{-1}\Bbb Z^m$.\footnote{We note that $A_j$ and $A$ are defined up to right multiplication with a matrix from $GL_m(\Bbb Z)$.}
Consider now a sequence of subgroups
$$
\frac{\det A}{\det A_1}\Bbb Z\subset \frac{\det A}{\det A_2} \Bbb Z\subset\cdots
$$
in $\Bbb R$.
Then  $\tau(S):=\bigcup_{j=1}^\infty \frac{\det A}{\det A_j}\Bbb Z$
 is a dense off-rational subgroup  of $\Bbb R$.
It is straightforward to verify that $\tau(S)$ does not depend on the choice of  the sequence $(A_j)_{j=1}^\infty$.

Suppose we are given a  sequence $\Gamma=(\Gamma_j)_{j=1}^\infty$ of lattices
$\Gamma_1\supset\Gamma_2\supset\cdots$ in $H_3(\Bbb R)$.
Then $S_{\Gamma}:=\bigcup_{j=1}^\infty p(\Gamma_j)^*$
is an off-rational subgroup of $\Bbb R^2$
and $\xi_{\Gamma}:=\bigcup_{j=1}^\infty\xi_{\Gamma_j}^{-1}\Bbb Z$ is an off-rational subgroup in $\Bbb R$.
If $\{c(t)\mid t\in\Bbb R\}\cap\bigcap_{j=1}^\infty\Gamma_j=\{1\}$ then $\xi_{\Gamma}$ is dense in $\Bbb R$.
Thus, in view of Theorem~3.1, if $T$ is free then $\xi_\Gamma$ is dense in $\Bbb R$.

\proclaim{Proposition 5.4} We have $\tau( S_\Gamma)\supset\xi_\Gamma$.
\endproclaim
\demo{Proof}
Let $p(\Gamma_j)=A_j(\Bbb Z^2)$ for some $A_j\in GL_2(\Bbb R)$, $j\in\Bbb N$.
Since $k_{\Gamma_j}=\frac{|\det A_{j}|
}{\xi_{\Gamma_j}}$ by Remark~1.6, we obtain that
$\frac 1{\det A_{j}}\Bbb Z\supset\frac 1{\xi_{\Gamma_j}}\Bbb Z$ for each $j$.
\qed
\enddemo

\proclaim{Proposition 5.5}
Given an off-rational subgroup $S$ in $\Bbb R^2$
and an off-rational subgroup $\xi$ in $\Bbb R$ such that
$\tau(S)\supset\xi$, there is a sequence $\Gamma$ of lattices $\Gamma_1\supset\Gamma_2\supset\cdots$ in $H_3(\Bbb Z)$
such that
$S_\Gamma=S$ and $\xi_\Gamma=\xi$.
If $S$ is dense then
$\bigcap_{j=1}^\infty p(\Gamma_j)=\{0\}$.
If, in addition, $\xi$ is dense in $\Bbb R$ then
$\bigcap_{j=1}^\infty\Gamma_j=\{1\}$.
\endproclaim
\demo{Proof}
Let us choose for each $j>0$, a  matrix $A_j\in GL_2(\Bbb R)$ with $\det A_j>0$ and a  positive real $\zeta_j$ such that
$(A_1^*)^{-1}\Bbb Z^2\subset (A_2^*)^{-1}\Bbb Z^2\subset\cdots$,
$\bigcup_{j=1}^\infty(A_j^*)^{-1}\Bbb Z^2=S$,
 $\zeta_1^{-1}\Bbb Z\subset\zeta_2^{-1}\Bbb Z\subset\cdots$ and
$\bigcup_{j=1}^\infty\zeta_j^{-1}\Bbb Z=\xi$.
We then have
$A_1\Bbb Z^2\supset A_2\Bbb Z^2\supset\cdots$
and
$\zeta_1\Bbb Z\supset\zeta_2\Bbb Z\supset\cdots$
Since $\tau(S)\supset\xi$, we can assume without loss of generality (passing to a subsequence of $(A_j)_{j=1}^\infty$ if necessary) that $\det (A_j^*) ^{-1}| \zeta_j^{-1}$ or, equivalently,
$\zeta_j | \det A_j$ for each $j$.
We now let
$$
\Gamma_j:=\{c(t_3)b(t_2)a(t_1)\mid (t_1,t_2)\in A_j(\Bbb Z^2), t_3\in \zeta_j\Bbb Z\}.
$$
Then $\Gamma_j$ is a lattice in $H_3(\Bbb R)$ with $\xi_{\Gamma_j}=\zeta_j$,
$p(\Gamma_j)=A_j(\Bbb Z^2)$ and $k_{\Gamma_j}=\det A_j/\zeta_j$.
We let $\Gamma:=(\Gamma_j)_{j=1}^\infty$.
Of course,  $\Gamma_1\supset\Gamma_2\supset\cdots$, $\xi_\Gamma=\xi$ and $S_\Gamma=S$.

If $S$ is dense in $\Bbb R^2$, it follows that
$\bigcap_{j=1}^\infty p(\Gamma_j)=\{0\}$.
If $\xi$ is dense in $\Bbb R$, it follows that
$\bigcap_{j=1}^\infty \zeta_j\Bbb Z=\{0\}$ and
therefore $\{c(t)\mid t\in\Bbb R\}\cap\bigcap_{j=1}^\infty\Gamma_j=\{1\}$.
 It follows that
$\bigcap_{j=1}^\infty\Gamma_j=\{1\}$.
\qed
\enddemo

\definition{Definition 5.6} Two $H_3(\Bbb R)$-odometers $T$ and $T'$ are called {\it f-isomorphic} if
they are associated with some sequences $(\Gamma_j)_{j=1}^\infty$ and $(\Gamma_j')_{j=1}^\infty$ (respectively) of lattices in $H_3(\Bbb R)$ such that $\Gamma_j$ and $\Gamma_j'$ are conjugate in $H_3(\Bbb R)$ for each $j$.
\enddefinition

\proclaim{Theorem 5.7} Let $\Gamma=(\Gamma_j)_{j=1}^\infty$ and $\Gamma'=(\Gamma'_j)_{j=1}^\infty$ be two  sequences of lattices in $H_3(\Bbb R)$ such that $\Gamma_1\supset\Gamma_2\supset\cdots$ and $\Gamma'_1\supset\Gamma'_2\supset\cdots$.
Let $T$ denote the odometer associated to $\Gamma$ and let $T'$ denote the odometer associated to $\Gamma'$.
Then $T$ and $T'$ are f-isomorphic if  and only if $S_\Gamma=S_{\Gamma'}$
and $\xi_\Gamma=\xi_{\Gamma'}$.
\endproclaim

\demo{Proof}
Of course, if $T$ and $T'$ are f-equivalent then $S_\Gamma=S_{\Gamma'}$
and $\xi_\Gamma=\xi_{\Gamma'}$.
Let us prove the converse.
We first set $j_1:=1$.
Since  $S_\Gamma=S_{\Gamma'}$ and $\xi_\Gamma=\xi_{\Gamma'}$,
there are $l_1>1$  and $j_2>j_1$ such that $p(\Gamma_{j_1})^*\subset p(\Gamma'_{l_1})^*\subset p(\Gamma_{j_2})^*$
and $\xi_{\Gamma_{j_1}}^{-1}\Bbb Z \subset\xi_{\Gamma'_{l_1}}^{-1}\Bbb Z\subset \xi_{\Gamma_{j_2}}^{-1}\Bbb Z$.\footnote{This is true because all the groups $p(\Gamma_j)^*$, $p(\Gamma_j')^*$, $\xi_{\Gamma_j}^{-1}\Bbb Z$, $\xi_{\Gamma_j'}^{-1}\Bbb Z$, $j\in\Bbb N$, are finitely generated.}
Hence $p(\Gamma_{j_1})\supset p(\Gamma'_{l_1})\supset p(\Gamma_{j_2})$, $\xi_{\Gamma_{j_1}}|\xi_{\Gamma'_{l_1}}$ and $\xi_{\Gamma'_{l_1}}|\xi_{\Gamma_{j_2}}$.
We now let $\widetilde\Gamma_{j_1}:=\Gamma_{j_1}\cap p^{-1}(p(\Gamma'_{l_1}))$.
Then $\widetilde\Gamma_{j_1}$ is a lattice in $H_3(\Bbb R)$,\footnote{$\Gamma_{j_1}$ and $\widetilde \Gamma'_{j_1}$ have the same $p$-fibers over $p(\Gamma'_{l_1})$.}
$\widetilde\Gamma_{j_1}\supset\Gamma_{j_2}$,
$p(\widetilde\Gamma_{j_1})=p(\Gamma'_{l_1})$ and $\xi_{\widetilde\Gamma_{j_1}}=\xi_{\Gamma_{j_1}}$.
 According to Lemma~1.5 there is a lattice $\widehat\Gamma_{j_1}$ such that $\widetilde\Gamma_{j_1}\supset\widehat\Gamma_{j_1}\supset\Gamma_{j_2}$ with
$p(\widehat\Gamma_{j_1})=p(\Gamma'_{l_1})$ and $\xi_{\widehat\Gamma_{j_1}}=\xi_{\Gamma'_{l_1}}$.
It follows from Proposition~1.3 and Remark~1.6 that
$\widehat\Gamma_{j_1}=q_1\Gamma'_{l_1}q_1^{-1}$ for some element $q_1\in H_3(\Bbb R)$.
In the next step we argue in a similar way to find integers $j_3>j_2$ and $l_2>l_1$, a lattice $\widehat\Gamma_{j_2}\subset H_3(\Bbb R)$ and an element $q_2\in H_3(\Bbb R)$ such that
$\Gamma_{j_2}\supset\widehat\Gamma_{j_2}\supset\Gamma_{j_3}$
and $\widehat\Gamma_{j_2}=q_2\Gamma'_{l_2}q_2^{-1}$.
 Continuing this infinitely many times we obtain a sequence
$(\widehat\Gamma_{j_s})_{s=1}^\infty$ of lattices in $H_3(\Bbb R)$,
 a subsequence $(\Gamma'_{l_s})_{s=1}^\infty$ of $\Gamma'$
and a sequence $(q_s)_{s=1}^\infty$ of elements in $H_3(\Bbb R)$ such that
$$
\align
&\Gamma_{j_1}\supset\widehat\Gamma_{j_1}
\supset\Gamma_{j_2}\supset\widehat\Gamma_{j_2}\supset\cdots
\text{ and }\tag5-1\\
&\widehat\Gamma_{j_s}=q_s\Gamma'_{l_s}q_s^{-1}
\text{ with $q_s\in H_3(\Bbb R)$ for each $s>1$}.
\endalign
$$
It follows from \thetag{5-1} that $T$ is isomorphic to the $H_3(\Bbb R)$-odometer associated with $(\widehat\Gamma_{j_s})_{s=1}^\infty$.
Of course, $T'$ is isomorphic  to the $H_3(\Bbb R)$-odometer associated with $(\Gamma'_{l_s})_{s=1}^\infty$.
Hence $T$ and $T'$ are f-isomorphic.
\qed
\enddemo

\comment
We now define a mapping $d_s:H_3(\Bbb R)/\Gamma'_{l_s}\to H_3(\Bbb R)/\widehat\Gamma_{j_s}$ by setting $d_s(g\Gamma'_{l_s}):=gq_s^{-1}\widehat\Gamma_{j_s}$.
It is easy to see that $d_s$ is a well defined  $H_3(\Bbb R)$-equivariant isomorphism.
It follows from \thetag{5-2} that the following infinite diagram
$$
\CD
H_3(\Bbb R)/\Gamma'_{l_1} @<<<  H_3(\Bbb R)/\Gamma'_{l_2} @<<< \cdots\\
@V{d_1}VV            @VV{d_2}V\\
H_3(\Bbb R)/\widehat\Gamma_{j_1} @<<<  H_3(\Bbb R)/\widehat\Gamma_{j_2} @<<< \cdots\\
\endCD
$$
commutes.
Hence $T$ and $T'$ are isomorphic.
\endcomment

\proclaim{Corollary 5.8} The Heisenberg odometers
 $T$ and $T'$ are f-isomorphic if and only if the Koopman unitary representations of $H_3(\Bbb R)$ associated with them are unitarily equivalent.
\endproclaim

\demo{Proof}
Let $T$ be associated with a sequence $\Gamma$ of nested lattices in $H_3(\Bbb R)$ and let $T'$ be associated with a sequence $\Gamma'$ of nested lattices in $H_3(\Bbb R)$.
It follows from Theorem~5.2 that if the unitary representations of $H_3(\Bbb R)$ associated with $T$ and $T'$ are unitarily equivalent
then $S_\Gamma=S_{\Gamma'}$ and $\xi_{\Gamma}=\xi_{\Gamma'}$.
Applying Theorem~5.7 we conclude that $T$ and $T'$ are f-isomorphic.
The converse is obvious.
\qed
\enddemo

We now consider Cartesian products of Heisenberg odometers.
First of all we  recall the Clebsch-Gordan decompositions for the tensor products of two irreducible unitary representations of $H_3(\Bbb R)$ (see \cite{Ki}):
$$
\align
\pi_{\alpha,\beta}\otimes\pi_{\alpha',\beta'} &=\pi_{\alpha+\alpha',\beta+\beta'},\quad\text{for all }\alpha,\beta,\alpha',\beta',\\
\pi_\gamma\otimes\pi_{\gamma'}&=\bigoplus_{1}^\infty\pi_{\gamma+\gamma'}\quad\text{whenever }\gamma\ne -\gamma',\\
\pi_\gamma\otimes\pi_{-\gamma}&=\int_{\Bbb R^2}\pi_{\alpha,\beta}\,d\alpha\, d\beta\quad\text{for all }  \gamma\\
\pi_\gamma\otimes\pi_{\alpha,\beta}&=\pi_\gamma\quad\text{for all }\alpha,\beta,\gamma.
\endalign
$$

\proclaim{Proposition 5.9}
Let $T$ and $T'$ be two  Heisenberg odometers associated with the nested sequences of lattices $\Gamma=(\Gamma_j)_{j=1}^\infty$ and $\Gamma'=(\Gamma'_j)_{j=1}^\infty$ in $H_3(\Bbb R)$ respectively.
Then
\roster
\item"\text{(i)}"
$T\times T'$ is ergodic  if and only if  $S_\Gamma\cap S_{\Gamma'}=\{0\}$.
\item"\text{(ii)}"
$T\times T'$ is ergodic and has discrete maximal spectral type if and only if  $S_\Gamma\cap S_{\Gamma'}=\{0\}$ and $\xi_{\Gamma}\cap\xi_{\Gamma'}=\{0\}$.
In this case the Koopman unitary representation $U_{T\times T'}$ of $H_3(\Bbb R)$ decomposes into irreducible representations as follows
$$
U_{T\times T'}=\bigoplus_{(\alpha,\beta)\in S_\Gamma+S_{\Gamma'}}\pi_{\alpha,\beta}
\oplus\bigoplus_{0\ne \gamma\in\xi_\Gamma+\xi_{\Gamma'}}\bigoplus_1^\infty\pi_{\gamma}.\tag5-2
$$
\item"\text{(iii)}" $T\times T'$ is not isomorphic to any  Heisenberg odometer.
\footnote{As follows from the proof below, $T\times T'$ even is not spectrally equivalent to any Heisenberg odometer.}
\endroster
\endproclaim

\demo{Proof} (i) and (ii) follow immediately from the Clebsch-Gordan decompositions, Theorem~4.1 and the fact that the  invariant vectors of the Koopman unitary representation of $H_3(\Bbb R)$ form the subspace corresponding to the component $\pi_{0,0}$.

 (iii) Suppose that $T\times T'$ is a Heisenberg odometer.
Then it is ergodic and has  discrete maximal spectral type.
Hence by (ii),  $\xi_\Gamma\cap\xi_{\Gamma'}=\{0\}$.
This implies that the sum $\xi_\Gamma+\xi_{\Gamma'}$ is not off-rational.\footnote{Indeed, if $\theta_1 Q_1\cap\theta_2 Q_2=0$ for some non-zero reals $\theta_1,\theta_2$ and nontrivial subgroups $Q_1,Q_2\subset\Bbb Q$ then  $\theta_1/\theta_2\not\in\Bbb Q$. If, furthermore, $\theta_1Q_1+\theta_2Q_2=\theta_3Q_3$ for a non-zero $\theta_3\in\Bbb R$ and nontrivial subgroup $Q_3\subset \Bbb Q$ then we get that $\theta_3/\theta_2$ is rational and irrational at the same time.}
On the other hand,  by \thetag{5-2} and Theorem~4.1, $\xi_\Gamma+\xi_{\Gamma'}=\xi_{\Gamma''}$ for a nested sequence $\Gamma''$ of lattices in $H_3(\Bbb R)$ that determines the Heisenberg odometer $T\times T'$.
Since $\xi_{\Gamma''}$ is always off-rational, we get a contradiction.
\qed
\enddemo

Our next purpose is to show that Heisenberg odometers (even the normal ones) are not isospectral.  

Let $T$ be the Heisenberg odometer associated with a decreasing sequence $(\Gamma_j)_{j+1}^\infty$ of lattices in $H_3(\Bbb Z)$.
Let $(X,\mu)$ be the space of $T$ and fix $x\in X$.
Then $x=(g_j\Gamma_j)_{j\in\Bbb Z}$ with $g_j\Gamma_j=g_{j+1}\Gamma_j$ for all $j\in\Bbb N$.
Of course, we have $g_1\Gamma_1g_1^{-1}\supset g_2\Gamma_2{g_2^{-1}}\supset\cdots$.
Let $T^{(x)}$ denote
 the odometer associated with the sequence  $(g_j\Gamma_jg_j^{-1})_{j=1}^\infty$ of lattices in $H_3(\Bbb R)$.
 Denote by $X^{(x)}$ the space of this odometer.
 Then the map
 $$
 i^{(x)}:X\ni (h_1\Gamma_1,h_2\Gamma_2,\dots)\mapsto(h_1g_1^{-1}\cdot g_1\Gamma_1g_1^{-1},h_2g_2^{-1}\cdot g_2\Gamma_2g_2^{-1},\dots)\in X^{(x)}
 $$
 is a topological conjugacy of $T$ with $T^{(x)}$ such that $i^{(x)}x=(g_j\Gamma_jg_j^{-1})_{j\in\Bbb N}\in X^{(x)}$.
 Below we call $i^{(x)}$ a {\it changing the origin in $X$}.

 \proclaim{Lemma 5.10}
 Let $T$ and $T'$ be $H_3(\Bbb R)$-odometers  associated with the sequences $(\Gamma_j)_{j\in\Bbb N}$ and $(\Gamma_j)_{j\in\Bbb N}$ respectively.
 If $T'$ is a factor (in the measure theoretic sense) of $T$ then there are  sequences $(g_j)_{j\in\Bbb N}$ and $(g_j')_{j\in\Bbb N}$ in $H_3(\Bbb R)$ with $g_j\Gamma_j=g_{j+1}\Gamma_{j}$ and $g_j'\Gamma_j'=g_{j+1}'\Gamma_{j}'$ for all $j\in\Bbb N$ such that for each $n>0$, there is $m>0$ with  $g_m\Gamma_mg_m^{-1}\subset g_n'\Gamma_n' (g_n')^{-1}$.
 \endproclaim
 
 \demo{Proof} Let $(X,\mu)$ and $(X',\mu')$ be the spaces of $T$ and $T'$.
 Let $\phi:X\to X'$ be the Borel factor map intertwining $T$ with $T'$.
 Denote by
 $p_n$ the canonical projection of $X$ onto $H_3(\Bbb R)/\Gamma_n$
  and denote by $q_n$ the canonical projection of $X'$ onto $H_3(\Bbb R)/\Gamma_n'$.
 By Lusin's theorem, for  each $n\in \Bbb N$, there is a closed subset $X_n\subset X$ with $\mu(X_n)>1-2^{-25n}$
 such that the map $q_n\circ\phi$ is continuous on $X_n$.
 For each $n$, let us fix three elements $\gamma_{n,1},\gamma_{n,2}$ and $\gamma_{n,3}$ that generate $\Gamma_n$.
 It follows from the pointwise ergodic theorem that for each $n\in\Bbb N$ and $i\in\{1,2,3\}$ there is a subset $Y_n\subset X_n$ of ``large'' measure such that for each $x\in Y_n$, there exists  a limit
 $$
 \lim_{N\to\infty}\frac{\#\{j\in\{1,\dots,N\}\mid \gamma^j_{n,i}x\in X_n\}}{N}>0.9.
 $$
 The ``large'' means that the series $\sum_{n=1}^\infty (1-\mu(Y_n))$
converges. 
Hence by the Borel-Cantelli lemma, the intersection $\bigcap_{n>M}(X_n\cap Y_n)$ is non-empty
for some $M\in\Bbb N$.
Without loss of generality we may assume that $M=1$.
Changing the origin in $X$, i.e. replacing $T$ with an isomorphic odometer, we may assume that 
$$
e:=(\Gamma_1,\Gamma_2,\dots)\in \bigcap_{n>M}(X_n\cap Y_n).
$$
In a similar way, changing the origin in $X'$ we may assume that 
$$
\phi(e)=e':=(\Gamma_1',\Gamma_2',\dots).
$$
Under these assumptions we will show
 that
 for each $n>0$ there is $m>0$ with  $\Gamma_m\subset \Gamma_n'$.
 
 Given $n$, we fix a small open neighborhood $\widetilde U$ of $\Gamma_n'$ 
 in $H_3(\Bbb Z)/\Gamma_n'$. 
 Then $U:=q_n^{-1}(\widetilde U)$ is an open neighborhood of 
 $e'$ in $X'$.
 Therefore there is an open neighborhood $V$ of $e$ in $X$
such that  
$$
\phi^{-1}(U)\supset V\cap X_n\ne\emptyset.\tag5-3
$$
By definition of the topology on $X$,  we may assume without loss of generality that there are $m>0$ and an open neighborhood $\widetilde V$ of $\Gamma_m$ in $H_3(\Bbb R)/\Gamma_m$ such that $V=p_m^{-1}(\widetilde V)$.
For each $j\in\Bbb N$ and $i\in\{1,2,3\}$, we have $\gamma_{m,i}^j\Gamma_m=\Gamma_m$.
Hence $\gamma^j_{m,i}e\in V$.
Let $J_N:=\{j\in\{1,\dots,N\}\mid \gamma_{m,i}^je\in X_n\}$.
Then for each $j\in J_N$, we have
$ \gamma_{m,i}^je\in V\cap X_n$.
It follows now from \thetag{5-3} and the fact that $\phi$ is equivariant that 
$$
\gamma_{m,i}^je'=\gamma_{m,i}^j\phi(e)=\phi(\gamma_{m,i}^je)\in U.
$$
Hence $\gamma_{m,i}^j\Gamma_n'\in\widetilde U$ for each $j\in J_N$.
Since $\#(J_N)/N>0.9$ for all  large $N$ and $\widetilde U$ is small, it is easy to deduce that
$\gamma_{m,i}\in\Gamma_n'$.
Indeed, $d:=p(\gamma_{m,i})p(\Gamma_n')\in \Bbb R^2/p(\Gamma_n')$.
As usual, $p$ stands for the standard projection of $H_3(\Bbb R)$ onto $\Bbb R^2$.
Take a weak limit point  $\kappa$ of the sequence 
$\frac 1N\sum_{j=1}^N\delta_{\gamma_{m,i}^j\Gamma_n'}$
of measures on $H_3(\Bbb R)/\Gamma_n'$.
The natural projection
 $$
 v:H_3(\Bbb R)/\Gamma_n'\to H_3(\Bbb R)/p^{-1}(p(\Gamma_n'))=\Bbb T^2
 $$ 
 maps $\kappa$ to the limit of the sequence  $\frac 1N\sum_{j=1}^N\delta_{d^j}$.
This limit is the Haar measure, say $\lambda_A$, of the compact monothetic subgroup $A$ of $\Bbb T^2$ generated by the element $d$.
Since we assumed that $\widetilde U$ is small, the projection $v(\widetilde U)$ is also a small (in diameter) neighborhood of the unit in $\Bbb T^2$.
We thus obtain that  $\lambda_A(v(\widetilde U))>0.8$.  (Of course, we may assume without loss of generality that the boundary of $v(\widetilde U)$ is $\lambda_A$-null.) This is only possible if $A$ is trivial, i.e.
$\gamma_{m,i}\in c(t)\Gamma_n'$ for some $t\in [0,\xi_{\Gamma_n'})$.
Thus $\kappa$ is supported by the circle $\Bbb T=v^{-1}(p^{-1}(p(\Gamma_n')))$, which is the fiber of the fibration $v$ over the unit of $\Bbb T^2$.
Since the intersection of $\widetilde U$ with this fiber is small in diameter, we apply a similar reasoning as above to conclude that $t=0$.

Therefore $\Gamma_m\subset\Gamma_n'$.
\qed
 \enddemo

We now provide an example of f-isomorphic but non-isomorphic odometers.

\example{Example  5.11(\rom{cf. \cite{Li--Ug, Example 4.9}})}
Fix a sequence of natural numbers $k_1<k_2<\cdots$ such that $k_1=1$ and $k_n(k_n+1)=k_{n+1}$ for each $n$.
Let 
$$
\align
\Gamma_n
&:=\{c(k_nj_3)b(k_nj_2)a(k_nj_1)\mid j_1,j_2,j_3\in\Bbb Z\}
\text{ and }\\
\Gamma_n'
&:=\{c(k_nj_3+j_1)b(k_nj_2)a(k_nj_1)\mid j_1,j_2,j_3\in\Bbb Z\},
\endalign
$$
 $n\in\Bbb N$.
It is easy to see that
$H_3(\Bbb Z)=\Gamma_1\supset\Gamma_2\supset\cdots$
and
$H_3(\Bbb Z)=\Gamma_1'\supset\Gamma_2'\supset\cdots$.
Moreover, $\Gamma_n$ and $\Gamma_n'$ are both normal in $H_3(\Bbb Z)$ for each $n\in\Bbb N$.
Denote by $T$ and $T'$ the $H_3(\Bbb R)$-odometers associated with the sequences
$(\Gamma_n)_{n=1}^\infty$
and $(\Gamma_n')_{n=1}^\infty$  respectively.
They are normal.
Since
$b(1/k_n)\Gamma_n'b(-1/k_n)=\Gamma_n$ for each $n\in\Bbb N$, we obtain that
$T$ and $T'$ are f-isomorphic.
If they were isomorphic then by Lemma~5.10,  there is $h\in H_3(R)$ such that for each $n>0$, there exists $m>0$ with
$\Gamma_m'\subset h\Gamma_nh^{-1}$.
However
$\xi_{\Gamma_m'}=1$ for each $m$ but $\xi_{h\Gamma_nh^{-1}}\to\infty$ as $n\to\infty$, a contradiction.
We also examine certain ``symmetry'' property   for  $T$ and $T'$.
Recall that  a measure preserving action $R$  of
$H_3(\Bbb R)$
is called {\it symmetric \cite{Da4}} if $R$ is isomorphic to $R\circ\theta$, where $\theta$ denotes the flip in $H_3(\Bbb R)$.
Since $\theta(H_3(\Bbb Z))=H_3(\Bbb Z)$, the symmetric $H_3(\Bbb Z)$-actions are defined
in a similar way.
It is easy to see that  $T$ is symmetric  because $\theta(\Gamma_n)=\Gamma_n$ for each $j$.
It is straightforward to verify that
$$
\theta(\Gamma_n')=a(-1/k_n)\Gamma_na(1/k_n)=
\{c(k_nj_3+j_2)b(k_nj_2)a(k_nj_1)\mid j_1,j_2,j_3\in\Bbb Z\}
$$ 
for each $n$.
We now show that $T'$ is asymmetric.
If $T'$ and $T'\circ\theta$ were isomorphic then by Lemma~5.10, 
there is $h\in H_3(R)$ such that for each $n>0$, there exists $m>0$ with
$\Gamma_m'\subset h\theta(\Gamma_n')h^{-1}$.
However it is easy to see  that 
$$
\inf\{t\in\Bbb R\mid c(t)a(r)\in\Gamma_m'\text{ for some }r\in\Bbb R'\}=1
$$
while 
$$
\inf\{t\in\Bbb R\mid c(t)a(r)\in\theta(\Gamma_n')\text{ for some }r\in\Bbb R'\}\to\infty\quad\text{ as $n\to\infty$},
$$
a contradiction.
Hence $T'$ is not symmetric.
We also note that $T'$
 is  f-isomorphic to $T'\circ\theta$ because
 $\theta(\Gamma_n')=a(-1/k_n)b(1/k_n)\Gamma_n'b(-1/k_n)a(1/k_n)$.
  Hence the Koopman representations $U_{T'}$ and $U_{T'}\circ\theta$ are unitarily equivalent.
\endexample

We conclude this section with a ``coordinate-free'' definition of  Heisenberg odometers.

\proclaim{Proposition 5.12} Let $T$ be an ergodic action of $H_3(\Bbb R)$ on a standard probability space $(X,\goth B,\mu)$.
Let $\goth F_1\subset\goth F_2\subset\cdots$ be a sequence of factors of $T$ such that the union
$\bigcup_{j=1}^\infty\goth F_j$ is dense in $\goth B$ and the restriction $T\restriction\goth F_j$ is isomorphic to a transitive $H_3(\Bbb R)$-action on the 3-dimensional torus.
Then $T$ is a Heisenberg odometer.
\endproclaim
\demo{Proof}
Since the dynamical system $T\restriction\goth F_j$ is transitive, it is isomorphic to the homogeneous space $H_3(\Bbb R)/\Gamma_j$ for a  co-compact subgroup $\Gamma_j\subset H_3(\Bbb R)$.
Since $H_3(\Bbb R)/\Gamma_j$ is isomorphic to the 3-dimensional torus, $\Gamma_j$ is a lattice.
Thus the dynamical system $(X,\mu,T)$ is the projective limit of the following sequence of $H_3(\Bbb R)$-homogeneous spaces
$$
H_3(\Bbb R)/\Gamma_1 @< q_1<< H_3(\Bbb R)/\Gamma_2 @< q_2<<\cdots.
$$
Hence there is a sequence $(g_n)_{n=1}^\infty$ of elements in $H_3(\Bbb R)$
such that $\Gamma_j\supset g_j^{-1}\Gamma_{j+1}g_j$
 and $q_j(g\Gamma_{j+1})=gg_j \Gamma_j$ for each $j$.
 We now let
 $$
 \Gamma_j':=g_1^{-1}g_2^{-1}\cdots g_{j-1}^{-1}\Gamma_j g_{j-1}\cdots g_2g_1
 $$
 and define  an equivariant  isomorphism  $r_j:H_3(\Bbb R)/\Gamma_j\to H_3(\Bbb R)/\Gamma_j'$
by setting
$r_j(g\Gamma_j):=gg_{j-1}\cdots g_2g_1\Gamma_j'$.
Then we obtain a nested sequence
$\Gamma_1'\supset\Gamma_2'\supset\cdots$ of lattices in $H_3(\Bbb R)$.
Moreover,  the following infinite diagram commutes
$$
\CD
H_3(\Bbb R)/\Gamma_1 @< q_1<< H_3(\Bbb R)/\Gamma_2 @< q_2<<\cdots \\
@V{r_1}VV @V{r_1}VV\\
H_3(\Bbb R)/\Gamma_1' @<<< H_3(\Bbb R)/\Gamma_2' @<<<\cdots,
\endCD
$$
where the horizontal arrows in the bottom line
denote the natural projections.
Hence $T$ is isomorphic to the Heisenberg odometer associated with the sequence $(\Gamma_j')_{j=1}^\infty$.
\qed
\enddemo

\head 6.
Joining structure of Heisenberg odometers
\endhead

In this section we describe the 2-fold ergodic self-joinings, as dynamical systems, for all Heisenberg odometers.
We recall that a  2-fold  self-joining of a $G$-action $T=(T_g)_{g\in G}$ on a standard probability space $(X,\goth B,\mu)$ is a measure $\lambda$ on the product space $(X\times X,\goth B\otimes\goth B)$ such that $\lambda$ is invariant  under the diagonal action $(T_g\times T_g)_{g\in G}$ and the projections of $\lambda$ to each of the two coordinates are $\mu$.
The corresponding dynamical system $(X\times X, \lambda, (T_g\times T_g)_{g\in G})$ is also called a 2-fold self-joining of $T$.
A 2-fold self joining of $T$ is called off-diagonal if it is supported on the graph of a transformation of $(X,\mu)$ commuting with $T$.
In this section we  show that 2-fold self-joinings of Heisenberg odometers need not be Heisenberg odometers themselves, even in the transitive case.
In particular, there is an ergodic non-transitive 2-fold self-joining of a transitive Heisenberg odometer.
This is in  sharp contrast with the actions with pure point spectrum.
Indeed, each transitive action $T$ of a locally compact second countable group $G$ with pure point spectrum can be represented as an action by rotations (via an onto group homomorphism from $G$ to $K$) on a homogeneous space $K/H$, where $K$ is a compact second countable group and $H$ is a subgroup of $K$ \cite{Ma3}.
It is a standard argument  (see \cite{dJRu} for example) that each ergodic 2-fold self-joining of $T$ considered as a dynamical system   is a factor of an ergodic 2-fold self-joining  of the transitive $G$-action by rotations on $K$.
All such self-joinings are off-diagonal, isomorphic to the original $G$-action on $K$, and hence transitive themselves.
Since each factor of  a transitive system is transitive, it follows that each ergodic 2-fold self-joining of a transitive action with pure point spectrum is transitive.
Similar reasoning shows that the family of
   $G$-actions with pure point spectrum
is also closed under the operation of taking ergodic 2-fold self-joinings.
\footnote{For each  $G$-action $T$ with pure point spectrum there are a compact group $K$, a closed subgroup $H\subset K$ and a continuous homomorphism $\phi:G\to K$ with dense range in $K$ such that $T_gkH=\phi(g)kH$ for all $g\in G$ and $k\in K$ \cite{Ma3}. Hence, by \cite{dJRu}, every ergodic 2-fold self-joining of $T$ is a factor of an ergodic 2-fold self-joining of the $G$-action by rotations (via $\phi$) on $K$. The latter joining is off-diagonal; it is isomorphic as a dynamical system to the $G$-action by rotations (via $\phi$) on $K$. It remains to note that a factor of each action with pure point spectrum is an action with pure point spectrum.}

{\bf (A) Self-joinings of transitive Heisenberg odometers.} Let $\Gamma$ be a lattice in $H_3(\Bbb R)$.
In view of Proposition~1.3, we may assume without
loss of generality  that
$\Gamma=\{c(n/k)b(m)a(l)\mid n,m,l\in\Bbb Z\}$
for some $k\in\Bbb N$.
Every element $g\in H_3(\Bbb R)$ can be written uniquely as
$g=c(t_3)b(t_2)a(t_1)\gamma$
for some $\gamma\in\Gamma$ and $0\le t_3<1/k$, $0\le t_2<1$ and
$0\le t_1<1$.
Then the quotient space $H_3(\Bbb R)/\Gamma$ is a 3-torus
$$
\Bbb  T^3=\{( t_1,t_2,t_3)\mid
 0\le t_1<1,0\le t_2<1\text{ and }0\le t_3<1/k
\}.
$$
Our purpose now is to describe  the ergodic 2-fold self-joinings of the quotient (transitive) action
$T=(T_g)_{g\in H_3(\Bbb R)}$ of $H_3(\Bbb R)$
on this space.
We write it in  skew product form as follows:
$$
T_g(y,z)=(p(g)\cdot y,\alpha(g,y)+z),\tag6-1
$$
 where $(y,z)\in Y\times Z:=(\Bbb R/\Bbb Z)^2\times(\Bbb R/k^{-1}\Bbb Z)$, the symbol ``$\cdot$'' denotes the usual action of $\Bbb R^2$ on $Y$ by rotations and $\alpha:H_3(\Bbb R)\times Y\to
Z$ is the corresponding cocycle.
Let $\widetilde\Gamma:=p^{-1}(p(\Gamma))=\{c(t)b(m)a(n)\mid t\in\Bbb R,n,m\in\Bbb Z\}$.
Since $Y=H_3(\Bbb R)/\widetilde\Gamma$, the map $s:Y\ni(y_1,y_2)\mapsto b(y_2)a(y_1)\in H_3(\Bbb R)$, $0\le y_1,y_2<1$, is a cross-section of the quotient map $H_3(\Bbb R)\to H_3(\Bbb R)/\widetilde\Gamma$.
Denote by $h_s:H_3(\Bbb R)\times Y\to \widetilde\Gamma$ the corresponding choice-cocycle.
A direct calculation shows that
$$
\align
h_s(c(t_3),y) & = c(t_3),\\
h_s(b(t_2),y) & = b([t_2+y_2])c(-y_1[t_2+y_2]),\\
h_s(a(t_1),y) & = a([t_1+y_1])c(t_1y_2)
\endalign
$$
for all $t_1,t_2,t_3\in\Bbb R$ and $y=(y_1,y_2)$, $0\le y_1,y_2<1$,
where ``$[.]$'' denotes the integer part.
 By Remark~2.0 (in which $N=\widetilde\Gamma$),  we may assume without loss of generality that
 $\alpha=\phi\circ h_s$ for a group homomorphism $\phi:\widetilde\Gamma\to Z$.
 We now compute $\phi$.
 For $\widetilde\gamma=c(t)b(m)a(n)$ with $t\in\Bbb R$ and $n,m\in\Bbb Z$,
 we have $T_{\widetilde\gamma}\Gamma=c(t)\Gamma$.
 Applying \thetag{6-1},
 we obtain
 $$
 T_{\widetilde\gamma}(\widetilde\Gamma,0)=(\widetilde\Gamma, \alpha(\widetilde\gamma,\widetilde\Gamma))=(\widetilde\Gamma,t+k^{-1}\Bbb Z).
 $$
 Thus, we obtain that $\phi(c(t)b(m)a(n))=t+k^{-1}\Bbb Z$.
 This yields
  $$
\align
\alpha(c(t_3),y) & = t_3+k^{-1}\Bbb Z,\\
\alpha(b(t_2),y) & = -y_1[t_2+y_2]+k^{-1}\Bbb Z,\\
\alpha(a(t_1),y) & = t_1y_2+k^{-1}\Bbb Z
\endalign
$$
for all $t_1,t_2,t_3\in\Bbb R$ and $y=(y_1, y_2)$, $0\le y_1,y_2<1$.
For $d=(d_1,d_2)\in Y$, $0\le d_1,d_2<1$, we denote by $\alpha\circ d$ the cocycle 
$$
H_3(\Bbb R)\times Y\ni(g,z)\mapsto
\alpha\circ d(g,y):=\alpha(g, d+y)\in Z.
$$
By Remark~2.0,  we can represent 
 the cocycle
 $$
 \alpha\times\alpha\circ d:H_3(\Bbb R)\times Y\to Z\times Z
 $$
 as $\psi\circ h_s$, where
 $\psi$ is a continuous  group homomorphism from $\widetilde\Gamma$ to $Z\times Z$.
 To compute $\psi$, we first note that
 $$
 h_s\circ d(g,y)=h_s(g,y+d)=s(p(g)\cdot(y+d))^{-1}gs(y+d).
 $$
 Substituting $y:=\widetilde\Gamma$ and $g:=c(t)b(m)a(n)$ we obtain
 that
 $$
 h_s\circ d(g,\widetilde\Gamma)=s(d+\widetilde\Gamma)^{-1}gs(d+\widetilde\Gamma)=c(nd_2-md_1)g.
 $$
 This formula, the formula for $\phi$ and Remark~2.0 yield that
 $$
\psi(c(t)b(m)a(n))=(t+k^{-1}\Bbb Z,t+nd_2-md_1+k^{-1}\Bbb Z).\tag6-2
$$
If $d$ is aperiodic, i.e. $d_1$ or $d_2$ is irrational, then
the subgroup $\psi(\widetilde\Gamma)$ is dense in $ Z\times Z$.
It now follows from Remark~2.1
the essential range $E( \alpha\times\alpha\circ d)$
of the cocycle $ \alpha\times\alpha\circ d$ is $Z\times Z$.
We note however that
$\psi(\widetilde\Gamma)\ne  Z\times Z$.
On the other hand, if $d$ is periodic,  i.e. $d_1$ and $d_2$ are both rational, then there is a positive integer $q=q(d)$ such that
$$
E( \alpha\times\alpha\circ d)=\psi(\widetilde\Gamma)=\bigcup_{j=0}^{q-1}\{(t+k^{-1}\Bbb Z,t+ j/({qk})+k^{-1}\Bbb Z\mid 0\le t<1\}.\tag6-3
$$
We denote this subgroup of $Z\times Z$ by $D_q$.

Let $\Delta_d$ denote the measure on $Y\times Y$ sitting on the subset $\{(y, d+y)\mid y\in Y\}$ and projecting on  the Haar measure on $Y$ along each of the two coordinate projections.
Given a closed subgroup $\Lambda$ in $Z\times Z$, we denote by $\lambda_\Lambda$ the Haar measure on $\Lambda$.
We consider it as a measure on $ Z\times Z$.
Given $z\in Z$, we denote by $\lambda_\Lambda\circ z$ the image of $\lambda_\Lambda$ under the rotation
$Z\times Z\ni(z_1,z_2)\mapsto (z_1,z_2+z)\in Z\times Z$.

\proclaim{Proposition 6.1}
The set $J^e_2(T)$ of all ergodic 2-fold self-joinings of $T$ is the union of two families as follows:
\footnote{Provided that a special  representation of $T$ as a skew product is chosen as above.}
$$
J^e_2(T)=\{\Delta_d\times\lambda_{Z\times Z}\mid  \text{ $d$ is aperiodic}\}\cup
\{\Delta_d\times\lambda_{D_{q(d)}}\circ z\mid \text{$d$ is periodic, $z\in Z$} \}.
$$
Every joining from the first family is a non-transitive dynamical system and every joining from the second family is a transitive dynamical system.
\endproclaim

\demo{Proof} Take $\kappa\in J^e_2(T)$.
Represent $T$ in the skew product form \thetag{6-1}.
Then $\kappa$ is a measure on the space  $Y\times Z\times Y\times Z$.
Project $\kappa$ onto $Y\times Y$.
The projection is an ergodic 2-fold self-joining on the transitive action of the 2-torus $Y$ on itself.
Hence there is $d\in Y$ such that this projection is $\Delta_d$.
It is a standard argument that the dynamical system $( (T_g\times T_g)_{g\in H_3(\Bbb R)},\kappa)$ can now be represented as the skew product over the same  (as $T$) base $Y$ but with another cocycle $\alpha\times\alpha\circ d$.
It remains to apply the argument preceding the statement of Proposition~6.1 to obtain the first assertion of the proposition.

Let $V$ denote the action of $\widetilde\Gamma$ on the compact group $\overline{\psi(\widetilde\Gamma)}$ by translations via $\psi$, i.e. an element $\gamma\in \widetilde\Gamma$ acts as the translation by $\psi(\gamma)$.
It is easy to see that the action  $( (T_g\times T_g)_{g\in H_3(\Bbb R)},\kappa)$ is isomorphic to Ind$_{\widetilde\Gamma}^{H_3(\Bbb R)}(V)$,
the second claim of the proposition follows from the fact whether or not
$\psi(\widetilde\Gamma)$ is closed and the following lemma.
\qed
\enddemo

\proclaim{Lemma 6.2}
Let $G$ be a locally compact second countable group and let $H$ be a closed subgroup of $G$.
A measure preserving action $V$ of $H$ is transitive if and only if the induced measure preserving action $\text{\rom{Ind}}_H^G(V)$ of $G$ is transitive.
\endproclaim
\demo{Proof}
If $V$ is non-transitive then, obviously,
$\text{\rom{Ind}}_H^G(V)$ is non-transitive.
Suppose now that $V$ is transitive.
Denote by $N$ a closed subgroup of $H$ such that $V$ is isomorphic to the natural $H$-action on the homogeneous space $H/N$.
Then $V$ is isomorphic to the induced action $\text{\rom{Ind}}_N^H(I)$, where $I$ denotes the trivial action of $N$ on a singleton.
We now have (see \cite{Ma1}) that
$$
\text{\rom{Ind}}_H^G(V)=\text{\rom{Ind}}_H^G(\text{\rom{Ind}}_N^H(I))=\text{\rom{Ind}}_N^G(I),
$$
i.e.  $\text{\rom{Ind}}_H^G(V)$ is isomorphic to the natural $G$-action on $G/N$.
\qed
\enddemo

\remark{Remark 6.3}
We note that if an action of an Abelian locally compact group is transitive then every ergodic 2-fold self-joining of this action is also transitive (and isomorphic to the original action).
As follows from Proposition~6.1, this is no longer true for the transitive actions of the Heisenberg group.
\endremark

We are going to show now that an ergodic 2-fold self-joining of a transitive odometer can be rather ``far from'' the Heisenberg odometers.
Though its maximal spectral type is discrete, the subgroup of 1-dimensional components of the spectrum is not  off-rational.

\example{Example 6.4}
Consider a particular case when $k=1$, i.e. $\Gamma=H_3(\Bbb Z)$.
Suppose that $d_1$ is irrational (see \thetag{6-2}).
Then the action $R=(R_g)_{g\in H_3(\Bbb R)}$ of $H_3(\Bbb R)$ defined on $Y\times Z\times Z$ by
$$
R_g(y,z):=(p(g)\cdot y,z+\psi(h_s(g,y))),\qquad y\in Y,z\in Z\times Z,\tag6-4
$$
is an ergodic 2-fold self-joining of $T$.
Denote by $U_R$  the Koopman unitary representations of $H_3(\Bbb R)$ generated by $R$.
We are going to decompose $U_R$ into irreducible components.
It follows from \thetag{6-4} that
$$
L^2(H_3(\Bbb R)/\Gamma)=\bigoplus_{n,m=1}^\infty L^2(Y)\otimes\eta_n\otimes\eta_m,
$$
and each subspace $L^2(Y)\otimes\eta_n\otimes\eta_m$ is $U_R$-invariant, where $\eta_n(t):=\exp{(2\pi int)}$, $t\in[0,1)$, $n\in\Bbb Z$.
Since
$$
R_{c(t)}(y,z)=(y,z+(t,t)+\Bbb Z^2)\qquad{\text{for all } }t\in\Bbb R,
$$
It follows that $U_R({c(t)})\restriction(L^2(Y)\otimes\eta_n\otimes\eta_m)
=\exp{(2\pi it(n+m))}\text{Id}$.
This yields
$$
U_R=\bigoplus_{\gamma\in\Bbb Z}\bigoplus_{r=1}^\infty\pi_{\gamma}\oplus
\bigoplus_{k\in\Bbb Z}U_R\restriction(L^2(Y)\otimes\eta_k\otimes\eta_{-k})\tag6-5
$$
We  note that
  the space $\Cal H_1:=\bigoplus_{k\in\Bbb Z}L^2(Y)\otimes\eta_k\otimes\eta_{-k}$ is the subspace of vectors fixed by $(U_R(c(t)))_{t\in\Bbb R}$.
It is straightforward to verify that $\Cal H_1=L^2(Y)\otimes\Cal T$, where
$$
\Cal T:=\{v\in L^2(Z\times Z)\mid v(z_1,z_2)=v(z_1+z,z_2+z)\text{ at a.a. }z_1,z_2\text{ for each }z\in Z \}.
$$
Hence $\Cal T=L^2(\goth F)$, where $\goth F$ is the sub-$\sigma$-algebra of sets in $Z\times Z$ measurable with respect to the map $Z\times Z\ni(z_1,z_2)\mapsto z_1-z_2\in Z$.
Thus, the $\sigma$-algebra $\goth B_Y\otimes\goth F$ is a factor
of $R$ such that $L^2(\goth B_Y\otimes\goth F)=\Cal H_1$.
In view of \thetag{6-4}, the restriction $\widehat R=(\widehat R_g)_{g\in G}$ of $R$ to this factor can be written as
$$
\widehat R_g(y,z)=(p(g)\cdot y,z+\widetilde\psi(h_s(g,y))),
$$
where $\widetilde\psi:\widetilde\gamma\to Z$ is the difference of the first and the second coordinate of $\widetilde\psi$, i.e., by \thetag{6-2},
$\widetilde\psi(c(t)b(m)a(n))=md_1-nd_2+\Bbb Z$, $t\in\Bbb R$, $n,m\in\Bbb R$.
This means that $\widehat R$ is induced from the the following action of $\widetilde\Gamma$ on $Z$ by rotations
$$
\widetilde\Gamma\times Z\ni(c(t)b(m)a(n),z)\mapsto z+ md_1-nd_2+\Bbb Z.
$$ 
This action has a pure point spectrum and hence $\widehat R$ has pure point spectrum.
It is straightforward to verify that 
$$
U_R\restriction\Cal H_1=\bigoplus_{\alpha,\beta\in\Omega_1}\pi_{\alpha,\beta},
$$
where
$\Omega_1= \Bbb Z^2+\pmatrix d_2 & 0\\
0 & -d_1\endpmatrix\Bbb Z^2$.
Since $d_1$ is irrational, $\Omega_1$ is not off-rational.
In view of \thetag{6-5},
$$
U_R=
\bigoplus_{\gamma\in\Bbb Z}\bigoplus_{r=1}^\infty\pi_{\gamma}\oplus
\bigoplus_{\alpha,\beta\in\Omega_1}\pi_{\alpha,\beta}.
$$
Since $\Omega_1$ is not off-rational, $R$ is not isomorphic  to any Heisenberg odometer.\footnote{In fact, we have shown that $R$ is not spectrally equivalent to any Heisenberg odometer.}

\endexample

{\bf (B)} 
{\bf Self-joinings of general Heisenberg odometers.}
We now consider a Heisenberg odometer $T$ associated to a sequence of latices  $\Gamma_1\supset\Gamma_2\supset\cdots$ in $H_3(\Bbb R)$.
Then every ergodic 2-fold self-joining of $T$ is the inverse limit of a sequence of ergodic 2-fold self-joinings of the transitive actions of $H_3(\Bbb R)$ on the homogeneous spaces $H_3(\Bbb R)/\Gamma_j$, $j\in \Bbb N$, viewed as factors of $T$.
Using this simple observation and Proposition~6.1 we can explicitly describe all ergodic 2-fold self-joinings of $T$.

The  $T$-action can be represented as a skew product.
The space of this action is the product $Y\times Z$ of two compact Abelian groups
$Y:=\projlim_{j\to\infty}\Bbb R^2/p(\Gamma_j)$ and $Z:=\projlim_{j\to\infty}Z_j$, where $Z_j:=\Bbb R/\xi_{\Gamma_j}\Bbb Z$.
The factor-space $Y$ is the factor of $T$ generated by all proper functions of $T$; the restriction of $T$ to this factor has pure point spectrum.
The corresponding cocycle of this factor with values in $Z$ is ergodic.\footnote{See Section~7 for details.}
Given $d\in Y$, we denote by $\Delta_d$ the image of the Haar measure on $Y$ under the map
$Y\ni y\mapsto(y,y+d)\in Y\times Y$.
Every element $d\in Y$ is a sequence $(d_j)_{j\in\Bbb N}$ of elements $d_j\in\Bbb R^2/p(\Gamma_j)$
such that $d_{j+1}$ maps to $d_j$ under the natural projection $\Bbb R^2/p(\Gamma_{j+1})\to\Bbb R^2/p(\Gamma_j)$ for each $j$.
In a similar way,
every element $z\in Z$ is a sequence $(z_j)_{j\in\Bbb N}$ of elements $z_j\in Z_j$
such that $z_{j+1}$ maps to $z_j$
 under the natural projection $Z_{j+1}\to
Z_j$ for each $j$.
If $d_j$ is periodic then we denote by  $D_j$ be the closed subgroup of $Z_j\times Z_j$ associated with $d_j$
in the way described in the part (A)  of this section (see \thetag{6-3}).
We note that $D_j$ contains the diagonal of  $Z_j\times Z_j$ as a subgroup of finite index.
Moreover, $D_{j+1}$ maps onto $D_j$ under the natural projection
$Z_{j+1}\to
Z_j$ for each $j$.
Hence a projective limit $D_d:=\projlim_{j\to\infty}Z_j$
is well defined.
It is a closed subgroup of $Z$.
As in Proposition~6.1 we will use the following notation.
Given a closed subgroup $\Lambda$ of $Z\times Z$, let $\lambda_\Lambda$ stand for the Haar measure on $\Lambda$.
Given $z\in Z$, let $\lambda_\Lambda\circ z$ denote the image of $\lambda_\Lambda$ viewed as a measure on $Z\times Z$ under the rotation $(z_1,z_2)\mapsto(z_1,z_2+z)$ of $Z\times Z$.

\proclaim{Theorem 6.5}
The set $J^e_2(T)$ of all ergodic 2-fold self-joinings of $T$ is the union of the following two families:\footnote{Provided that some special cocycle representing $T$ as a skew product over its maximal factor with pure point spectrum is chosen.}
$$
\align
J^e_2(T)= &\,
\{\Delta_d\times\lambda_{Z\times Z}\mid  \text{ $d=(d_j)_{j\in\Bbb N}$ with $d_j$ aperiodic for each $j$}\}\\
&\cup
\{\Delta_d\times\lambda_{D_{d}}\circ z\mid \text{$d=(d_j)_{j\in\Bbb N}$ with $d_j$ periodic for each $j$ and
$z\in Z$} \}.
\endalign
$$
\endproclaim

\head 7.   On spectral determinacy of Heisenberg odometers
\endhead

{\bf (A)} {\bf The case of transitive odometers.} Let $T$ be an ergodic action of $H_3(\Bbb R)$ on a standard probability space $(X,\mu)$.
Denote by $U_{T}$ the corresponding Koopman unitary representation of $H_{3}(\Bbb R)$.

\proclaim{Theorem 7.1}
If $U_{T}$ is unitarily equivalent to the Koopman unitary representation
generated by the  action $Q$ of  $H_{3}(\Bbb R)$ by translations on $H_{3}(\Bbb R)/\Gamma$ for a lattice $\Gamma$ in $H_3(\Bbb R)$ then $T$ is  isomorphic to $Q$.
\endproclaim

\demo{Proof}
Without loss of generality we may assume that there is $k\ge 1$ such that
$\Gamma=\{c(i_{3}/k)b(i_2)a(i_1)\mid i_1,i_2,i_3\in\Bbb Z\}$.
By Theorem~4.1,
$$
U_T=\bigoplus_{(\alpha,\beta)\in\Bbb Z^2}\pi_{\alpha,\beta}\oplus\bigoplus_{0\ne n\in\Bbb Z}\bigoplus_{j=1}^{|n|}\pi_{nk}\tag7-1
$$
up to the unitary equivalence.
Denote by $\Cal H$ the subspace of $L^2(X,\mu)$ where the unitary representation
 $\bigoplus_{(\alpha,\beta)\in\Bbb Z^2}\pi_{\alpha,\beta}$ is defined.
 Since each unitary representation $\pi_{\alpha,\beta}$ is 1-dimensional, it follows from the classical von Neumann theory of Abelian actions with pure point spectrum that there is a measure preserving factor map $q:(X,\mu)\to(\Bbb R^2/\Bbb Z^2,\lambda_{\Bbb R^2/\Bbb Z^2})$  intertwining, for each $g\in H_3(\Bbb R)$, the transformation $T_g$ with the rotation on $\Bbb R^2/\Bbb Z^2$ by $p(g)$
 and such that $\Cal H=\{f\circ q\mid f\in L^2(\Bbb R^2/\Bbb Z^2) \}$.
 Since $\pi_{nk}(c(t))v=e^{2\pi inkt}v$ for each vector $v$ of the Hilbert space where $\pi_{nk}$ is defined,
 it follows from \thetag{7-1} that $\Cal H$ is the subspace of vectors fixed by the unitary 1-parameter group $\{U_T(c(t))\mid t\in\Bbb R\}$.
Therefore by the von Neumann mean ergodic theorem,  $\frac 1{2N}\int_{-N}^NU_T(c(t))\,d t$ converges strongly as $N\to\infty$ to the orthogonal projection to $\Cal H$ which is the mathematical expectation to the sub-$\sigma$-aldebra of $\{T_{c(t)}\mid t\in\Bbb R\}$-invariant subsets.
Thus $q$ is the $\{T_{c(t)}\mid t\in\Bbb R\}$-ergodic decomposition and $(\Bbb R^2/\Bbb Z^2,\lambda_{\Bbb R^2/\Bbb Z^2})$ is the space of $\{T_{c(t)}\mid t\in\Bbb R\}$-ergodic components.
We note that $U_T(c(k^{-1}))=\text{Id}$.
Moreover, if $0<s<k^{-1}$ then the unitary operator $U_T(c(s))$ has no non-trivial fixed vectors in the orthocomplement to $\Cal H$ in $L^2(X)$.
Hence the restriction of $\{T_{c(t)}\mid t\in\Bbb R\}$ to almost every ergodic component is isomorphic to the periodic flow with the smallest period $k^{-1}$.
Thus $(X,\mu)$ splits into the product
$$
(X,\mu)=(\Bbb R^2/\Bbb Z^2,\lambda_{\Bbb R^2/\Bbb Z^2})\times(\Bbb R/k^{-1}\Bbb Z,\lambda_{\Bbb R/k^{-1}\Bbb Z})
$$
such that $q(y,z)=y$ and
$$
T_{c(t)}(y,z)=(y,t+z)\tag7-2
$$
for a.a. $y\in \Bbb R^2/\Bbb Z^2$ and $z\in\Bbb R/k^{-1}\Bbb Z$.
Thus there exists a Borel cocycle
$$
\delta:H_3(\Bbb R)\times\Bbb R^2/\Bbb Z^2\to\Bbb R/k^{-1}\Bbb Z
$$
such that
$$
T_g(y,z)=(p(g)\cdot y,\delta(g,y)+z) \qquad\text{for a.a. }(y,z)\in X.
$$
It follows that $T$ is a transitive action on a 3-dimensional torus.
Hence there is a lattice $\Gamma'$ in $H_3(\Bbb R)$
such that $T$ is isomorphic to a homogeneous $H_3(\Bbb R)$-action on $H_3(\Bbb R)/\Gamma'$.
Since the 1-dimensional spectrum of $T$ is $\{\pi_{\alpha,\beta}\mid (\alpha,\beta)\in\Bbb Z^2\}$ we deduce from Theorem~4.2 that  $p(\Gamma')=\Bbb Z^2$.
It follows from \thetag{7-2} that $\xi_{\Gamma'}=k^{-1}$.
Thus $p(\Gamma)=p(\Gamma')$ and $\xi_\Gamma=\xi_{\Gamma'}.$
Hence by Corollary~4.3, $\Gamma$ and $\Gamma'$ are conjugate, i.e. $T$ and $Q$ are isomorphic.
\qed
\enddemo

\remark{Remark 7.2}
\roster
\item
In fact, a stronger claim has been proved: if the maximal spectral types of $T$ and $Q$ are equivalent then $T$ and $Q$ are isomorphic.
\item
As  a byproduct, we have discovered the following interesting fact: let $T$ be a transitive Heisenberg odometer and let $U_T$ denote the corresponding Koopman unitary representation of $H_3(\Bbb R)$.
If $V$ is a unitary representation of $H_3(\Bbb R)$ with the same maximal spectral type as $U_T$ but with a different multiplicity function, then $V$ is not unitarily equivalent to the Koopman unitary representation of $H_3(\Bbb R)$ generated by any probability preserving action of $H_3(\Bbb R)$.
\endroster
\endremark

{\bf (B)} {\bf The general case.} Here we prove the main result of this section: the $H_3(\Bbb R)$-odometers are not spectrally determined.

\proclaim{Lemma 7.3}
Let $\xi=\bigcup_{j=1}^\infty d_j^{-1}\Bbb Z$ be an off-rational subgroup of $\Bbb R$ and let
$S=\bigcup_{j=1}^\infty A_j^{-1}\Bbb Z^2$ be an off-rational subgroup of $\Bbb R^2$.
We set $Z:=\projlim_{j\to\infty}\Bbb R/d_j\Bbb Z$ and
$Y:=\projlim_{j\to\infty}\Bbb R^2/A_j^{*}\Bbb Z^2$.
Denote by $\lambda_Z$ and $\lambda_Y$ the Haar measure on $Z$ and $Y$ respectively.
Denote by $\overline r_j:Z\to \Bbb R /d_j^{-1}\Bbb Z$ and $\overline q_j:Y\to \Bbb R^2/A_j^{*}\Bbb Z^2$
the canonical projections.
Let $\vartheta:\Bbb R\to Z$
and $\varphi:\Bbb R^2\to Y$ be continuous  homomorphisms defined by
$\overline r_j(\vartheta(t))=t+d_j^{-1}\Bbb Z$ and
$\overline q_j(\varphi(c))=c+A_j^{*}\Bbb Z^2$
for all $t\in\Bbb R$, $c\in\Bbb R^2$, $j\in\Bbb N$.
Let $T$ be an ergodic measure preserving action  on a probability space $(X,\mu)$.
Suppose that the associated Koopman unitary
representation $U_T$ of $H_3(\Bbb R)$ decomposes  into irreducibles as
$$
U_T=\bigoplus_{(\alpha,\beta)\in S}\pi_{\alpha,\beta}\oplus\bigoplus_{\gamma\in\xi}\bigoplus_{l=1}^{m_\gamma}\pi_\gamma\tag7-4
$$
for  a spectral  multiplicity map $m:\xi\ni\gamma\mapsto m_\gamma\in\Bbb N\cup\{\infty\}$.
Then $(X,\mu)$ is isomorphic to the product $(Y\times Z,\lambda_Y\times\lambda_Z)$ in such a way that
$$
T_g(y,z)=(D_gy,\alpha(g,y)+z),
$$
where  
$$
D:H_3(\Bbb R)\times Y\ni(g,y)\mapsto D_gy:=\varphi(p(g))+y\in Y
$$
is an ergodic action of $H_3(\Bbb R)$ on $Y$ with pure point spectrum and 
$$
\alpha:H_3(\Bbb R)\times Y\ni (g,y)\mapsto\alpha(g,y)\in Z
$$ 
is an ergodic  cocycle of $D$
 such that $\alpha(c(t),y)=\vartheta(t)$ for all $t\in\Bbb R$ and $y\in Y$.
\endproclaim

\demo{Proof}
We  have a decomposition of $L^2(X,\mu)$ corresponding to \thetag{7-4}:
$$
L^2(X,\mu)=\bigoplus_{(\alpha,\beta)\in S}\Cal H_{\alpha,\beta}\oplus\bigoplus_{\gamma\in\xi}\bigoplus_{l=1}^{m_\gamma}\Cal H_{\gamma,l}.
$$
We now define a unitary representation $W=(W(z))_{z\in Z}$ of $Z$ in $L^2(X,\mu)$ by setting
$$
W(z)v=\cases
v &\text{ if $v\in\Cal H_{\alpha,\beta}$, $(\alpha,\beta)\in S$ and }\\
e^{2\pi i\gamma\overline r_j(z)}v
&\text{ if $v\in\Cal H_{\gamma,l},\gamma\in\xi,1\le l\le m_\gamma$.}
\endcases
$$
We consider here $\overline r_j(z)$ as a real number from the segment $[0,d_j)$.
Then $W$ is continuous in the weak operator topology.
Moreover,
$$
W({\vartheta(t)})=U_T(c(t))\quad\text{ for all }t\in\Bbb R.\tag7-5
$$
Since the group Aut$(X,\mu)$ of all $\mu$-preserving transformations of $X$ is a closed subgroup the unitary group  $\Cal U(L^2(X,\mu))$ and
the subgroup $\vartheta(\Bbb R)$ is dense in $Z$, it follows that
$W(z)\in\text{Aut}(X,\mu)$ for all $z\in Z$.
Thus, $\{W(z)\mid z\in Z\}$ is a compact group of transformations commuting with the action $T$.
Therefore there are a probability space $(Y,\nu)$ and an action $D=(D_g)_{g\in H_3(\Bbb R)}$ of $H_3(\Bbb R)$ on it such that $(X,\mu)$ is isomorphic to the product
 $(Y\times Z,\nu\times\lambda_Z)$ such that
 $$
 T_g(y,z)=(D_gy,\alpha(g,y)+z)\quad{ and }\quad W(z')(y,z)=(y,z+z')
 $$
for all $g\in H_3(\Bbb R)$ and $z'\in Z$, where $\alpha:H_3(\Bbb R)\times Y\to Z$ is a measurable cocycle of $D$ (see, e.g. \cite{dJRu}).
It follows from \thetag{7-5} that $\alpha(c(t),y)=\vartheta(t)$ for all $t\in\Bbb R$.
Since
$$
\align
\bigoplus_{(\alpha,\beta)\in S}\Cal H_{\alpha,\beta}&=\{v\in L^2(X,\mu)\mid U_T(c(t))v
=v\text{ for all }t\in\Bbb R\}\\
&=\{v\in L^2(X,\mu)\mid W(z)v=v\text{ for all }z\in Z\},
\endalign
$$
it follows that $L^2(Y)=\bigoplus_{(\alpha,\beta)\in S}\Cal H_{\alpha,\beta}$ and the Koopman unitary representation of $H_3(\Bbb R)$ generated by $D$  decomposes into irreducible ones as $U_D=\bigoplus_{(\alpha,\beta)\in S}\pi_{\alpha,\beta}$.
Therefore without loss of generality we may  assume that $Y$ is a compact Abelian group  dual  to
$S$ and $D_gy=\varphi (p(g))+ y$ for all $y\in Y$ and $g\in H_3(\Bbb R)$.
\qed
\enddemo

Now we prove a converse to Lemma~7.3.

\proclaim{Lemma 7.4}
Let  $S$, $\xi$, $Y$, $Z$  and $D$ be as  in Lemma~7.3.
\roster
\item"\rom{(i)}"
Given a cocycle $\beta$ of $D$ with values in $Z$ such that
$$
\beta(c(t),y)=\vartheta(t) \qquad\text{for all $t\in\Bbb R$,}\tag7-6
$$
there is a  map $m:\xi\ni\gamma\mapsto m_\gamma\in\Bbb N\cup\{\infty\}$ such that
the Koopman unitary representation $U_{D^\beta}$ of $H_3(\Bbb R)$ generated by the $\beta$-skew product extension $D^\beta$ of $D$
 is unitarily equivalent to  \thetag{7-4}.
\item"\rom{(ii)}"
Given two cocyles $\beta,\beta'$ of $D$ with values in $Z$ satisfying \thetag{7-6}, the skew product extensions $D^\beta$ and $D^{\beta'}$ of $D$ are isomorphic if and only if there exist an element $y_0\in Y$ and a Borel map  $f:Y\to Z$ such that
$$
\beta(g,y+y_0)=f(y)+\beta'(g,y)-f(\varphi(p(g))+y)
$$
for each $g\in H_3(\Bbb R)$ at a.a. $y\in Y$.
\endroster
\endproclaim

\demo{Proof}
(i) We first  note that  $\xi$ is isomorphic to the dual $\widehat Z$ of $Z$.
 Hence we can identify  an element $\gamma$ of $\xi$ with a character $\xi_\gamma\in L^2(Z,\lambda_Z)$.
 Then we obtain a decomposition of $L^2(Y\times Z,\lambda_Y\times\lambda_Z)$ into an orthogonal sum of $U_{D^\beta}$-invariant subspaces:
 $$
 L^2(Y\times Z,\lambda_Y\times\lambda_Z)=\bigoplus_{\gamma\in\xi}L^2(Y,\lambda_Y)\otimes\xi_\gamma.
 $$
 It follows from the property of $\beta$ that
 $$
 U_{D^\beta}(c(t))(v\otimes\xi_\gamma)=\xi_\gamma(\vartheta(t)) v\otimes\xi_\gamma=
e^{2\pi i\gamma t}  v\otimes\xi_\gamma, \qquad v\in L^2(Y,\lambda_Y).
 $$
 Therefore $U_{D^\beta}\restriction (L^2(Y,\lambda_Y)\otimes\xi_\gamma)=\bigoplus_{1}^{m_\gamma}\pi_\gamma$ for some multiplicity $m_\gamma\in\Bbb N\cup\{\infty\}$ for each $\gamma\in\xi$, $\gamma\ne 0$.
 Of course,
$U_{D^\beta}\restriction (L^2(Y,\lambda_Y)\otimes 1)= \bigoplus_{(\alpha,\beta)\in S}\Cal H_{\alpha,\beta}$.
It follows that~\thetag{7-4} holds, as desired.

(ii) Since $D^\beta$ and $D^{\beta'}$ are isomorphic, it follows from a rather standard  description of the  centralizer of an ergodic compact group extension of an action by rotations on a compact Abelian group (see, e.g., \cite{Da1}) that there exist
an element $y_0\in Y$, a continuous  automorphism $\theta$ of $Z$ and a Borel map  $f:Y\to Z$ such that
$$
\beta(g,y+y_0)=f(y)+\theta(\beta'(g,y))-f(\varphi(p(g))+y)
$$
for each $g\in H_3(\Bbb R)$ at a.a. $y\in Y$.
Now \thetag{7-6} yields that $\theta(\vartheta( t))=\vartheta(t)$ for all $t\in\Bbb R$.
Since $\vartheta$ has a dense range in $Z$, it follows that $\theta$ is the identity, as desired.
\qed
\enddemo

Consider now the case where $T$ is the $H_3(\Bbb R)$-odometer associated with  a sequence $\Gamma_1\supset\Gamma_2\supset\cdots$ of lattices in $H_3(\Bbb R)$.
Denote by $X$ the space of $T$.
We set $Y_j:=\Bbb R^2/p(\Gamma_j)$ and $Z_j:=\Bbb R/\xi_{\Gamma_j} \Bbb Z$.
Recall that the homogeneous space $H_3(\Bbb R)/\Gamma_j$ is isomorphic to the product $Y_j\times Z_j$ in such a way that the natural $H_3(\Bbb R)$-action has a skew product form as follows
$$
g\cdot(y,z):=(p(g)\cdot y,\alpha_j(g,y)+z)
\tag7-7
$$
for all $(y,z)\in Y_j\times Z_j$ and $ g\in H_3(\Bbb R)$,
where $\alpha_j:H_3(\Bbb R)\times Y_j\to\Bbb R/\xi_{\Gamma_j}\Bbb Z$
is a cocycle satisfying
$$
\alpha_j(c(t),y)=t+\xi_{\Gamma_j}\Bbb Z
\tag7-8
$$
The natural projection $H_3(\Bbb R)/\Gamma_j\to H_3(\Bbb R)/\Gamma_{j-1}$  corresponds  to the mapping
$$
Y_j\times Z_j\ni(y,z)\mapsto(q_j(y),r_j(z+f_j(y)))\in Y_{j-1}\times Z_{j-1},
$$
where $q_j:Y_j\to Y_{j-1}$ and $r_j:Z_j\to Z_{j-1}$ are the natural projections and $f_j:Y_j\to Z_j$ is a Borel map.\footnote{This follows from the fact that the projection $H_3(\Bbb R)/\Gamma_j\to H_3(\Bbb R)/\Gamma_{j-1}$ intertwines the actions of the center of $H_3(\Bbb R)$ on these homogeneous  spaces while $Y_j$ is the space of  ergodic components for the action of the center.}
It is easy to see that the $H_3(\Bbb R)$-equivariance of the map $H_3(\Bbb R)/\Gamma_j\to H_3(\Bbb R)/\Gamma_{j-1}$ is equivalent to the following property
$$
\alpha_{j-1}(g,q_j(y))=r_j(f_j(p(g)\cdot y)+\alpha_j(g,y)-f_j(y))
$$
for all $g\in G$ and $y\in Y_j$.
Therefore, replacing $\alpha_j$ with a cohomologous cocycle, we may assume without loss of generality that
$$
\alpha_{j-1}(g,q_j(y))=r_j(\alpha_j(g,y))
\tag7-9
$$
for every $j$.
Of course, this cohomological change will not affect  property \thetag{7-8}.
We now set $Y:=\projlim_{j\to\infty} Y_j$
 and $Z:=\projlim_{j\to\infty}Z_j$ in the category of compact Abelian groups furnished with the Haar measures.
 Denote by $\overline{q}_j:Y\to Y_j$ and $\overline r_j:Z\to Z_j$ the canonical projections.
We note that $Z$ and $Y$
are the same objects as in Lemmata~7.3 and 7.4 if we set
$\xi:=\bigcup_{j=1}^\infty\xi_{\Gamma_j}^{-1}\Bbb Z$
and $S:=\bigcup_{j=1}^\infty p(\Gamma_j)^*$.
Recall that $D$ is an ergodic action of $H_3(\Bbb R)$ on $Y$ with pure point spectrum $S$ (see the statement of Lemma~7.3).
In view of \thetag{7-9}, a cocycle $\alpha:H_3(\Bbb R)\times Y\to Z$ of $D$ is well defined by
$$
\overline r_j(\alpha(g,y))=\alpha_j(g,\overline q_j(y)), \qquad j\in\Bbb N.\tag7-10
$$
It follows from \thetag{7-8} that $\alpha(c(t),y)=\vartheta(t)$ for all $t\in\Bbb R$ and $y\in Y$.
Since $X=\projlim_{j\to\infty}H_3(\Bbb R)/\Gamma_j$, we obtain that $X$ is measure theoretically isomorphic to $Y\times Z$ in such a way that  $T$ is the $\alpha$-skew product extension $D^\alpha$ of $D$.

\definition{Definition 7.5}
Let $Y$, $Z$ and $D$ be as above.
A Borel cocycle $\alpha$ of $D$ with values in $Z$ is called {\it finitary} if  there are cocycles $\alpha_j$
of the natural action of $H_3(\Bbb R)$ on $Y_j$ with values in $Z_j$ such that \thetag{7-10} is satisfied.
\enddefinition

\remark{Remark 7.6}
For $y_0\in Y$, we define a $Z$-valued cocycle $\beta$ of $D$ by setting $\beta(g,y):=\alpha(g,y+y_0)$.
If $\alpha$ is finitary then so is $\beta$.
\endremark

We now state one of the main results of this section.

\proclaim{Theorem 7.7}
Let $S$ be an off-rational subgroup of $\Bbb R^2$ and let $\xi$ be an off-rational subgroup of $\Bbb R$ such that $\tau(S)\supset\xi$.
Fix
 a nested sequence $\Gamma=(\Gamma_j)_{j=1}^\infty$ of lattices in $H_3(\Bbb R)$
 with $S_\Gamma=S$ and $\xi_\Gamma=\xi$ \footnote{It  exists by Proposition 5.5.}.
Let $T$ be an ergodic action of $H_3(\Bbb R)$ on a standard probability space.
Suppose that \thetag{7-4} holds for the associated Koopman unitary representation $U_T$ of $H_3(\Bbb R)$.
As in Lemma~7.3, represent $T$ as a skew product extension $D^\alpha$ of $D$ for a cocycle $\alpha$ of $D$ with values in $Z$.
Then $T$ is   an $H_3(\Bbb R)$-odometer (which is f-isomorphic to the $H_3(\Bbb R)$-odometer associated with $\Gamma$) if and only if $\alpha$ is cohomologous to a finitary cocycle.
\endproclaim

\demo{Proof}
$(\Rightarrow)$
If $T$ is a Heisenberg odometer then (as was shown above) there exists a representation $T$ as $D^{\alpha'}$ with a finitary cocycle $\alpha'$.
By Lemma~7.4(ii), $\alpha$ is cohomologous to a `rotation of $\alpha'$ by an element of $Y$'.
As was noted in Remark~7.6, the class of finitary cocycles is invariant under such rotations.
Thus,  $\alpha$ is cohomologous to a finitary cocycle.

$(\Leftarrow)$
Conversely, suppose that $\alpha$ is finitary.
Then  a sequence of cocycles $\alpha_j:H_3(\Bbb R)\times Y_j\to Z_j$
is well defined by \thetag{7-10}.
Since $\alpha(c(t),y)=\varphi(t)$ for all $t\in\Bbb R$ and $y\in Y$, it follows from \thetag{7-10}  that $\alpha_j$ satisfies \thetag{7-8} for each $j\in\Bbb N$.
This yields that the $\alpha_j$-skew product extension $H_3(\Bbb R)$ (see \thetag{7-7})
is transitive on $Y_j\times Z_j$.
Then by Proposition~5.12, $T$ is a Heisenberg odometer.
\qed
  \enddemo

\proclaim{Corollary 7.8} Under the notation of Theorem~7.7, if $S$ is closed in $\Bbb R^2$ then
$T$ is an $H_3(\Bbb R)$-odometer.
This means that the subclass of degenerate Heisenberg odometers is spectrally determined.
\endproclaim

We now consider  the problem of spectral determinacy  within  the entire class Heisenberg odometers.

\proclaim{Corollary 7.9}
Let $T$ be a non-degenerate Heisenberg odometer associated with a sequence $\Gamma=(\Gamma_j)_{j=1}^\infty$
of lattices in $H_3(\Bbb R)$.
Then there is an ergodic action
$R$ of $H_3(\Bbb R)$ such that
\roster
\item"---" $R$ has the same maximal spectral type as $T$ but
\item"---" $R$ is not isomorphic to $T$ (and hence to any $H_3(\Bbb R)$-odometer).
\endroster
\endproclaim

\demo{Proof}
By Lemma~7.3 and Theorem~7.7, we can represent $T$ as $D^\beta$ for a $Z$-valued finitary cocycle $\beta$ of $D$ satisfying  \thetag{7-6}.
Let $\alpha:\Bbb R^2\times Y\to Z$ be a cocycle of the underlying $\Bbb R^2$-odometer.
By analogy with Definition~7.5, we call $\alpha$ {\it finitary} if
$$
\overline r_j(\alpha(s,y))=\alpha_j(s,\overline q_j(y)), \qquad s\in\Bbb R^2,j\in\Bbb N
$$
for a sequence of cocycles $\alpha_j:\Bbb R^2\times Y_j\to Z_j$.
Since $T$ is non-degenerated, the underlying $\Bbb R^2$-odometer is not transitive.
Therefore there is a cocycle  $\alpha$ of it  which is not cohomologous to a finitary one.
Then we  consider a  $Z$-valued cocycle
$$
\beta'(g,y):=\beta(g,y)+\alpha(p(g),y)
$$ 
of $D$.
Of course, $\beta'$ satisfies \thetag{7-6}.
Therefore by Lemma~7.4(i), the maximal spectral type of the $\beta'$-skew product extension $D^{\beta'}$ of $D$ coincides with the maximal spectral type of $T$.
However $\beta'$ is not finitary.\footnote{This is because $\beta$ is finitary and $\alpha$ is not finitary.}
Therefore by Theorem~7.7,
 $D^{\beta'}$ is not isomorphic to an odometer.
\qed
\enddemo

To construct an explicit counter-example to the spectral determinacy of non-degenerate Heisenberg odometers we first restate  Theorem~7.7 in a more suitable (for this purpose) form.
Let $(\Gamma_j)_{j=1}^\infty$, $D$, $Z$ and $\alpha$ be the same objects as in Lemma~7.3 and Theorem~7.7.
We need some more notation.
Let $D'$ denote the $p(\Gamma_1)$-odometer associated with the  sequence $p(\Gamma_1)\supset p(\Gamma_2)\supset\cdots$ of cofinite subgroups in $p(\Gamma_1)$.
Let $Y':=\projlim_{j\to\infty}p(\Gamma_1)/p(\Gamma_{j})$ be the space of $D'$.
We note that $Y'$ is a compact Abelian group (furnished with Haar measure).
Let $\widetilde\Gamma_1:=p^{-1}(p(\Gamma_1))$.
Denote by $\widetilde D=(\widetilde D_\gamma)_{\gamma\in\widetilde\Gamma_1}$ the following $\widetilde\Gamma_1$-action on $Y'$:
$$
\widetilde D_\gamma y':=p(\gamma)\cdot y'.
$$
Then one can deduce easily from Corollary~2.5 that $D=\text{Ind}_{\widetilde\Gamma_1}^{H_3(\Bbb R)}(\widetilde D)$.
Therefore by Lemma~2.4(iii), there is a $Z$-valued cocycle $\beta$ of $\widetilde D$ such that $\alpha$ is cohomologous to the induced cocycle $\widetilde\beta$ of $D$.
Moreover, by Lemma~2.4(ii),
the cohomology class of $\beta$ is defined uniquely by a cohomology class of $\alpha$.
Let $\goth Y_j$ denote the factor of $D$ generated by the canonical projection $\overline{p}_j:Y\to Y_j$ and let  $\widetilde{\goth Y}_j$ denote the factor of $\widetilde D$ generated by the canonical projection $\overline{p}'_j:Y'\to p(\Gamma_1)/p(\Gamma_j)$.
It is straightforward to verify that $\goth Y_j$ is induced by $\widetilde{\goth Y}_j$.
Therefore Lemma~2.4(ii), (iv) yield that $\alpha$ is cohomologous to a finitary cocycle if and only if $\beta$ is cohomologous to a finitary cocycle, say $\beta'$, of $\widetilde D$, i.e. the map
$$
Y'\ni y'\mapsto\overline r_j(\beta'(\gamma,y'))\in Z_j
$$
is $\widetilde{\goth Y}_j$-measurable for each $\gamma\in\widetilde\Gamma_1$  and $j\in\Bbb N$.
Thus we obtain  the following corollary from Theorem~7.7.

\proclaim{Corollary 7.10} Using the notation of Theorem~7.7, $T$ is an $H_3(\Bbb R)$-odometer (associated with $\Gamma$) if and only if $\beta$ is cohomologous   to a finitary cocycle.
\endproclaim

We are now ready to construct an action of $H_3(\Bbb R)$ which is unitarily equivalent to a Heisenberg odometer but not isomorphic to it.

\example{Example 7.11}
Let $p_1,p_2,\dots$
be an infinite sequence of pairwise different primes, $p_1=1$.
We set
$$
\Gamma_j:=\{c(n)b(p_1\cdots p_jm)a(p_1\cdots p_jl)\mid n,m,l\in\Bbb Z\}.
$$
Then $\Gamma:=(\Gamma_j)_{j=1}^\infty$ is a decreasing sequence of lattices in $H_3(\Bbb Z)$.
Denote by $T$ the associated $H_3(\Bbb R)$-odometer.
As in Lemma~7.3, represent   $T$ as a skew product extension $T=D^\alpha$ for a cocycle $\alpha$ of $D$ with values in $Z$.
We note that in this case $Z=Z_j=\Bbb R/\Bbb Z$ for each $j\in\Bbb N$.
By Theorem~7.7, $\alpha$ is (up to a cohomological change) finitary.
It is straightforward to verify that the following holds
\roster
\item"(i)"
$Y'=\projlim_{j\to\infty}\Bbb Z^2/p_1\cdots p_j\Bbb Z^2$,
\item"(ii)" $\widetilde\Gamma_1=\{c(t)b(m)a(l)\mid t\in\Bbb R, m,l\in\Bbb Z\}$ and
\item"(iii)"
$\alpha$ is cohomologous to the cocycle induced from the cocycle $\beta$ of $\widetilde D$ given by
$\beta(\gamma,y')=t+\Bbb Z$  if $\gamma=c(t)b(m)a(n)\in\widetilde\Gamma_1$.
\endroster
Thus $\beta$ is finitary.
It follows from (i) and the mutual coprimeness  of $p_j$, $j\in \Bbb N$, that the dynamical system $(Y',D')$ is isomorphic to the infinite product $Y'=\bigotimes_{j\in\Bbb N}\Bbb Z^2/p_j\Bbb Z^2$ such that
$$
D'_\gamma(y_1,y_2,\dots)=(\gamma_1+y_1,\gamma_2+y_2,\dots)
$$
for all $\gamma\in p(\Gamma_1)=\Bbb Z^2$ and $y_j\in \Bbb Z^2/p_j\Bbb Z^2$, where $\gamma_j:=\gamma+p_j\Bbb Z^2$, $j\in\Bbb N$.
Partition $\Bbb N$ into two infinite disjoint subsets $J_1$ and $J_2$ and
let $Y'_1:=\bigotimes_{j\in J_1}\Bbb Z^2/p_j\Bbb Z^2$ and
$Y'_2:=\bigotimes_{j\in J_2}\Bbb Z^2/p_j\Bbb Z^2$.
Then $Y_1'$ and $Y_2'$ are natural (coordinate) factors of $D'$.
Denote the restriction of $D'$ to these factors by $D^1$ and $D^2$ respectively.
We see that $D^1$ and $D^2$ are both free and $D'$ is naturally isomorphic to the Cartesian product $D^1\times D^2$.

We now note that the group  $Y'$ acting on itself by rotations is the centralizer $C(D')$ of $D'$.
Since $D'$ has pure point spectrum, $\widetilde{ \goth Y}_j$ is the $\sigma$-algebra of subsets fixed by a compact subgroup of $Y'$.
We can find a family of elements $(S_j)_{j\in\Bbb N}$  in $C(D')$ such that the following two properties are satisfied:
\roster
\item"(iv)" $\widetilde {\goth Y}_j$ is  fixed by $S_j$, i.e. $S_j$ is trivial on $\widetilde {\goth Y}_j$.
\item"(v)"
If we represent $S_j$ as the product $S_j=S_{1,j}\times S_{2,j}$ with $S_{1,j}\in C(D^1)$ and $S_{2,j}\in C(D^2)$
then the
`joint' action
$$
\Bbb Z^2\times\bigoplus_{j=1}^\infty\Bbb Z\ni(\gamma,n_1,n_2,\dots)\mapsto D^1_\gamma S_{1,1}^{n_1}S_{1,2}^{n_2}\cdots
$$
of the group $\Bbb Z^2\times\bigoplus_{j=1}^\infty\Bbb Z$ on $Y'_1$ is free.
\endroster
By an auxiliary Lemma~7.12 (see the end of this section), there is a cocycle $\omega$ of $D^1$ with values in $Z$ such that the cocycle
$\overline{\omega}:=\bigotimes_{j\in J_1} \omega\circ S_{1,j}$  with values in $Z^{J_1}$ is ergodic and, moreover, the product
$\Bbb Z^2$-action $(D^1)^{\overline{\omega}}\times D^2$  on $Y_1'\times Z^{J_1}\times Y_2$ is ergodic.
In other words, the $Z^{J_1}$-valued cocycle $\overline{\omega}\otimes 1$ of $D'$ (recall that $D'=D^1\times D^2$) is ergodic.

{\bf Claim A.} For any $j$, the cocycle $\omega\otimes 1$ of $D'$ is not cohomologous to a $\goth Y_j$-measurable cocycle.
Indeed, if $\omega\otimes 1$  is cohomologous to a $\goth Y_j$-measurable cocycle then in view of (iv),
$\omega\otimes 1$ is cohomologous to the cocycle $(\omega\otimes 1)\circ S_j$.
Therefore the $Z^2$-valued cocycle $\omega\otimes 1\times(\omega\otimes 1)\circ S_j$ of $D'$ is not ergodic.
However it is easy to see that
$$
\omega\otimes 1\times(\omega\otimes 1)\circ S_j=(\omega\times\omega\circ S_{1,j})\otimes 1
$$
and the skew product  $(D')^{(\omega\times\omega\circ S_{1,j})\otimes 1}$ is a factor
of the skew product  $(D')^{\overline{\omega}\otimes 1}$ which is erdodic by the choice of $\omega$.
We get a contradiction which proves Claim~A.

We now define a $Z$-valued cocycle $\kappa$ of $\widetilde D$ by setting
$$
\kappa(\gamma,y'):=\beta(\gamma,y')+(\omega\otimes 1)(p(\gamma), y')
$$
for  each $\gamma\in\widetilde\Gamma_1$ and $y'\in Y'$.
Denote by $\widetilde\alpha$ the $Z$-valued cocycle of $D$ induced from $\kappa$.
Finally, we let $\widetilde T:=D^{\widetilde\alpha}$.
Then $\widetilde T$ is an ergodic $H_3(\Bbb R)$-action.

It follows from Claim~A  that $\kappa$ is not cohomologous to a finitary cocycle.
Now Corollary~7.10 yields that $\widetilde T$ is not an $H_3(\Bbb R)$-odometer.
On the other hand, by Lemma~7.4(i), the maximal spectral type of the Koopman unitary representation $U_{\widetilde T}$
coincides with the maximal spectral type of $U_T$.
Since $\xi_\Gamma=\Bbb Z$ and $\pi_n$ is contained in $U_T$ with infinite multiplicity for each $n\ne 0$ by Theorem~5.2(i), to prove that $T$ and $\widetilde T$ are unitarily equivalent
it remains to show that the multiplicity $\pi_n$ in $U_{\widetilde T}$ is also infinite for each $n\ne 0$.

Denote by $U$ the Koopman unitary representation of $\widetilde\Gamma_1$ generated by the skew product $ \widetilde D^\kappa$.
The space of $U$ is $L^2(Y'\times Z)$.
It decomposes into an orthogonal sum
$$
L^2(Y'\times Z)=\bigoplus_{\eta\in\widehat { Y_2'}}\bigoplus_{n\in\Bbb Z} L^2(Y_1')\otimes\eta\otimes\chi_n
$$
of   $U$-invariant subspaces $ L^2(Y_1')\otimes\eta\otimes\chi_n$.\footnote{The subspace $ L^2(Y_1')\otimes\eta\otimes\chi_n$ is invariant under $U$ because $\omega\otimes 1$ as a function of $(\gamma,y')=(\gamma,y'_1,y'_2)\in\widetilde \Gamma_1\times Y_1'\times Y_2'$ does not depend on $y_2'$ while $\beta$ does not depend on $y'$ at all.}
Since the Koopman unitary representation $U_{\widetilde T}$ is induced from $U$, we obtain the following decomposition of $L^2(Y\times Z)$ into $U_{\widetilde T}$-invariant subspaces:
$$
L^2(Y\times Z)=\bigoplus_{{\eta\in\widehat { Y_2'}} }\bigoplus_{n\in\Bbb Z}L^2(\Bbb R^2/\Bbb Z^2)\otimes L^2(Y_1')\otimes\eta\otimes\chi_n.
$$
Moreover, it is straightforward to verify that $U_{\widetilde T}(c(t))v=e^{2\pi int}v$ for each vector
$v\in L^2(\Bbb R^2/\Bbb Z^2)\otimes L^2(Y_1')\otimes\eta\otimes\chi_n.$
Therefore the restriction of $U_{\widetilde T}$  to the subspace $L^2(\Bbb R^2/\Bbb Z^2)\otimes L^2(Y_1')\otimes\eta\otimes\chi_n$ is a multiple of $\pi_n$ for each $n\ne 0$ and
$\eta\in\widehat { Y_2'}$.
Since the group $\widehat { Y_2'}$ is infinite, the multiplicity of  $\pi_n$ in $U_{\widetilde T}$ is also infinite, as desired.
\endexample

We now state and prove an auxiliary lemma used in Example 7.11.

\proclaim{Lemma 7.12}
Let $K$ be a compact metric group.
Let $R=(R_g)_{g\in G}$ and $Q=(Q_g)_{g\in G}$ be two ergodic actions of an infinite countable discrete Abelian group $G$  such that $R$ is free and the Cartesian product $R\times Q$ is ergodic.
Fix a countable infinite Abelian subgroup $A$ in the centralizer $C(R)$ of $R$ such that $A\cap\{R_g\mid g\in G\}=\{\text{\rom{Id}}\}$.
Then there is a cocycle $\omega$ of $R$ with values in $K$ such that the $K^A$-valued cocycle
$\overline\omega:=\bigotimes_{a\in A} \omega\circ a$ of $R$ is ergodic and the Cartesian product $R^{\overline\omega}\times Q$ is ergodic.
\endproclaim
\demo{Proof}
To prove this lemma we will use the orbit theory of amenable  actions and their cocycles (see \cite{Sc}, \cite{BeG}, \cite{GoSi}, \cite{Da1}).
Let $(X,\mu)$ denote the space of $R$.
We denote by $\Cal R\subset X\times X$ the $R$-orbit equivalence relation.
Recall that the {\it full group} $[\Cal R]$ of $\Cal R$ is the group of all $\mu$-preserving invertible transformations $F$ of $X$ such that  $Fx$ belongs to the $\Cal R$-class of $x$ for a.a. $x\in X$.
Choose
\roster
\item "(i)"
a weakly mixing free action $V=(V_g)_{g\in G}$ of $G$ such that  the orbit equivalence relation of $V$ is $\Cal R$ and
\item"(ii)"
a cocycle $\tau$ of $V$ with values in $K$ such that the skew product $V^\tau$ is weakly mixing.
\endroster
Let $V^A=(V^A(g))_{g\in G}$ denote the product $G$-action on the  product space $(X^A,\mu^A)$, i.e.
$(V^A_gx_b):=V_gx_b$ for each $x=(x_b)_{b\in A}\in X^A$ and $g\in G$.
We consider also a shift  action $S=(S_a)_{a\in A}$ of $A$ on $X^A$:
$$
(S_ax)_{b}:=x_{b-a},\qquad a,b\in A,\quad x=(x_b)_{b\in A}\in X^A.
$$
Of course,
\roster
\item"(iii)" $S_a\in C(V^A)$  and $S_a\not\in \{V^A_b\mid b\in A\}$ for each $ a\in A\setminus\{0\}$.
\endroster
We define a cocycle $\tau^A$ of $V^A$ with values in $K$  by setting
$\tau^A(g,x):=\tau(g,x_0)$ for all $g\in G$ and $x\in X^A$.
Consider the following cocycle $\overline{\tau^A}$ of $V^A$ with values in $K^A$:
$$
\overline{\tau^A}=\bigotimes_{a\in A}{\tau^A}\circ S_a^{-1}.
$$
We claim that $\overline{\tau^A}$ is ergodic.
For that we will show that the $\overline{\tau^A}$-skew product extension of $V^A$ is ergodic.
Take $g\in G$, $x=(x_b)_{b\in A}\in X^A$ and $y=(y_b)_{b\in A}\in K^A$.
Then $(x,y)=(x_a,y_a)_{a\in A}\in (X\times K)^A$ and
$$
\multline
((V^A)^{\overline{\tau^A}}_g(x,y))_a=(V^A_gx, \overline{\tau^A}(g,x)+y)_a
=(V_gx_a,\tau^A(g,S_a^{-1} x)+y_a ) \\=(V_gx_a,\tau(g, x_a)+y_a)=V^\tau_g(x_a,y_a).
\endmultline
$$
 Thus, $(V^A)^{\overline{\tau^A}}=(V^\tau)^A$, where $(V^\tau)^A:=\bigotimes_{a\in A}V^\tau$.
 Since the skew product $V^\tau$ is weakly mixing, $(V^A)^{\overline{\tau^A}}$ is ergodic, i.e.
 ${\overline{\tau^A}}$ is ergodic, as claimed.

 Denote by $\Cal V$ the $V^A$-orbit equivalence relation on $X^A$.
 It is standard to verify that (iii) implies that $\{S_a\mid a\in A\}\cap[\Cal V]=\{\text{Id}\}$.
 In a similar way, we deduce from the condition of the lemma that $ A\cap[\Cal R]=
 \{\text{Id}\}$.
 Therefore we may apply the main result from \cite{BeG}: there is a measure preserving isomorphism $L:X\to X^A$ and a map $A\ni a\mapsto\gamma_a \in  [\Cal V]$ such that
 $(L\times L)(\Cal R)=\Cal V$
  and $LaL^{-1}=\gamma_aS_a^{-1}$ for all $a\in A$.
  
  We note that $\tau^A$ extends naturally to $[\Cal V]$.
  Indeed, for each $\gamma\in[\Cal V]$, there is a unique (mod 0) partition $X^A=\bigsqcup_{g\in G}Z_g$ of $X^A$
into Borel subsets $Z_g$, $g\in G$, such that $\gamma z=(V^A)_gz$ whenever $z\in Z_g$.
We now set $\tau^A(\gamma,z):=\tau^A(g,z)$ whenever $z\in Z_g$.
Since $V^A$ is free, $\tau^A(\gamma,.)$ is well defined.
It is easy to verify that 
$$
\tau^A(\gamma_1\gamma_2,z)=\tau^A(\gamma_1,\gamma_2z)\tau^A(\gamma_2,z)
$$
at a.e. $z$ for all $\gamma_1,\gamma_2\in[\Cal V]$.
 We now define a cocycle $\omega$ of $R$ with values in $K$ by setting
  $$
  \omega(g,x):=\tau^A(LR_gL^{-1}, Lx)\qquad\text{for all $g\in G$ and $x\in X$.}
  $$
 For each $a\in A$, we have
 $$
 \align
 \omega\circ a(g,x) & =\omega(g,ax)\\
 &=\tau^A(LR_gL^{-1},Lax)\\
 &=\tau^A(LR_gL^{-1},\gamma_aS_a^{-1}Lx)\\
 &=f(R_gx)+\tau^A(S_a^{-1}LR_gL^{-1}S_a, S_a^{-1}Lx)-f(x),
 \endalign
 $$
 where $f(x):=\tau^A(\gamma_a,S_a^{-1}Lx)$.
 Thus,  $\omega\circ a$  is cohomologous to  $\tau^A\circ (S_a^{-1}L)=(\tau^A\circ S_a^{-1})\circ L$.
 This implies that the cocycle $ \overline\omega$ of $R$, given by
  $$
   \overline\omega(g,x):=\bigotimes_{a\in A} \omega(g, ax)
   $$ 
   for all $x\in X$, $g\in G$,
 is cohomologous to the cocycle 
 $$
 \bigg(\bigotimes_{a\in A} \tau^A \circ S_a^{-1}\bigg)\circ L
 =\overline{\tau^A}\circ L.
 $$
 Since $\overline{\tau^A}$ is ergodic,  $\overline{\omega}$ is ergodic too.
 Thus the first claim of the lemma is proved.

 We now show how to deduce the second claim of the lemma   from the first one.
 Given an ergodic action $W$ of $G$, we denote by $\Lambda_d(W)$ the discrete spectrum of $W$.
 Then $\Lambda_d(W)$ is a countable subgroup of $\widehat G$.
 It is  know that $R^{\overline\omega}\times Q$ is ergodic if and only if
 $\Lambda_d(R^{\overline\omega})\cap\Lambda_d(Q)=\{1\}$.
 Of course, $\Lambda_d(R)\subset \Lambda_d(R^{\overline\omega} )$.
 Since $R\times Q$ is ergodic, $\Lambda_d(R)\cap\Lambda_d(Q)=\{1\}$.
 Therefore to prove that $R^{\overline\omega}\times Q$ is ergodic
 it suffices to show that $\Lambda_d(R)=\Lambda_d(R^{\overline\omega})$.
Suppose that there is $\kappa\in\Lambda_d(R^{\overline\omega})\setminus\Lambda_d(R)$.
 Then for some character $\chi\in \widehat {K^A}$, the $\Bbb T$-valued cocycle $\chi\circ\overline\omega$ of $R$ is cohomologous to the cocycle $G\times X\ni(g,x)\mapsto \kappa(g)\in\Bbb T$.
 Since $\widehat {K^A}=\bigoplus_{a\in A}\widehat K$,
 there is a finite non-empty subset $A_0\subset A$ and a map $A_0\ni a\mapsto\chi_a\in\widehat K\setminus\{1\}$ such that
 $\chi(k)=\prod_{a\in A_0}\chi_a(k_a)$ for each $k=(k_a)_{a\in A}\in K^A$.
 Since $A$ is infinite,
 there is $b\in A$ such that $A_0\cap(b+A_0)=\emptyset$.
 For each subset $B\subset A$, we define by  $\overline\omega_B$ the cocycle
 $ \bigotimes_{a\in B}\omega\circ a$ of $R$ with values in $K^B$.
 Since $\overline\omega$ is ergodic, $\overline\omega_{B}$ is also ergodic.
 We note that
 $$
 \overline\omega_{A_0\cup (b+A_0)}=\overline\omega_{A_0}\times\overline\omega_{b+A_0}=
 \overline\omega_{A_0}\times\overline\omega_{A_0}\circ S_b.
 $$
Therefore the $\Bbb T^2$-valued cocycle 
$$
\chi\circ\overline\omega_{A_0}\times\chi\circ\overline\omega_{A_0}\circ S_b=\chi\circ\overline\omega\times\chi\circ\overline\omega\circ S_b
$$ 
of $R$ is also ergodic.
However the later cocycle is cohomologous to the  cocycle 
$$
G\times X\ni(g,x)\mapsto (\kappa(g),\kappa(g))
$$
  taking values in the diagonal of the torus $\Bbb T^2$.
  Hence this cocycle is not ergodic, a contradiction.
  \qed
\enddemo

\head 8. On $H_3(\Bbb Z)$-odometers
\endhead

Let $\Gamma_1\supset\Gamma_2\supset\cdots$ be a decreasing sequence of lattices (i.e. cofinite subgroups) in $H_3(\Bbb Z)$.
Denote by $T=(T_g)_{g\in H_3(\Bbb Z)}$ the associated $H_3(\Bbb Z)$-odometer.
Let $(X,\mu)$ be the space of this odometer.
We call $T$ {\it normal} if $\Gamma_j$ is normal in $H_3(\Bbb Z)$ for each $j$.\footnote{There is a difference in nomenclature used in our paper and \cite{Li-Ug}. 
By an $H_3(\Bbb Z)$-odometer the authors of \cite{Li-Ug} mean a free normal $H_3(\Bbb Z)$-odometer.
An $H_3(\Bbb Z)$-odometer satisfying the condition $\bigcap_{j=1}^\infty\Gamma_n=\{1\}$ is called a $H_3(\Bbb Z)$-{\it subodometer} there. 
The $H_3(\Bbb Z)$-odometers which do not satisfy this condition were not considered there.
We prefer to use the term ``normal'' to make it consistent with the well known {\it normality} concept  
  introduced by R.~Zimmer in \cite{Zi} (see also \cite{Fa}).}
  
If $T$ is normal then $X$ is a compact totally disconnected group and $\mu$ is the normalized Haar measure on $X$.
Indeed, we obtain a sequence 
$$
H_3(\Bbb Z)/\Gamma_1\leftarrow H_3(\Bbb Z)/\Gamma_2\leftarrow\cdots
$$ 
of finite groups $H_3(\Bbb Z)/\Gamma_j$ and canonical onto homomorphisms such that $X=\projlim_{j\to \infty}H_3(\Bbb Z)/\Gamma_j$.
Moreover, a group homomorphism $\varphi:H_3(\Bbb Z)\to X$ is well defined by
the formula $\varphi(g)=(\varphi(g)_j)_{j=1}^\infty$, where $\varphi(g)_j:=g\Gamma_j$.
Of course, $\varphi(H_3(\Bbb Z))$ is dense in $X$.
It is easy to see that $T_gx=\varphi(g)x$ for all $g\in H_3(\Bbb Z)$ and $x\in X$.
Hence $T$ has a pure point spectrum in the sense of \cite{Ma3}.
Moreover, $T$ is normal in the sense of \cite{Zi}.
Hence  \cite{Fa, Corollary~2} yields the following.

\proclaim{Corollary 8.1}
The normal $H_3(\Bbb Z)$-odometers  are isospectral.
\endproclaim

We also note that every ergodic 2-fold self-joining of a normal $H_3(\Bbb Z)$-odometer $T$ is off-diagonal.
Hence $T$ is 2-fold simple \cite{dJRu}.

Let $L_j$ denote the left regular representation of $H_3(\Bbb R)/\Gamma_j$.
Let $\Cal I_j$ stand for the unitary dual of $H_3(\Bbb R)/\Gamma_j$.
\footnote{Since we do not distinguish between unitarily equivalent representation, we consider
 elements of  $\Cal I_j$ as   irreducible unitary representations of 
$H_3(\Bbb R)/\Gamma_j$ rather then their unitary equivalence classes.}
It is well known that (up to the unitary equivalence)
$L_j=\bigoplus_{\tau\in I_j}\bigoplus_1^{d_\tau}\tau$, where $d_\tau$ is the dimension of $\tau$.
In particular, $\#(H_3(\Bbb R)/\Gamma_j)=\sum_{\tau\in\Cal I_j}d_\tau^2$.
Moreover, $\#\Cal I_j$ equals the cardinality of the set of congugacy classes in $H_3(\Bbb R)/\Gamma_j$.
Denote by $U_T$  the  Koopman unitary representation of $H_3(\Bbb Z)$ generated by $T$.
The canonical projection $X\to H_3(\Bbb R)/\Gamma_j$ generates an embedding $L^2(H_3(\Bbb R)/\Gamma_j)\subset X$.
Therefore we obtain an increasing sequence
$$
L^2(H_3(\Bbb R)/\Gamma_1)\subset L^2(H_3(\Bbb R)/\Gamma_2)\subset\cdots
$$
of $U_T$-invariant subspaces whose union is dense in $L^2(X)$ and such that the restriction    $U_T\restriction
L^2(H_3(\Bbb R)/\Gamma_j)$ is unitarily equivalent to $L_j\circ p_j$,
where $p_j:H_3(\Bbb R)\to H_3(\Bbb R)/\Gamma_j$ is the canonical projection.
This implies that 
$$
\{\tau\circ p_1\mid \tau\in \Cal I_1\}\subset\{\tau\circ p_2\mid \tau\in \Cal I_2\}\subset\cdots
$$
and we obtain the following decomposition of $U_T$ into the orthogonal sum of irreducible finite dimensional unitary representations of $H_3(\Bbb Z)$.

\proclaim{Theorem 8.2}
Let $\Cal I_T:=\bigcup_{j\in\Bbb N}\{\tau\circ p_j\mid \tau\in \Cal I_j\}$ and $d_\iota$ is the dimension of $\iota$.
Then we have
$$
U_T=\bigoplus_{\iota\in\Cal I_T}\bigoplus_{1}^{d_\iota}\iota.
$$
\endproclaim

An explicit  computation of $\Cal I_T$ in terms of the sequence $(\Gamma_j)_{j=1}^\infty$
was done in \cite{Li-Ug}.

\proclaim{Corollary 8.3}
Two normal $H_3(\Bbb Z)$-odometers $T$ and $R$ are (measure theoretically) isomorphic  if and only if $\Cal I_T=\Cal I_R$.
\endproclaim

We now provide an example of  non-isomorphic normal  $H_3(\Bbb Z)$-odometers $V$ and $V'$ such that the Koopman representations of $H_3(\Bbb R)$ generated by the $H_3(\Bbb R)$-odometers  associated with the same sequences of lattices (as $V$ and $V'$ respectively) are unitarily equivalent.\footnote{We note that the Koopman representations of $H_3(\Bbb Z)$ generated by $V$ and $V'$ are not unitarily equivalent by Corollary 8.1. Thus we obtain an example of unitarily non-equivalent Koopman unitary representations $U_V$ and $U_{V'}$ of $H_3(\Bbb Z)$ with pure point spectrum such that the induced unitary representations Ind$_{H_3(\Bbb Z)}^{H_3(\Bbb R)}(V)$
and Ind$_{H_3(\Bbb Z)}^{H_3(\Bbb R)}(V')$ of $H_3(\Bbb R)$ are unitarily equivalent.}
We proceed it with a lemma which was utilized in \cite{Li--Ug} with reference to \cite{CoPe}.
We provide an alternative proof.

\proclaim{Lemma 8.4}
Let $V$ be the  $H_3(\Bbb Z)$-odometer associated with a sequence $(\Gamma_j)_{j=1}^\infty$ of normal lattices in $H_3(\Bbb Z)$ such that $\bigcap_{j=1}^\infty\Gamma_j=\{1\}$.
Let $\Sigma$ be  a lattice in $H_3(\Bbb Z)$ such that the homogeneous space $H_3(\Bbb Z)/\Sigma$ is  a factor, say $\goth F$, of $V$.
Then there is $j_0>0$ such that $\Sigma\supset \Gamma_{j_0}$.
\endproclaim

\demo{Proof}
Since $V$ is 2-fold simple, we apply the Veech theorem \cite{dJRu}:
there is a compact subgroup $K$ in the centralizer $C(V)$ of $V$ such that $\goth F$
is the $\sigma$-algebra  Fix\,$K$ of subsets fixed by $K$.
Since $\goth F$ is finite, $K$ is open.
We recall that the space of $V$ is a compact group $X=\projlim_{j\to\infty}H_3(\Bbb Z)/\Gamma_j$.
Hence the centralizer $C(V)$ is isomorphic to $X$ acting on itself by right translations.
\footnote{We recall that $V$ acts on $X$ by left translations.}
Let $\Lambda_j$ be the kernel of the canonical projection $q_j:X\to H_3(\Bbb Z)/\Gamma_j$.
Then $\Lambda_j$ is on open subgroup in $X$, $\Lambda_1\supset\Lambda_2\supset\cdots$
and $\bigcap_{j=1}^\infty\Lambda_j=\{1\}$.
Therefore there is $j_0>0$ such that $K\supset\Lambda_{j_0}$.
Hence $\goth F\subset\text{Fix}\,\Lambda_{j_0}$.
However $\text{Fix}\,\Lambda_{j_0}$ is exactly the factor of $V$  determined by $q_{j_0}$.
Thus we obtain that $H_3(\Bbb Z)/\Sigma$ is a factor of $H_3(\Bbb Z)/\Gamma_{j_0}$.
Since $\Gamma_{j_0}$ is normal in $H_3(\Bbb Z)$, it follows that
 $\Sigma\supset\Gamma_{j_0}$.
 \qed
\enddemo

\example{Example 8.5 (\rom{cf. \cite{Li--Ug, Example 4.9}})}
Let $(\Gamma_n)_{n=1}^\infty$ and 
$(\Gamma_n')_{n=1}^\infty$
be  as in  Example~5.11.
Denote by $V$ and $V'$ the $H_3(\Bbb Z)$-odometers associated with the sequences
$(\Gamma_n)_{n=1}^\infty$
and $(\Gamma_n')_{n=1}^\infty$  respectively.
They are normal.
It was shown in \cite{Li--Ug, Example 4.9} that $T$ and $T'$ are not isomorphic.
\footnote{This fact follows also from Example~5.11 because if $V$ and $V'$ are isomorphic then the induced $H_3(\Bbb R)$-actions are also isomorphic. However these induced actions are the $H_3(\Bbb R)$-odometers associated with $(\Gamma_n)_{n=1}^\infty$
and $(\Gamma_n')_{n=1}^\infty$  respectively.
As was  shown in Example~5.11, these odometers are not isomorphic.}
Indeed, otherwise in view of Lemma~8.4, there is $n>0$ such that
$\Gamma'_n\subset\Gamma_2$.
However we see that $c(1)a(k_n)\in\Gamma_n'$ but $c(1)a(k_n)\not\in\Gamma_2$, a contradiction.
On the other hand, 
$b(1/k_n)\Gamma_n'b(-1/k_n)=\Gamma_n$ for each $n\in\Bbb N$.
Hence the $H_3(\Bbb R)$-odometers  $T$ and $T'$ associated with $(\Gamma_n)_{n=1}^\infty$
and $(\Gamma_n')_{n=1}^\infty$ are f-isomorpfic. 
By Corollary~5.8, the Koopman unitary representations of $H_3(\Bbb R)$ generated by them are unitarily equivalent.
We also examine the ``symmetry'' property   for $V$ and $V'$.
Since $\theta(H_3(\Bbb Z))=H_3(\Bbb Z)$, the symmetric $H_3(\Bbb Z)$-actions are defined
in a similar way as for the $H_3(\Bbb R)$-actions.
It is easy to see that $V$ is symmetric.
It is straightforward to verify that
$$
\theta(\Gamma_n')=a(-1/k_n)\Gamma_na(1/k_n)=
\{c(k_nj_3+j_2)b(k_nj_2)a(k_nj_1)\mid j_1,j_2,j_3\in\Bbb Z\}
$$ 
for each $n$.
If $V'$ and $V'\circ\theta$ were isomorphic then by Lemma~8.4, there is $n>0$ such that
$\theta(\Gamma_n')\subset\Gamma_2'$.
However we see that $c(k_n+1)b(k_n)\in\theta(\Gamma_n')$ but
$c(k_n+1)b(k_n)\not\in\theta(\Gamma_2')$.
This contradiction yields that $V'$ is not symmetric.
  \endexample

\head 9. Concluding remarks and open problems
\endhead

In view of Corollary 8.1, we  ask the following natural question.

\remark{Problem 9.1} Whether the non-normal $H_3(\Bbb Z)$-odometers are isospectral?
\endremark

Let $G$ be a locally compact second countable group and let $\Gamma$ be a lattice in $G$.

\definition{Definition 9.2}
We call two  probability preserving actions $V$ and $V'$ of $\Gamma$  {\it flow equivalent} if
the induced $G$-actions $\text{Ind}_\Gamma^G(V)$ and $\text{Ind}_\Gamma^G(V')$ are isomorphic.
\enddefinition

In the case $G=\Bbb R$ and $\Gamma=\Bbb Z$ this definition corresponds to the classical concept of flow equivalence in topological dynamics.
As was observed, e.g. in \cite{DaLe, Proposition~1.3}, if $G$ is Abelian then two ergodic actions of $\Gamma$ are flow equivalent if and only if they are isomorphic.

\remark{Problem 9.3} Whether two flow equivalent  ergodic actions of $H_3(\Bbb Z)$  are isomorphic? 
Whether two flow equivalent   $H_3(\Bbb Z)$-odometers  are isomorphic? 
\endremark

In this connection we make the following remark.

\remark{Remark 9.4}
 Let $V$ and $V'$ be two ergodic probability preserving actions of $H_3(\Bbb Z)$.
Suppose that they are flow equivalent.
Let $(Z,\kappa)$ be the space of $V=(V_\gamma)_{\gamma\in H_3(\Bbb Z)}$ and let $(Z',\kappa')$ be the space of $V'=(V'_{\gamma})_{\gamma\in H_3(\Bbb Z)}$.
Denote the homogeneous space $H_3(\Bbb R)/H_3(\Bbb Z)$ by $Y$.
The translation of $y\in Y$ by $g\in H_3(\Bbb R)$ is denoted by $g* y$.
Let $\lambda$ stand for the Haar measure on $Y$.
Then
 $(Y\times Z,\lambda\times\kappa)$ is the space of Ind$_{H_3(\Bbb Z)}^{H_3(\Bbb R)}(V)$ and
 $(Y\times Z',\lambda\times\kappa')$ is the space of Ind$_{H_3(\Bbb Z)}^{H_3(\Bbb R)}(V')$.
 Thus $Y$ is a  factor of both $V$ and $V'$.
Let $R$ be a measure preserving isomorphism of $Y\times Z$ onto $Y\times Z'$ that conjugates Ind$_{H_3(\Bbb Z)}^{H_3(\Bbb R)}(V)$ with Ind$_{H_3(\Bbb Z)}^{H_3(\Bbb R)}(V')$.
Suppose, in addition, that
$$
\text{$R$ maps the factor $Y$ of $V$ to the factor $Y$ of $V'$.}\tag{9-1}
$$
Then
there
 exists a measurable field $Y\ni y\mapsto R_y$ of isomorphisms $R_y$ from    $(Z,\kappa)$ to $(Z',\kappa')$ and a transformation $Q$ of $Y$ commuting with the $H_3(\Bbb R)$-action on $Y$ such that
$R(y,z)=(Qy,R_yz)$ for a.a. $(y,z)\in Y\times Z$.
By \cite{Da2, Lemma~2.1}, $Q$ is a translation by an element of the normalizer of $H_3(\Bbb Z)$ in $H_3(\Bbb R)$.
A straightforward  verification shows that this normalizer is $\{c(t)b(m)a(n)\mid t\in\Bbb R,n,m\in\Bbb Z\}$.
Thus there is $t_0\in\Bbb R$ such that $Qy=c(t_0) * y$ for $\lambda$-a.a. $y\in Y$.
Since $R$ conjugates the induced $H_3(\Bbb R)$-actions, we obtain that
$$
(c(t_0)g*y,R_{g*y}V_{h_s(g,y)}z)= (gc(t_0)*y,V'_{h_s(g,c(t_0)*y)}R_{y}z)
$$
for  $(\lambda\times\kappa)$-a.a. $(y,z)\in Y\times Z$ for each $g\in H_3(\Bbb R)$.
Here $h_s:H_3(\Bbb R)\times Y\to H_3(\Bbb Z)$ is the choice cocycle  corresponding to a Borel cross-section $s:Y\to H_3(\Bbb R)$.
By Fubini's theorem,  there is $y_0\in Y$ such that
$$
R_{y_0}V_{h_s(g,y_0)}z= V'_{h_s(g,c(t_0)*y_0)}R_{y_0}z\tag9-2
$$
for $\kappa$-a.a. $z\in Z$ and all $g$ from the stability group $S_{y_0}:=\{g\in H_3(\Bbb R)\mid g*y_0=y_0\}$ of $y_0$.
We used here that $S_{y_0}$ is countable (in fact, it is conjugate to $H_3(\Bbb Z)$) to guarantee that the
``$\kappa$-a.a. $z\in Z$'' does not depend on $g\in S_{y_0}$.
Next, we have
$$
\align
h_s(g,c(t_0)*y_0) & =s(gc(t_0)*y_0)^{-1}gs(c(t_0)*y_0)\\
& =s(c(t_0)*y_0)^{-1}gs(c(t_0)*y_0)\\
& =r_0^{-1}s(y_0)^{-1}gs(y_0)r_0\\
&=r_0^{-1}h_s(g,y_0)r_0
\endalign
$$
for some $r_0\in H_3(\Bbb Z)$ and all $g\in S_{y_0}$.
We used the fact that 
$$
s(c(t_0)*y_0)\in c(t_0)s(y_0)H_3(\Bbb Z).
$$
Now \thetag{9-2} implies that
$$
R_{y_0}V_{h_s(g,y_0)}R_{y_0}^{-1}= V'_{r_0^{-1}h_s(g,y_0)r_0}=
(V'_{r_0})^{-1} V'_{h_s(g,y_0)} V'_{r_0}.
$$
Since  $\{h_s(g,y_0)\mid g\in S_{y_0}\}=H_3(\Bbb Z)$, it follows that $V$ and $V'$ are isomorphic.
\endremark

Thus Problem~9.3 reduces to the following question:

\remark{Problem 9.5} Suppose that
Ind$_{H_3(\Bbb Z)}^{H_3(\Bbb R)}(V)$ and Ind$_{H_3(\Bbb Z)}^{H_3(\Bbb R)}(V')$ are isomorphic.
Is there an isomorphism satisfying \thetag{9-1}?
\endremark

Let $V$ and $V'$ be the  $H_3(\Bbb Z)$-odometers associated with sequences $\Gamma_1\supset\Gamma_2\supset\cdots$ and
$\Gamma_1'\supset\Gamma_2'\supset\cdots$ of lattices in $H_3(\Bbb Z)$ respectively.
Let $T$ and $T'$  be the $H_3(\Bbb R)$-odometers associated with the same sequences of lattices (considered as lattices in $H_3(\Bbb R)$).
It follows from  Corollary~2.5 that  $V$ and $V'$ are flow equivalent if and only if $T$ and $T'$ are isomorphic.
The concept of f-isomorphism for the $H_3(\Bbb R)$-odometers motivated the following definition.

\definition{Definition 9.6}
Let $T$ and $T'$ be two ergodic actions of a locally compact second countable group $G$ on standard probability spaces $(X,\goth B,\mu)$ and $(X',\goth B',\mu')$ respectively. 
We say that $T$ and $T'$ are {\it f-isomorphic} if there are increasing sequences $\goth F_1\subset\goth F_2\subset\cdots$ and $\goth F_1'\subset\goth F_2'\subset\cdots$ of factors  of $T$ and $T'$ respectively such that $\bigvee_{j=1}^\infty\goth F_j=\goth B$,
$\bigvee_{j=1}^\infty\goth F_j'=\goth B'$ and $T\restriction\goth F_j$ is isomorphic to $T'\restriction\goth F'_j$ for each $j\in\Bbb N$.
We say that $T$ and $T'$ are {\it F-isomorphic} if there is a finite sequence $R_1,\dots,R_n$ of ergodic $G$-actions such that $R_1=T$, $R_n=T'$ and $R_j$ is f-isomorphic to $R_{j+1}$ for each $j=1,\dots, n-1$.
Of course,  f-isomorphism is a symmetric relation on the set of ergodic $G$-actions and  F-isomorphism is  the smallest equivalence relation  majorizing  f-isomorphism.
\enddefinition

We note that  F-isomorphism is weaker than  isomorphism.

\example{Example 9.7}
Let $T$ be a $G$-action with MSJ (see \cite{dJRu}).
Let $(X,\goth B,\mu)$ be the space of $T$.
Denote by $T\odot T$  the symmetric factor of the Cartesian product $T\times T$, i.e. the restriction of $T\times T$ to the $\sigma$-albegra of subsets invariant under the involution $(x,y)\mapsto(y,x)$.
Then it is easy to verify that the $G$-actions $R:=T\times T\times\cdots$ and
$R':=(T\odot T)\times T\times T\times\cdots$ are f-isomorphic.
Indeed, let $\goth F_n:=\goth B^{\otimes n}\otimes (\goth B\odot\goth B)\otimes\goth N\otimes\goth N\otimes\cdots$ 
and $\goth F'_n:=(\goth B\odot\goth B)\otimes \goth B^{\otimes n}\otimes\goth N\otimes\goth N\otimes\cdots$.
Then $(\goth F_n)_{n=1}^\infty$ is an increasing sequence of factors of $R$ with
$\bigvee_{n=1}^\infty\goth F_n$ being the entire $\sigma$-algebra of $R$;
$(\goth F_n')_{n=1}^\infty$ is an increasing sequence of factors of $R'$ with
$\bigvee_{n=1}^\infty\goth F_n'$ being the entire $\sigma$-algebra of $R'$; and $R\restriction\goth F_n$ is isomorphic to $R'\restriction\goth F_n'$ for each $n$.
On the other hand,  $R$ and $R'$ are not isomorphic (see \cite{dJRu}).
\endexample

\remark{Problem 9.8}
\roster
\item
Give an example of two ergodic $G$-actions which are F-isomorphic but not f-isomorphic.
\item
We note that in Example 9.7 the actions $R$ and $R'$ are weakly isomorphic.
Are there F-isomorphic $G$-actions which are not  weakly isomorphic?
 \item
Are there F-isomorphic  $G$-actions which are not Markov
quasi-equivalent? \footnote{See \cite{FPS} for the definitions.}
\endroster
\endremark

\enddocument